\documentclass[11pt,a4paper]{article}
\usepackage[T1]{fontenc}
\usepackage[utf8]{inputenc}
\usepackage{cite}
\usepackage[left=2cm,right=2cm,top=2cm,bottom=2cm]{geometry}
\usepackage{amsmath,amssymb,amsthm,amsfonts}
\usepackage{mathrsfs}
\usepackage{enumitem}
\usepackage[hidelinks]{hyperref}
\theoremstyle{plain}
\newtheorem{theorem}{Theorem}[section]
\newtheorem{proposition}[theorem]{Proposition}
\newtheorem{corollary}[theorem]{Corollary}
\newtheorem{lemma}[theorem]{Lemma}
\theoremstyle{definition}
\newtheorem{definition}[theorem]{Definition}
\newtheorem{problem}[theorem]{Problem}
\newtheorem{example}[theorem]{Example}
\theoremstyle{remark}
\newtheorem{remark}[theorem]{Remark}
\newcommand{\R}{\mathbb{R}}
\newcommand{\Z}{\mathbb{Z}}
\newcommand{\C}{\mathbb{C}}
\newcommand{\Hq}{\mathbb{H}}
\newcommand{\N}{\mathbb{N}}
\newcommand{\Cl}{\mathbb{C}_{m}}
\newcommand{\Cln}{\mathbb{C}_{2n}}
\newcommand{\Hqn}{\mathbb{H}_{4n}}
\newcommand{\CM}{\mathbb{CM}_{2\times 2}}
\newcommand{\HM}{\mathfrak{H}}
\newcommand{\dZ}{\partial_{Z}}
\newcommand{\dZd}{\partial_{Z^{\dagger}}}
\newcommand{\dX}{\partial_{X}}
\newcommand{\dXt}{\partial_{X^{|}}}
\newcommand{\Xt}{X^{|}}
\newcommand{\Yt}{Y^{|}}
\newcommand{\Zd}{Z^{\dagger}}
\newcommand{\Dop}{\mathcal{D}_{(Z,\Zd)}}
\newcommand{\DopK}{\mathcal{D}^{K}_{(Z,\Zd)}}
\newcommand{\Ekernel}{\boldsymbol{\mathcal{E}}}
\newcommand{\Gmat}{G^{1}_{2}}
\newcommand{\Fmat}{F^{1}_{2}}
\newcommand{\Xmat}{\Xi^{1}_{2}}
\newcommand{\Phimat}{\Phi^{1}_{2}}
\newcommand{\Jmat}{\mathfrak{J}^{1}_{2}}
\newcommand{\Hmat}{H^{1}_{2}}
\newcommand{\Zero}{\mathcal{O}}
\newcommand{\lam}{\lambda_{n}}
\newcommand{\Cauchy}{\mathcal{C}_{\Gamma}}
\newcommand{\Teo}{\mathcal{T}_{\Omega}}
\newcommand{\Hilb}{\mathcal{H}_{\Gamma}}
\newcommand{\circu}{\operatorname{circ}}
\newcommand{\dimH}{\dim_{H}}
\newcommand{\tr}{\operatorname{tr}}

\begin{document}

\title{Hermitian Clifford Analysis on Euclidean Domains with Smooth and Fractal
Boundaries: \\[2mm] {\large A Survey of Matrix Cauchy and Hilbert Transforms, Boundary Value Problems, Helmholtz-Type Dirac Operators and Duality.\par}}
\author{
Ricardo Abreu-Blaya\\[0.5mm]
{\small Facultad de Matem\'aticas, Universidad Aut\'onoma de Guerrero, Chilpancingo de los Bravo, M\'exico}\\{\small Email: rabreublaya@yahoo.es}\\[3mm]
Juan Bory-Reyes\thanks{Corresponding author}\\[0.5mm]
{\small SEPI-ESIME Zacatenco, Instituto Polit\'ecnico Nacional, Ciudad de M\'exico, M\'exico}\\{\small Email: juanboryreyes@yahoo.com}\\[3mm]
Carlos Daniel Tamayo-Castro\\[0.5mm]
{\small SEPI-ESIME Zacatenco, Instituto Polit\'ecnico Nacional, Ciudad de M\'exico, M\'exico}\\
{\small Facultad de Ciencias, Universidad Nacional Aut\'onoma de M\'exico, Ciudad de M\'exico, M\'exico}\\{\small Email: cdtamayoc@gmail.com}
}

\date{}
\maketitle

\begin{abstract}
Hermitian Clifford analysis is a refinement of Euclidean Clifford analysis obtained by imposing a complex structure on $\R^{2n}$ and considering the simultaneous null solutions (the \emph{$h$-monogenic} functions) of two mutually adjoint Hermitian Dirac operators.
This survey collects, in a unified notation and without proofs, the main results obtained along this line of research between 2008 and 2026.  We review the matrix Borel--Pompeiu, Koppelman and Cauchy formulae, first for domains with Liapunov or Ahlfors--David regular boundary and then for Jordan domains whose boundary is fractal in the sense of Mandelbrot; the associated Cauchy transform, Hilbert transform, Hardy projections and their mapping properties, including H\"older norm estimates in terms of the $d$-sum of a Whitney decomposition; the jump decomposition,
Dirichlet and Riemann boundary value problems in each setting, including the linear Riemann problems with variable circulant coefficients and their nonlinear counterparts. We also review a duality theory for the module of germs of $H$-monogenic matrix functions, together with its Euclidean model (the annihilators of harmonic differential forms) and its algebraic prerequisites (a Poincar\'e lemma and a Mittag--Leffler theorem for circulant matrices); the function theory associated with a Hermitian version of the Helmholtz equation, the corresponding $\bar\partial$-problem and its application to a formal two-dimensional Maxwell-type system; and finally, the parallel development in the quaternionic Hermitian setting, where four Dirac operators and circulant $(4\times 4)$ matrices take over.  
\end{abstract}

 

% ==================================================================
\section{Introduction}
% ==================================================================

\subsection{From Euclidean to Hermitian Clifford analysis}
Clifford analysis is a higher-dimensional function theory built around
the Dirac operator
\[
\dX=\sum_{j=1}^{m}e_{j}\,\partial_{x_{j}},
\]
a first-order, vector-valued, rotation-invariant differential operator
which factorizes the Laplacian, $\Delta_{m}=-\dX^{2}$.  Its null
solutions, the \emph{monogenic} functions, simultaneously generalize
holomorphic functions of one complex variable and refine classical
harmonic analysis, in the sense that monogenicity implies harmonicity
componentwise but not conversely. For a deeper discussion of this
theory we refer the reader to \cite{BDS1982, DSS1992, GHS2008}.

Hermitian Clifford analysis arose as a further refinement of this
picture.  One restricts to even dimension $m=2n$ and introduces a
complex structure $J$ on $\R^{2n}$; the orthogonal Dirac operator
$\dX$ then splits into two \emph{Hermitian Dirac operators} $\dZ$ and
$\dZd$, which are invariant under the action of the unitary group
$\mathrm{U}(n)$ instead of $\mathrm{SO}(2n)$.  Functions annihilated by
$\dZ$ and $\dZd$ simultaneously are called \emph{Hermitian monogenic},
or $h$-monogenic for short.  The systematic development of this theory
was carried out in \cite{BBSES2007I,BDSchS2008};
holomorphic functions of several complex variables and the classical
Dolbeault complex appear as particular instances, which explains part
of the interest of the subject.

By now the subject possesses a considerable internal structure of its
own: a Fischer decomposition of homogeneous polynomials into
$h$-spherical monogenics, a Cauchy--Kovalevskaya extension theory whose
constraints on the initial data motivated the weaker $h$-submonogenic
system, an axiomatic radial algebra, and a complete algebraic analysis
of the underlying overdetermined system; for all of which we refer
to the recent panoramic account of Sabadini and Sommen
\cite{SabSom2025}.  Two facts from that circle of ideas will be needed
here.  First, the free resolution of the module associated with the
Hermitian Dirac system is linear of length $n$ and the cohomology
vanishing which usually produces duality theorems \emph{fails} in this
case (\cite{DES2009}); consequently every duality
statement in the Hermitian setting has to be proved directly, as is
done in \cite{A47,CSS2012} and reviewed in
Section~\ref{sec:duality}.  Second, and in the same spirit, a Hartogs-type extension 
phenomenon does \emph{not} hold for $h$-monogenic
functions \cite{SabSom2025,DES2009}.  It is worth stressing already at
this point that both obstructions dissolve at the level of circulant
matrices: the matrix Dirac operator introduced below has trivial first
syzygies, whence a Poincar\'e lemma and a Mittag--Leffler theorem
\cite{CSS2012}.

The refinement, however, comes at a price.  Since $h$-monogenicity is
an overdetermined system, the naive attempt to write a Cauchy integral
formula fails: neither $\dZ$ nor $\dZd$ possesses a fundamental
solution in the usual sense; the functions
\[
\mathcal{E}(Z)=\frac{2}{\sigma_{2n}}\frac{Z}{|Z|^{2n}},\qquad
\mathcal{E}^{\dagger}(Z)=\frac{2}{\sigma_{2n}}\frac{\Zd}{|Z|^{2n}}
\]
built out of the Euclidean Cauchy kernels are \emph{not} fundamental
solutions of $\dZ$ and $\dZd$ respectively.  The decisive observation,
which organizes essentially all of the work surveyed here, is that the
circulant $(2\times 2)$ matrix
\[
\Ekernel=\begin{pmatrix}
\mathcal{E} & \mathcal{E}^{\dagger}\\[1mm]
\mathcal{E}^{\dagger} & \mathcal{E}
\end{pmatrix}
\]
\emph{is} a fundamental solution of the circulant matrix Dirac operator
\[
\Dop=\begin{pmatrix}
\dZ & \dZd\\[1mm]
\dZd & \dZ
\end{pmatrix}.
\]
Null solutions of $\Dop$ are called \emph{$H$-monogenic} matrix
functions.  Their theory possesses all the classical integral
representation formulae, and $h$-monogenicity is recovered as the
special case of a diagonal matrix.  This \emph{circulant matrix
approach}, introduced in \cite{BDKDS2008,BDKDSS2009}, is the language
in which the whole survey is written.

\subsection{The three axes of the survey}

The papers reviewed here develop the theory along three largely
independent axes, which nevertheless interact in interesting ways.

\medskip
\noindent
\textbf{(i) Geometry of the boundary.}
The starting point is a Jordan domain $\Omega\subset\R^{2n}$ with
boundary $\Gamma$.  When $\Gamma$ is Liapunov or Ahlfors--David
regular, the Cauchy integral over $\Gamma$ makes sense as a genuine
surface integral and the classical machinery applies.  When
$\dimH(\Gamma)>2n-1$, that is, when $\Gamma$ is fractal in the sense of
Mandelbrot, boundary integration loses its meaning.  The remedy,
inspired by the Gauss--Green theorem for fractal boundaries of Harrison
and Norton \cite{HN1992}, is to \emph{define} the Cauchy integral
through a Whitney extension and a Teodorescu transform over $\Omega$,
and to control everything by the $d$-summability of $\Gamma$ and the
$d$-sum of a Whitney decomposition of $\Omega$.  This programme is
carried out in \cite{A25,A35,A36,ABBM2012,A49,A61,A65}.

\medskip
\noindent
\textbf{(ii) The differential operator.}
Besides the Hermitian Dirac operators themselves, one may perturb them
by constants so as to factorize the Helmholtz operator
$\Delta_{2n}+k^{2}$ rather than the Laplacian.  This leads to a
function theory for \emph{Helmholtz $H$-monogenic} matrix functions
\cite{A45}, to its fractal counterpart \cite{A61}, and, most recently,
to the analysis of the associated $\bar\partial$-problem and its
application to Maxwell's equations \cite{A152}.  In the same vein one
may iterate the matrix Dirac operator instead of perturbing it, which
leads to the Riemann problems for $H$-$2$-monogenic matrix functions of
\cite{GF2014,Gu2019}, or keep the operator and make the transmission
condition nonlinear, as in \cite{GuLiu2023}.

\medskip
\noindent
\textbf{(iii) The algebraic structure.}
Replacing the complex structure by a quaternionic one produces
\emph{quaternionic Hermitian Clifford analysis}, centred on the
simultaneous null solutions of \emph{four} Dirac operators in
$\R^{4n}$; the circulant $(2\times 2)$ matrices are then replaced by
circulant $(4\times 4)$ matrices.  Cauchy integral formulae, Cauchy and
Hilbert transforms and boundary value problems have been established in
this setting in \cite{A51,A56,A55,A65,BSLS2020}.

\medskip
Cutting across the three axes is a recurrent methodological pattern,
which we have tried to make visible in the organization of the text:

\begin{enumerate}[label=(\alph*)]
\item one first identifies the correct matrix Cauchy kernel and proves
a Borel--Pompeiu formula;
\item the Cauchy integral is then defined and its boundary behaviour
analysed, producing Plemelj--Sokhotski formulae and a matrix Hilbert
transform $\Hilb$;
\item $\Hilb$ turns out to be bounded and involutive, whence the Hardy-type projections $P^{\pm}=\tfrac12(\mathbf{1}\pm\Hilb)$ split the
boundary data space;
\item every boundary value problem considered---jump, Dirichlet,
Riemann---is then solved by a criterion expressed exactly in terms of
$\Hilb$, $P^{\pm}$, or of the Euclidean Cauchy-type integrals out of
which they are assembled.
\end{enumerate}

\subsection{Scope, and what is not covered}

The survey is organized around integral representation formulae, the
operators built out of them, and the boundary value problems they
solve; the geometry of the boundary is the leading thread.  A number of
adjacent chapters of Hermitian Clifford analysis are therefore left
aside: the Fischer decomposition and the theory of $h$-spherical
monogenics, the Cauchy--Kovalevskaya extension problem and the
$h$-submonogenic system, Hermitian Clifford--Hermite polynomials and
the associated wavelet transforms, the axiomatic radial algebra, the
algebraic analysis of the Hermitian system, and the generalizations of
the framework to superspace and to the symplectic, discrete and
split-quaternionic settings.  Close to our subject, but also outside
it, are the Laurent expansions and $L^{2}$ boundary values of He and Huang
\cite{HH2024} and the half Dirichlet problems of Ku, K\"ahler, Wang \cite{KKW2013}.
For all of this, and for a far more complete bibliography, the reader is referred to \cite{SabSom2025}, 
which allows us to keep our list of references short.

\subsection{Organization of the paper}

Section~\ref{sec:prelim} fixes the algebraic and geometric
preliminaries.  Section~\ref{sec:integral} collects the integral
representation formulae, first in the regular and then in the fractal
case.  Section~\ref{sec:hilbert} is devoted to the Cauchy transform,
the Hilbert transform, the Hardy projections and the H\"older norm
estimates.  Section~\ref{sec:jump} treats jump and decomposition
problems, Section~\ref{sec:dirichlet} the Dirichlet-type problems in
H\"older and in Sobolev scales, and Section~\ref{sec:duality} the
duality theory, first in its Euclidean model (the annihilators of
harmonic differential forms), then at the matricial level and finally at the
scalar level.  Section~8 develops the theory of perturbed and iterated Hermitian Dirac operators: the Hermitian Helmholtz theory, Riemann problems with constant and variable coefficients, for $H$-$2$-monogenic matrix functions and with nonlinear transmission conditions, the fractal case and, in Subsection~8.6, the $\bar\partial$-problem with its application to a two-dimensional Maxwell-type system. Section~9 presents the quaternionic Hermitian setting. Throughout, proofs are omitted; each statement carries a precise
reference to the paper in which it was established.  We have
systematically unified notation across sources; in particular the
spelling ``Hermitian'' is used throughout, although several of the
original papers write ``Hermitean''.\footnote{Titles in the references are quoted as published.}

% ==================================================================
\section{Preliminaries}\label{sec:prelim}
% ==================================================================

\subsection{Clifford algebras and the Dirac operator}

Let $(e_{1},\dots,e_{m})$ be an orthonormal basis of Euclidean space
$\R^{m}$ and let $\Cl$ ($\R_{0,m}$) be the complex (real) Clifford algebra constructed over
$\R^{m}$, whose non-commutative multiplication is governed by
\begin{equation}\label{eq:cliffordrules}
e_{j}^{2}=-1,\quad j=1,\dots,m,\qquad
e_{j}e_{k}+e_{k}e_{j}=0,\quad j\neq k .
\end{equation}
For $A=\{j_{1},\dots,j_{k}\}\subset\{1,\dots,m\}$ with
$j_{1}<\dots<j_{k}$ one puts $e_{A}=e_{j_{1}}\cdots e_{j_{k}}$, and
$e_{\emptyset}=1$; thus every $\tau\in\Cl$ may be written as
$\tau=\sum_{A}\tau_{A}e_{A}$ with $\tau_{A}\in\C$ ($\tau_{A}\in\R$), and
$\dim_{\C}\Cl=2^{m}$.  

The \emph{Hermitian conjugate} of $\tau$ is
\[\tau^{\dagger}=\sum_{A}\tau^{c}_{A}\,\overline{e_{A}},\]
where the bar denotes the main anti-involution of the real Clifford
algebra, for which $\overline{e_{j}}=-e_{j}$, and $\tau^{c}_{A}$ is
the complex conjugate of $\tau_{A}$.  The associated Hermitian inner
product and norm are $\langle\tau,\mu\rangle=[\tau^{\dagger}\mu]_{0}$ and
$|\tau|^{2}=[\tau^{\dagger}\tau]_{0}=\sum_{A}|\tau_{A}|^{2}$.

Euclidean space is embedded in $\Cl$ by identifying
$(x_{1},\dots,x_{m})$ with the real Clifford vector
$X=\sum_{j=1}^{m}e_{j}x_{j}$, for which
$X^{2}=-\langle X,X\rangle=-|X|^{2}$.  The Fischer dual of $X$ is the
Dirac operator $\dX=\sum_{j}e_{j}\partial_{x_{j}}$, which factorizes
the Laplacian, $\Delta_{m}=-\dX^{2}$.  A continuously differentiable
$\Cl$-valued function $g$ on an open set $\Omega\subset\R^{m}$ is
\emph{(left) monogenic} in $\Omega$ if $\dX g=0$ there.  The
fundamental solution of $\dX$ is the Euclidean Cauchy kernel
\begin{equation}\label{eq:Ekernel}
E(X)=-\frac{1}{\sigma_{m}}\frac{X}{|X|^{m}},\qquad X\in\R^{m}\setminus\{0\},
\end{equation}
$\sigma_{m}$ denoting the surface area of the unit sphere in $\R^{m}$.

\subsection{The Hermitian setting}\label{subsec:hermitian}

The transition to the Hermitian framework is based on the introduction
of a \emph{complex structure}, i.e.\ an element $J\in\mathrm{SO}(m)$
with $J^{2}=-\mathbf{1}_{m}$.  Such an element forces $m$ to be even,
so from now on $m=2n$.  A possible realization in terms of the chosen
basis is
\[
J[e_{2j-1}]=-e_{2j},\qquad J[e_{2j}]=e_{2j-1},\qquad j=1,\dots,n .
\]
Associated with $J$ one considers the \emph{twisted} Clifford vector
and the \emph{twisted} Dirac operator
\[
\Xt=J[X]=\sum_{j=1}^{n}\bigl(e_{2j-1}x_{2j}-e_{2j}x_{2j-1}\bigr),
\qquad
\dXt=\sum_{j=1}^{n}\bigl(e_{2j-1}\partial_{x_{2j}}
-e_{2j}\partial_{x_{2j-1}}\bigr),
\]
so that $\dXt$ also factorizes the Laplacian,
$\Delta_{2n}=-\dXt^{2}$.  Introducing the Hermitian Clifford variables
and the Hermitian Dirac operators
\[
Z=\tfrac12\bigl(X+i\Xt\bigr),\qquad
\Zd=-\tfrac12\bigl(X-i\Xt\bigr),
\]
\[
\dZ=-\tfrac14\bigl(\dX-i\,\dXt\bigr),\qquad
\dZd=\tfrac14\bigl(\dX+i\,\dXt\bigr),
\]
one obtains the fundamental factorization of the Laplacian in the
Hermitian setting,
\begin{equation}\label{eq:hermlaplace}
\Delta_{2n}=4\bigl(\dZ\dZd+\dZd\dZ\bigr).
\end{equation}
The operators $\dZ$ and $\dZd$ commute with the action of the group
$\mathrm{Spin}_{J}(2n)=\{s\in\mathrm{Spin}(2n)\mid ss_{J}=s_{J}s\}$,
which realizes the unitary group $\mathrm{U}(n)$ inside the Clifford
algebra; one says that they are invariant under the action of the
unitary group.

\begin{definition}[\cite{BBSES2007I,BBSES2007II}]\label{def:hmono}
A continuously differentiable $\Cln$-valued function $g$ on an open set
$\Omega\subset\R^{2n}$ is called \emph{(left) Hermitian monogenic},
or \emph{$h$-monogenic}, in $\Omega$ if and only if it satisfies there
the system
\[
\dZ g=0=\dZd g,
\qquad\text{or, equivalently,}\qquad
\dX g=0=\dXt g .
\]
Right $h$-monogenicity is defined similarly, and functions which are
both left and right $h$-monogenic are called \emph{two-sided}
$h$-monogenic.
\end{definition}

The second, equivalent, formulation makes it plain that
$h$-monogenicity is a genuine refinement of monogenicity: an
$h$-monogenic function is monogenic with respect to \emph{both} Dirac
operators $\dX$ and $\dXt$, hence in particular harmonic.

\medskip
\noindent
\textbf{The Witt basis and the spinor picture.}
The algebraic mechanism behind the above splitting is best displayed
through the projection operators
$\pi^{\pm}=\tfrac12\bigl(\mathbf{1}_{2n}\pm iJ\bigr)$, whose images
$W^{\pm}=\pi^{\pm}(\R^{2n})$ are maximally isotropic subspaces of
$\C^{2n}$, so that $\C^{2n}=W^{+}\oplus W^{-}$.

\begin{definition}[Witt basis,
\cite{BBSES2007I,BDSchS2008,SabSom2025}]\label{def:witt}
The \emph{Witt basis} of $\Cln$ associated with $J$ is the family
$\{\mathfrak{f}_{j},\mathfrak{f}^{\dagger}_{j}\}_{j=1}^{n}$ given by
\[
\mathfrak{f}_{j}=\pi^{+}[e_{2j-1}]=\tfrac12\bigl(e_{2j-1}-ie_{2j}\bigr),
\qquad
\mathfrak{f}^{\dagger}_{j}=-\pi^{-}[e_{2j-1}]
=-\tfrac12\bigl(e_{2j-1}+ie_{2j}\bigr),
\]
$j=1,\dots,n$.  Its elements obey the Grassmann and duality relations
\[
\mathfrak{f}_{j}\mathfrak{f}_{k}+\mathfrak{f}_{k}\mathfrak{f}_{j}=0,
\qquad
\mathfrak{f}^{\dagger}_{j}\mathfrak{f}^{\dagger}_{k}
+\mathfrak{f}^{\dagger}_{k}\mathfrak{f}^{\dagger}_{j}=0,
\qquad
\mathfrak{f}_{j}\mathfrak{f}^{\dagger}_{k}
+\mathfrak{f}^{\dagger}_{k}\mathfrak{f}_{j}=\delta_{jk},
\]
in particular
$\mathfrak{f}_{j}^{2}=\bigl(\mathfrak{f}^{\dagger}_{j}\bigr)^{2}=0$.
\end{definition}

Writing $z_{j}=x_{2j-1}+ix_{2j}$ for the complex coordinates
associated with $J$, and $\partial_{z_{j}}$, $\partial_{z^{c}_{j}}$ for
the corresponding Cauchy--Riemann operator and its conjugate, one has
the compact expressions
\[
Z=\sum_{j=1}^{n}\mathfrak{f}_{j}z_{j},\qquad
\Zd=\sum_{j=1}^{n}\mathfrak{f}^{\dagger}_{j}z^{c}_{j},
\qquad
\dZ=\sum_{j=1}^{n}\mathfrak{f}^{\dagger}_{j}\partial_{z_{j}},\qquad
\dZd=\sum_{j=1}^{n}\mathfrak{f}_{j}\partial_{z^{c}_{j}},
\]
together with $X=Z-\Zd$, $\Xt=\tfrac1i(Z+\Zd)$,
$\dX=2(\dZd-\dZ)$ and $\dXt=\tfrac2i(\dZd+\dZ)$.  The Witt relations
give at once the isotropy
\[
Z^{2}=\bigl(\Zd\bigr)^{2}=0,\qquad
\dZ^{2}=\dZd^{2}=0,
\]
which is precisely what prevents $\dZ$ and $\dZd$ from admitting
fundamental solutions of the usual type, and which forces the passage
to matrices carried out in Section~\ref{subsec:circulant}.

Setting $I_{j}=\mathfrak{f}_{j}\mathfrak{f}^{\dagger}_{j}$ one obtains
mutually commuting self-adjoint idempotents; with $I=I_{1}\cdots I_{n}$
the \emph{complex spinor space} is
$\mathbb{S}=\Cln I\cong\Lambda^{\dagger}_{n}I$, where
$\Lambda^{\dagger}_{n}=\mathrm{Alg}_{\C}
\{\mathfrak{f}^{\dagger}_{1},\dots,\mathfrak{f}^{\dagger}_{n}\}$, and it
decomposes as
$\mathbb{S}=\bigoplus_{j=0}^{n}\mathbb{S}^{j}$ with
$\mathbb{S}^{j}=\bigl(\Lambda^{\dagger}_{n}\bigr)^{(j)}I$.  The
Hermitian Dirac operators shift this grading:
$\dZ$ maps $\mathbb{S}^{j}$ into $\mathbb{S}^{j+1}$ and $\dZd$ maps
$\mathbb{S}^{j}$ into $\mathbb{S}^{j-1}$, with the convention
$\mathbb{S}^{-1}=\mathbb{S}^{n+1}=0$.  This elementary observation is
the seed of the algebraic analysis recalled in
Subsection~\ref{subsec:algebraic}.

\begin{proposition}[\cite{SabSom2025}]\label{prop:spinsplit}
Let $g=\sum_{j=0}^{n}g_{j}$ be the decomposition of an
$\mathbb{S}$-valued function into its homogeneous spinor parts,
$g_{j}$ taking values in $\mathbb{S}^{j}$.  Then $\dZ g=0$ is
equivalent to the system $\dZ g_{j}=0$, $j=0,\dots,n$, and likewise for
$\dZd$.
\end{proposition}

\begin{example}[\cite{SabSom2025}]\label{ex:elementary}
For $F$ holomorphic in an open subset of the $z_{j}$-plane, the
function $g=\mathfrak{f}^{\dagger}_{j}F(z_{j})$ is $h$-monogenic;
similarly $g=\mathfrak{f}_{j}F(z^{c}_{j})$ is $h$-monogenic for $F$
anti-holomorphic.  Both statements are immediate from the Witt relations of 
Definition~\ref{def:witt}.
\end{example}

\begin{proposition}[\cite{BDSchS2008,SabSom2025}]\label{prop:topspinor}
Let $g$ take values in the top piece $\mathbb{S}^{n}$, so that
$g=g_{n}\,\mathfrak{f}^{\dagger}_{1}\cdots\mathfrak{f}^{\dagger}_{n}I$
with $g_{n}$ scalar valued.  Then $g$ is $h$-monogenic if and only if
$g_{n}$ is holomorphic in $(z_{1},\dots,z_{n})$.
\end{proposition}

Proposition~\ref{prop:topspinor} is the precise sense in which
Hermitian Clifford analysis contains the function theory of several
complex variables; as we shall recall in
Remark~\ref{rem:MB}, the Cauchy formula of Section~\ref{sec:integral}
reduces, on this component, to the classical Martinelli--Bochner
formula.

\subsection{The circulant matrix framework}\label{subsec:circulant}

Starting from the pair $(E,E^{|})$ of fundamental solutions of $\dX$
and $\dXt$, where
\[
E(X)=-\frac{1}{\sigma_{2n}}\frac{X}{|X|^{2n}},\qquad
E^{|}(X)=-\frac{1}{\sigma_{2n}}\frac{\Xt}{|X|^{2n}},
\]
one introduces
\begin{equation}\label{eq:hermkernels}
\mathcal{E}=-\bigl(E+iE^{|}\bigr)=\frac{2}{\sigma_{2n}}\frac{Z}{|Z|^{2n}},
\qquad
\mathcal{E}^{\dagger}=E-iE^{|}=\frac{2}{\sigma_{2n}}\frac{\Zd}{|Z|^{2n}},
\end{equation}
where, as in \cite{BDKDSS2009}, $|Z|=|\Zd|:=|X|$.

As already emphasized, $\mathcal{E}$ and $\mathcal{E}^{\dagger}$ are
\emph{not} fundamental solutions of $\dZ$ and $\dZd$.  However, the
following holds.

\begin{theorem}[\cite{BDKDSS2009,A21,A35}]\label{thm:fundsol}
Let $\Dop$ and $\Ekernel$ be the circulant $(2\times 2)$ matrices
\[
\Dop=\begin{pmatrix}\dZ&\dZd\\ \dZd&\dZ\end{pmatrix},
\qquad
\Ekernel=\begin{pmatrix}
\mathcal{E}&\mathcal{E}^{\dagger}\\
\mathcal{E}^{\dagger}&\mathcal{E}\end{pmatrix},
\qquad
\boldsymbol{\delta}=\begin{pmatrix}\delta&0\\0&\delta\end{pmatrix},
\]
$\delta$ being the Dirac distribution.  Then
$\Dop\,\Ekernel=\boldsymbol{\delta}$, i.e.\ $\Ekernel$ is a fundamental
solution of $\Dop$ in the matricial sense.
\end{theorem}

\begin{remark}\label{rem:conventions}
The labelling of the pair $(\dZ,\dZd)$ is not uniform in the
literature and several sources interchange the two symbols; all
formulae in this survey are written in the convention just fixed, which
is the one of \cite{BDKDSS2009,A45,CSS2012,SabSom2025}.  With it one
has, on $\R^{2n}\setminus\{0\}$,
\[
\dZd\,\mathcal{E}=0,\qquad \dZ\,\mathcal{E}^{\dagger}=0,
\]
for the kernels $\mathcal{E},\mathcal{E}^{\dagger}$ of
\eqref{eq:hermkernels}, while $\dZ\mathcal{E}$ and
$\dZd\mathcal{E}^{\dagger}$ do not vanish; this pair of identities is
what makes the matrix $\Ekernel$ of Theorem~\ref{thm:fundsol} a
fundamental solution, and it will be used repeatedly in
Section~\ref{sec:duality}.
\end{remark}

We denote by $\CM$ the set of circulant $(2\times 2)$ matrices over
$\Cln$.  It is closed under sums and products and, defining the
Hermitian conjugate $A^{\dagger}$ entrywise, $\CM$ becomes a right
linear associative algebra over $\Cln$ endowed with an anti-involution.
Given $\Cln$-valued functions $g_{1},g_{2}$ we write
\[
\Gmat=\begin{pmatrix}g_{1}&g_{2}\\ g_{2}&g_{1}\end{pmatrix},
\qquad
G_{0}=\begin{pmatrix}g&0\\ 0&g\end{pmatrix},
\]
$\Zero$ standing for the zero matrix and $\mathbf{1}$ for the identity
of $\CM$.

\begin{definition}[\cite{A21,A35,A45}]\label{def:Hmono}
The matrix function $\Gmat$ is called \emph{(left) $H$-monogenic} in
$\Omega$ if $\Dop[\Gmat]=\Zero$ in $\Omega$; right $H$-monogenicity is
defined by $\Gmat\,\Dop=\Zero$.  We write $\HM_{l}(\Omega)$, resp.\ $\HM_{r}(\Omega)$, for the corresponding spaces.
\end{definition}

Explicitly, left $H$-monogenicity of $\Gmat$ reads
\begin{equation}\label{eq:Hmonosystem}
\dZ g_{1}+\dZd g_{2}=0,\qquad
\dZd g_{1}+\dZ g_{2}=0 .
\end{equation}

In general $H$-monogenicity of $\Gmat$ does \emph{not} imply
$h$-monogenicity of the entries $g_{1},g_{2}$.  The crucial exception
is the diagonal case: $G_{0}$ is $H$-monogenic if and only if $g$ is
$h$-monogenic.  Moreover, from
$4\,\Dop\,\Dop^{\dagger}=4\,\Dop^{\dagger}\Dop=\boldsymbol{\Delta}
:=\mathrm{diag}(\Delta_{2n},\Delta_{2n})$
one sees that every $H$-monogenic matrix function of class $C^{2}$ is
harmonic, whence its entries are harmonic in the usual sense.

Put $\Pi=\bigl(\begin{smallmatrix}0&1\\1&0\end{smallmatrix}\bigr)$ and $\mathcal D^\dagger_{(Z,Z^\dagger)}=\bigl(\begin{smallmatrix}\partial_{Z^\dagger}&\partial_Z\\ \partial_Z&\partial_{Z^\dagger}\end{smallmatrix}\bigr)$. Then $\mathcal D^\dagger_{(Z,Z^\dagger)}=\Pi\,\mathcal D_{(Z,Z^\dagger)}=\mathcal D_{(Z,Z^\dagger)}\Pi$, so that 
\[
\bigl(\mathcal D_{(Z,Z^\dagger)}\bigr)^2=\tfrac14\,\Delta_{2n}\,\Pi.
\]
Throughout the paper we abbreviate
\begin{equation}\label{eq:lambdan}
\lam:=(-1)^{n}(2i)^{n}
\end{equation}
a constant which arises from the translation between the Euclidean and
the Hermitian volume and surface elements and which appears in
virtually every integral formula below\footnote{In \cite{BDKDSS2009}, where $J[e_j]=-e_{n+j}$ and $z_j=x_j+ix_{n+j}$, this constant reads $(-1)^{\frac{n(n+1)}{2}}(2i)^{n}$; the additional sign $(-1)^{\frac{n(n-1)}{2}}$ comes from reordering $dx_1\wedge dx_{n+1}\wedge\dots\wedge dx_n\wedge dx_{2n}$ into the standard orientation $dx_1\wedge\dots\wedge dx_{2n}$. Several of the papers reviewed here (e.g.\ \cite{A35,A45}) use the adjacent-pair complex structure of Subsection~\ref{subsec:hermitian} but keep the constant of \cite{BDKDSS2009}. With the adjacent-pair structure and the standard orientation the correct value is $(-2i)^{n}$; the two values coincide for $n\equiv 0,1\pmod 4$ and differ by a sign for $n\equiv 2,3\pmod 4$. All formulae below are stated with \eqref{eq:lambdan}.}.

\subsection{Function spaces}\label{subsec:fspaces}

Continuity, differentiability and integrability of a matrix function
$\Gmat\in\CM$ are understood entrywise; thus $\Gmat\in C^{k}(E)$,
$C^{0,\nu}(E)$, $L^{p}(E)$, $W^{k}_{2}(E)$ mean that each entry belongs
to the corresponding scalar space.  Equivalently, introducing the
non-negative function
\[
\|\Gmat(X)\|=\max\bigl\{|g_{1}(X)|,|g_{2}(X)|\bigr\},
\]
these classes may be described by the traditional conditions
\[
|\Gmat|_{\nu,E}
:=\sup_{X,Y\in E,\,X\neq Y}
\frac{\|\Gmat(X)-\Gmat(Y)\|}{|X-Y|^{\nu}}<+\infty,
\qquad
\|\Gmat\|_{p}
:=\bigl(\int_E\|G_2^1(X)\|^p\,dV(X)\bigr)^{1/p}<+\infty,
\]
and the norm of $\Gmat\in C^{0,\nu}(E)$ is
$\|\Gmat\|_{\nu,E}=\max_{X\in E}\|\Gmat(X)\|+|\Gmat|_{\nu,E}$.
The Hermitian Hilbert--Sobolev spaces $W^{k}_{2}(\Omega)$ are equipped
with the $\CM$-valued inner product
\begin{equation}\label{eq:innerprod}
\bigl\langle \mathfrak{U}^{1}_{2},\mathfrak{V}^{1}_{2}\bigr\rangle_{L_{2}}
=\int_{\Omega}\Bigl(\mathfrak{U}^{1}_{2}\Bigr)^{\dagger}\mathfrak{V}^{1}_{2}\,dV ,
\end{equation}
used in Section~\ref{sec:dirichlet}.

\subsection{Boundary regularity: from Liapunov to fractal}
\label{subsec:geometry}

Let $\Omega\subset\R^{2n}$ be a \emph{Jordan domain}, i.e.\ a bounded
oriented connected open set whose boundary $\Gamma$ is a compact
topological surface; we put $\Omega^{+}=\Omega$ and
$\Omega^{-}=\R^{2n}\setminus\overline{\Omega}$, both assumed connected.
For $n=1$ this reduces to the classical notion of a Jordan domain in
the plane.  Three levels of regularity of $\Gamma$ occur in the papers
under review.

\medskip
\noindent
\textbf{(L) Liapunov / piecewise smooth / Lipschitz.}  The unit normal
$\eta(X)$ exists almost everywhere.  This is the setting in which
the classical Plemelj--Sokhotski formulae are available and is used in
\cite{A45,A51,A55,A56}.

\medskip
\noindent
\textbf{(AD) Ahlfors--David regular.}  $\Gamma$ is
\emph{Ahlfors--David regular} if there is $c_{AD}>0$ with
\begin{equation}\label{eq:AD}
c_{AD}^{-1}r^{2n-1}\le
\mathcal{H}^{2n-1}\bigl(\Gamma\cap \overline{B}(X,r)\bigr)
\le c_{AD}\,r^{2n-1}
\end{equation}
for all $X\in\Gamma$ and $0<r\le\operatorname{diam}\Gamma$, where
$\mathcal{H}^{s}$ denotes $s$-dimensional Hausdorff measure.  If only
the left-hand inequality is required (together with
$\mathcal{H}^{2n-1}(\Gamma)<\infty$) one speaks of \emph{left}
Ahlfors--David regularity, abbreviated $(\ell)$AD-regularity.  Here the
normal vector in the sense of Federer \cite{Federer} replaces the
classical one.  This is the setting of \cite{A19,A21}.

\medskip
\noindent
\textbf{(F) Fractal, controlled by $d$-summability.}  Let $\mathcal{K}\subset
\R^{2n}$ be compact and let $N_{\mathcal{K}}(\tau)$ be the least number of
$\tau$-balls needed to cover $\mathcal{K}$.

\begin{definition}[\cite{HN1992}]\label{def:dsummable}
The compact set $\mathcal{K}$ is \emph{$d$-summable} if the improper integral
\[
\int_{0}^{1}N_{\mathcal{K}}(\tau)\,\tau^{\,d-1}\,d\tau
\]
converges.
\end{definition}

\begin{lemma}[\cite{HN1992}]\label{lem:dsummable}
\textup{(i)} Any $d$-summable set $\mathcal{K}$ has box dimension
$\alpha(\mathcal{K})\le d$.  \textup{(ii)} If $\alpha(\mathcal{K})<d$ then $\mathcal{K}$ is
$d$-summable.  \textup{(iii)} If $\mathcal{K}$ is $d$-summable then it is
$(d+\varepsilon)$-summable for every $\varepsilon>0$.
\end{lemma}

One assumes throughout the fractal part of the theory that
\begin{equation}\label{eq:dimhyp}
2n-1\le \dimH(\Gamma)\le\alpha(\Gamma)<2n,
\end{equation}
which includes the case where $\Gamma$ is fractal in the sense of
Mandelbrot, i.e.\ $2n-1<\dimH(\Gamma)$.  By
Lemma~\ref{lem:dsummable} there is then always $d\in(2n-1,2n)$ such
that $\Gamma$ is $d$-summable.

The technical device that makes the analysis possible is the
\emph{Whitney decomposition} $\mathcal{W}$ of $\Omega$, obtained by
starting from the mesh $\mathcal{M}_{1}$ of unit cubes of the lattice
$\Z^{2n}$ meeting $\Omega$, recursively bisecting, and setting
\[
\mathcal{W}_{1}=\{Q\in\mathcal{M}_{1}\mid
\text{all neighbouring cubes of }Q\text{ belong to }\Omega\},
\]
\[
\mathcal{W}_{k}=\{Q\in\mathcal{M}_{k}\mid
\text{all neighbouring cubes of }Q\text{ belong to }\Omega,\ \text{and }
\nexists\,Q^{*}\in\mathcal{W}_{k-1}\text{ with }Q\subset Q^{*}\},
\quad k\ge 2 .
\]

\begin{lemma}[\cite{HN1992}]\label{lem:dsum}
If $\Omega$ is a Jordan domain of $\R^{2n}$ with $d$-summable boundary
$\Gamma$, then the \emph{$d$-sum}
\[
s(d)=\sum_{Q\in\mathcal{W}}|Q|^{\,d}
\]
of the Whitney decomposition $\mathcal{W}$ of $\Omega$ is finite.
\end{lemma}

The link between the two notions is provided by the following
statement, which is what allows the fractal results to contain the
classical ones as particular cases.

\begin{lemma}[\cite{A49}]\label{lem:ADdsum}
If $\Gamma$ is $(\ell)$AD-regular, then $\Gamma$ is
$(2n-1+\varepsilon)$-summable for every $\varepsilon>0$; moreover
$s(d)\le c\,\mathcal{H}^{2n-1}(\Gamma)/\varepsilon$ for
$d=2n-1+\varepsilon$.
\end{lemma}

Finally, we shall repeatedly need the Whitney extension theorem in its
matricial form.

\begin{theorem}[Whitney extension theorem, \cite{Stein}, matrix form
\cite{A35,A61,A65}]\label{thm:whitney}
Let $\mathcal{K}\subset\R^{2n}$ be compact and $\Gmat\in C^{0,\nu}(\mathcal{K})$.  Then
there exists a compactly supported matrix function $\varUpsilon[{G}^{1}_{2}]$
satisfying
\begin{enumerate}[label=\upshape(\roman*)]
\item $\varUpsilon[{G}^{1}_{2}]\big|_{E}=\Gmat$;
\item $\varUpsilon[{G}^{1}_{2}]\in C^{1}\bigl(\R^{2n}\setminus\mathcal{K}\bigr)$;
\item $\bigl\|\Dop\varUpsilon[{G}^{1}_{2}](X)\bigr\|
\le c\,\operatorname{dist}(X,\mathcal{K})^{\nu-1}$ for
$X\in\R^{2n}\setminus \mathcal{K}$.
\end{enumerate}
Any such extension is called a \emph{Whitney-type extension} of
$\Gmat$.
\end{theorem}

% ==================================================================
\section{Integral representation formulae}\label{sec:integral}
% ==================================================================

\subsection{Differential forms: the Hermitian Stokes and Cauchy
theorems}\label{subsec:forms}

Before turning to the operators which carry the analysis on
$\Gamma$, it is useful to record the version of the basic identities
which is phrased in the language of differential forms, since it is the
one needed in the duality theory of Section~\ref{sec:duality}.
Following \cite{RSS2002,BDKDSS2009,SabSom2025} one introduces the
$(n-1,n)$- and $(n,n-1)$-forms
\[d\sigma_{Z}=\sum_{j=1}^{n}\mathfrak{f}^{\dagger}_{j}
\bigl(dz_{1}\wedge dz^{c}_{1}\bigr)\wedge\dots\wedge
\bigl(\widehat{dz_{j}}\wedge dz^{c}_{j}\bigr)\wedge\dots\wedge
\bigl(dz_{n}\wedge dz^{c}_{n}\bigr),\]
\[d\sigma_{\Zd}=\sum_{j=1}^{n}\mathfrak{f}_{j}
\bigl(dz_{1}\wedge dz^{c}_{1}\bigr)\wedge\dots\wedge
\bigl(dz_{j}\wedge\widehat{dz^{c}_{j}}\bigr)\wedge\dots\wedge
\bigl(dz_{n}\wedge dz^{c}_{n}\bigr),\]
the hat indicating that the corresponding factor is omitted, together
with the volume form
\[dW(Z,\Zd)=\bigl(dz_{1}\wedge dz^{c}_{1}\bigr)\wedge\dots\wedge
\bigl(dz_{n}\wedge dz^{c}_{n}\bigr),\]
which will reappear in \eqref{eq:teo}.

\begin{theorem}[Hermitian Clifford--Stokes theorems,
\cite{RSS2002,BDKDSS2009,SabSom2025}]\label{thm:hstokes}
Let $\Omega$ be open in $\R^{2n}$, let $f,g\in C^{1}(\Omega,\Cln)$ and
let $\Sigma\subset\Omega$ be a $2n$-dimensional compact oriented
differentiable manifold with smooth boundary $\partial\Sigma$.  Then
\[\int_{\partial\Sigma}f\,d\sigma_{Z}\,g=\int_{\Sigma}\bigl[(f\dZ)g+f(\dZ g)\bigr]\,dW(Z,\Zd),\]
\[\int_{\partial\Sigma}f\,\bigl(-d\sigma_{\Zd}\bigr)\,g=\int_{\Sigma}\bigl[(f\dZd)g+f(\dZd g)\bigr]\,dW(Z,\Zd).\]
\end{theorem}

\begin{theorem}[Hermitian Cauchy theorems,
\cite{BDKDSS2009,SabSom2025}]\label{thm:hcauchythm}
Under the hypotheses of Theorem~\ref{thm:hstokes}, if $g$ is $h$-monogenic in $\Omega$ then
\[\int_{\partial\Sigma}d\sigma_{Z}\,g=0,
\qquad
\int_{\partial\Sigma}d\sigma_{\Zd}\,g=0.\]
\end{theorem}

As was stressed in the introduction, no Cauchy \emph{formula} for
$h$-monogenic functions can be extracted from these identities; the
correct statement lives at the matrix level.  

From now on $Y$ and $\Yt$ denote Clifford vectors associated with
points of $\Omega^{\pm}$, and their Hermitian counterparts are
$V=\tfrac12(Y+i\Yt)$, $V^{\dagger}=-\tfrac12(Y-i\Yt)$.

Writing
\[\Ekernel_{\sharp}=\frac{1}{\lam}
\begin{pmatrix}
\mathcal{E}&\mathcal{E}^{\dagger}\\
\mathcal{E}^{\dagger}&\mathcal{E}
\end{pmatrix},
\qquad
d\Sigma_{(Z,\Zd)}=
\begin{pmatrix}
d\sigma_{Z}&-d\sigma_{\Zd}\\
-d\sigma_{\Zd}&d\sigma_{Z}
\end{pmatrix},\]
one has the following statement, which is the differential-form incarnation of Corollary~\ref{cor:cauchy} below.

\begin{theorem}[Matrix Cauchy formula in form language,
\cite{BDKDSS2009,CSS2012,SabSom2025}]\label{thm:cauchyforms}
Let $\Omega$ be open in $\R^{2n}$, let $\Sigma\subset\Omega$ be as in
Theorem~\ref{thm:hstokes} and let $Y$ lie in the interior of $\Sigma$.
If $\Gmat\in\HM_{l}(\Omega)$ then
\[\Gmat(Y)=\int_{\partial\Sigma}\Ekernel_{\sharp}(Z-V)\, d\Sigma_{(Z,\Zd)}\,\Gmat(X),\]
while for $\Gmat\in\HM_{r}(\Omega)$
\[\Gmat(Y)=\int_{\partial\Sigma}\Gmat(X)\,d\Sigma_{(Z,\Zd)}\, \Ekernel_{\sharp}(Z-V).\]
In particular the formula applies to a diagonal matrix $G_{0}$, that
is, to an $h$-monogenic function $g$.
\end{theorem}

\begin{remark}\label{rem:MB}
Two comments on normalizations and on the scalar content of
Theorem~\ref{thm:cauchyforms}.  (i) The kernel $\Ekernel_{\sharp}$
differs from the kernel $\Ekernel$ of Theorem~\ref{thm:fundsol} by the
factor $\frac{1}{\lam}$; this is exactly the bookkeeping which, when
$d\Sigma_{(Z,\Zd)}$ is rewritten in terms of the surface measure
$N_{(Z,\Zd)}\,d\mathcal{H}^{2n-1}$ of Subsection~\ref{subsec:regular},
produces the constant $\lam$ of \eqref{eq:lambdan} in all the formulae
below.  (ii) If $g$ takes values in the top spinor space
$\mathbb{S}^{n}$, so that by Proposition~\ref{prop:topspinor} it is
$h$-monogenic precisely when its scalar coefficient is holomorphic in
$(z_{1},\dots,z_{n})$, then Theorem~\ref{thm:cauchyforms} reduces to
the classical Martinelli--Bochner formula of several complex variables
\cite{BDKDSS2009,SabSom2025}; compare with Remark~\ref{rem:qMB} in the
quaternionic setting.
\end{remark}

\subsection{The regular case}\label{subsec:regular}
The two basic integral operators of the theory are the \emph{Hermitian
Cauchy integral}
\begin{equation}\label{eq:cauchyint}
\Cauchy\bigl[\Gmat\bigr](Y)
=\int_{\Gamma}\Ekernel(Z-V)\,N_{(Z,\Zd)}\,\Gmat(X)\,
d\mathcal{H}^{2n-1},
\qquad Y\in\Omega^{\pm},
\end{equation}
and the \emph{Hermitian Teodorescu transform}
\begin{equation}\label{eq:teo}
\Teo\bigl[\Fmat\bigr](Y)
=-\int_{\Omega}\Ekernel(Z-V)\,\Fmat(X)\,dW(Z,\Zd),
\qquad Y\in\R^{2n}.
\end{equation}
Here the circulant matrix
\[
N_{(Z,\Zd)}=\begin{pmatrix}N&-N^{\dagger}\\ -N^{\dagger}&N\end{pmatrix}
\]
contains, up to a constant factor, the Hermitian projections $N$ and
$-N^{\dagger}$ of the unit normal vector $\eta(X)$ at $X\in\Gamma$,
$d\mathcal{H}^{2n-1}=\mathrm{diag}(d\mathcal{H}^{2n-1},
d\mathcal{H}^{2n-1})$, and the volume element $dW(Z,\Zd)$, reflecting
integration in the underlying complex planes, is defined through
\begin{equation}\label{eq:dVdW}
dV(X)=\Bigl(\frac{i}{2}\Bigr)^{n}dW(Z,\Zd),\qquad\text{i.e.}\qquad dW(Z,\Zd)=\lam\,dV(X),
\end{equation}
where $dV=dx_{1}\wedge dx_{2}\wedge\dots\wedge dx_{2n}$ carries the standard orientation.

Finally, $\boldsymbol{\chi}_{\Omega}
=\mathrm{diag}(\chi_{\Omega},\chi_{\Omega})$ denotes the matricial
characteristic function of $\Omega$; writing the right-hand side of an
integral formula as $\lam\,\boldsymbol{\chi}_{\Omega}\Gmat$ is thus a
compact way of distinguishing the two cases $Y\in\Omega^{+}$ and
$Y\in\Omega^{-}$, the value being $\lam\Gmat(Y)$ in the first and $0$
in the second.

\begin{theorem}[Hermitian Borel--Pompeiu formula,
\cite{BDKDSS2009,A35,A45}]\label{thm:BP}
Let $\Gmat\in C^{1}(\overline{\Omega})$, $\Omega$ a domain with
piecewise $C^{\infty}$ smooth boundary $\Gamma$.  Then
\[
\Cauchy\bigl[\Gmat\bigr](Y)+\Teo\bigl[\Dop\Gmat\bigr](Y)
=\lam\,\boldsymbol{\chi}_{\Omega}(Y)\,\Gmat(Y),
\qquad Y\in\R^{2n}\setminus\Gamma .
\]
\end{theorem}

\begin{corollary}[Hermitian Cauchy formula, \cite{BDKDSS2009,A45}]
\label{cor:cauchy}
If moreover $\Gmat\in\HM(\Omega)$, then
$\Cauchy[\Gmat](Y)=\lam\,\boldsymbol{\chi}_{\Omega}(Y)\Gmat(Y)$ for
$Y\in\R^{2n}\setminus\Gamma$.
\end{corollary}

The next statement expresses the fundamental property of the matricial
Teodorescu transform, namely that it is an algebraic right inverse of
$\Dop$ up to the factor $\lam$.

\begin{theorem}[\cite{A36,A45}]\label{thm:teoinverse}
If $\Gmat\in C^{1}(\Omega)\cap C(\overline{\Omega})$, then
$\Dop\,\Teo[\Gmat](Y)=\lam\,\boldsymbol{\chi}_{\Omega}(Y)\Gmat(Y)$.
Moreover \cite{A36}:
\begin{enumerate}[label=\upshape(\roman*)]
\item if $\Gmat\in L^{p}(\Omega)$ with $p>2n$, then
$\Teo[\Gmat]\in C^{0,\frac{p-2n}{p}}(\R^{2n})$;
\item the displayed inversion formula already holds for
$\Gmat\in L^{1}(\Omega)$ in the sense of distributions.
\end{enumerate}
\end{theorem}

\subsection{Fractal boundaries: the Cauchy integral revisited}

When $\dimH(\Gamma)>2n-1$ the surface integral \eqref{eq:cauchyint}
ceases to be meaningful.  The way out, proposed in \cite{A35}, is to
\emph{take the Borel--Pompeiu formula as a definition}: one replaces
$\Gmat$ by a Whitney-type extension and moves the boundary integral to
the left-hand side.

\begin{definition}[Fractal Hermitian Cauchy integral, \cite{A35}]
\label{def:fractalcauchy}
Let $d\in(2n-1,2n)$ and $d-2n+1<\nu\le 1$, and let $\Omega$ be a Jordan
domain with $d$-summable boundary $\Gamma$.  For
$\Gmat\in C^{0,\nu}(\Gamma)$ define
\begin{equation}\label{eq:fractalcauchy}
\Cauchy^{*}\bigl[\Gmat\bigr](Y)
=\lam\,\boldsymbol{\chi}_{\Omega}(Y)\,\varUpsilon[{G}^{1}_{2}](Y)
-\Teo\bigl[\Dop\varUpsilon[{G}^{1}_{2}]\bigr](Y),
\qquad Y\in\R^{2n}\setminus\Gamma,
\end{equation}
where $\varUpsilon[{G}^{1}_{2}]$ is a Whitney-type extension of $\Gmat$.
\end{definition}

\begin{proposition}[\cite{A35}]\label{prop:welldefined}
The matrix function \eqref{eq:fractalcauchy} is well defined for every
$Y\in\R^{2n}\setminus\Gamma$ and does not depend on the particular
choice of the Whitney-type extension of $\Gmat$.
\end{proposition}

\begin{remark}[\cite{A35}]
Definition~\ref{def:fractalcauchy} remains valid under the condition
$\nu>\alpha(\Gamma)-2n+1$ for any $\Gamma$ of box dimension
$\alpha(\Gamma)$.  When $\Gamma$ is sufficiently regular, for
instance Ahlfors--David regular, the expression
\eqref{eq:fractalcauchy} reduces to the classical Cauchy integral
\eqref{eq:cauchyint}, as follows from Theorem~\ref{thm:BP}.
\end{remark}

With this definition the whole family of integral representation
formulae survives in the fractal setting.

\begin{theorem}[Hermitian Borel--Pompeiu formula on fractal domains,
\cite{A35}]\label{thm:fractalBP}
Let $\Omega$ be a Jordan domain with $d$-summable boundary $\Gamma$,
$d\in(2n-1,2n)$, and let $\Gmat\in C^{1}(\overline{\Omega})$.  Then
\[
\Cauchy^{*}\bigl[\Gmat\bigr](Y)+\Teo\bigl[\Dop\Gmat\bigr](Y)
=\lam\,\boldsymbol{\chi}_{\Omega}(Y)\,\Gmat(Y).
\]
\end{theorem}

\begin{theorem}[Hermitian Koppelman formula, \cite{A35}]
\label{thm:koppelman}
Under the hypotheses of Theorem~\ref{thm:fractalBP},
\[
\Cauchy^{*}\bigl[\Gmat\bigr](Y)+\Teo\bigl[\Dop\Gmat\bigr](Y)
+\Dop\,\Teo\bigl[\Gmat\bigr](Y)
=2\lambda_n\,
\boldsymbol{\chi}_{\Omega}(Y)\,\Gmat(Y).
\]
\end{theorem}

\begin{theorem}[Hermitian Cauchy formula on fractal domains, \cite{A35}]
\label{thm:fractalcauchyformula}
Under the hypotheses of Theorem~\ref{thm:fractalBP}, if $\Gmat$ is
moreover $H$-monogenic in $\Omega$, then
$\Cauchy^{*}[\Gmat](Y)=\lam\,\boldsymbol{\chi}_{\Omega}(Y)\Gmat(Y)$.
\end{theorem}

% ==================================================================
\section{The Hermitian Hilbert transform}\label{sec:hilbert}
% ==================================================================

\subsection{Plemelj--Sokhotski formulae and the Hilbert transform}
\label{subsec:PS}

In the Euclidean setting one has two Cauchy-type integrals at one's disposal:
${C}_{\Gamma}$ and ${C}^{|}_{\Gamma}$, associated with
$\dX$ and $\dXt$, together with the corresponding singular integral operators
$S_{\Gamma}$ and $S^{|}_{\Gamma}$; their Plemelj--Sokhotski formulae
read
\begin{equation}\label{eq:euclPS}
\bigl({C}^{\pm}_{\Gamma}f\bigr)(X)
=\tfrac12\bigl(S_{\Gamma}f(X)\pm f(X)\bigr),
\qquad
\bigl({C}^{|,\pm}_{\Gamma}f\bigr)(X)
=\tfrac12\bigl(S^{|}_{\Gamma}f(X)\pm f(X)\bigr).
\end{equation}
The matricial \emph{singular Cauchy integral}, equivalently the
\emph{Hermitian Hilbert transform}, is the circulant combination
\begin{equation}\label{eq:hilbertdef}
\Hilb=\frac12
\begin{pmatrix}
S_{\Gamma}+S^{|}_{\Gamma} & -S_{\Gamma}+S^{|}_{\Gamma}\\[1mm]
-S_{\Gamma}+S^{|}_{\Gamma} & S_{\Gamma}+S^{|}_{\Gamma}
\end{pmatrix},
\end{equation}
its action on a matrix function being given by matrix multiplication
followed by the action of the entries as operators.  (In \cite{A21}
this operator is denoted $\mathcal{S}$; we use the notation $\Hilb$ of
\cite{A45,A59} throughout.)

\begin{theorem}[Matricial Plemelj--Sokhotski formulae,
\cite{A21,A45}]\label{thm:PS}
Let $\Gamma$ be $(\ell)$AD-regular and let $\Gmat\in C^{0,\nu}(\Gamma)$, $0<\nu<1$.  Then the continuous boundary values of the
Hermitian Cauchy integral exist and are given, at every $U\in\Gamma$,
by
\[
\Cauchy^{\pm}\bigl[\Gmat\bigr](U)
=\lim_{\substack{Y\to U\\ Y\in\Omega^{\pm}}}
\Cauchy\bigl[\Gmat\bigr](Y)
=\lam\left(\pm\tfrac12\,\Gmat(U)
+\tfrac12\,\Hilb\bigl[\Gmat\bigr](U)\right).
\]
\end{theorem}

\begin{theorem}[\cite{A21,A45}]\label{thm:hilbertprops}
The Hermitian Hilbert transform $\Hilb$ enjoys the following
properties.
\begin{enumerate}[label=\upshape(\roman*)]
\item $\Hilb$ is a bounded linear operator on $C^{0,\nu}(\Gamma)$,
$0<\nu<1$: there is $c_{\Gamma}>0$, depending only on $\Gamma$, with
$\|\Hilb[\Gmat]\|_{\nu,\Gamma}\le c_{\Gamma}\|\Gmat\|_{\nu,\Gamma}$.
\item $\Hilb$ is an involution on $C^{0,\nu}(\Gamma)$, i.e.\
$\Hilb^{2}=\mathbf{1}$.
\end{enumerate}
\end{theorem}

Both statements admit an $L^{p}$ counterpart \cite{A45}. If $\Gamma$ is a Lipschitz hypersurface and
$\Gmat\in L^{p}(\Gamma)$, $1<p<\infty$, then the Cauchy integral is
understood in the Lebesgue sense, the limits in Theorem~\ref{thm:PS}
exist almost everywhere, and $\Hilb$ is bounded on $L^{p}(\Gamma)$ by
Calder\'on--Zygmund theory together with the density of
$C^{0,\nu}(\Gamma)$ in $L^{p}(\Gamma)$.

\subsection{Hardy projections and trace characterizations}
\label{subsec:hardy}

Theorem~\ref{thm:hilbertprops} immediately yields the Hardy-type
projections
\begin{equation}\label{eq:projections}
P^{+}_{\Gamma}=\tfrac12\bigl(\mathbf{1}+\Hilb\bigr),
\qquad
P^{-}_{\Gamma}=\tfrac12\bigl(\mathbf{1}-\Hilb\bigr).
\end{equation}

\begin{theorem}[\cite{A21}]\label{thm:projections}
$P^{+}_{\Gamma}$ and $P^{-}_{\Gamma}$ are mutually complementary
projection operators:
\[
\bigl(P^{\pm}_{\Gamma}\bigr)^{2}=P^{\pm}_{\Gamma},\qquad
P^{+}_{\Gamma}+P^{-}_{\Gamma}=\mathbf{1},\qquad
P^{+}_{\Gamma}P^{-}_{\Gamma}=P^{-}_{\Gamma}P^{+}_{\Gamma}=\Zero .
\]
Consequently one has the direct decomposition
\[
C^{0,\nu}(\Gamma)=P^{+}_{\Gamma}\bigl[C^{0,\nu}(\Gamma)\bigr]
\oplus P^{-}_{\Gamma}\bigl[C^{0,\nu}(\Gamma)\bigr],
\qquad 0<\nu<1,
\]
so that each $\Gmat\in C^{0,\nu}(\Gamma)$ decomposes uniquely into
components belonging to the two summands.
\end{theorem}

The projections characterize exactly which H\"older continuous matrix
functions on $\Gamma$ are traces of $H$-monogenic functions.

\begin{theorem}[\cite{A21,A47}]\label{thm:traces}
Let $\Gmat\in C^{0,\nu}(\Gamma)$, $0<\nu<1$.
\begin{enumerate}[label=\upshape(\roman*)]
\item $\Gmat$ is the boundary value of a matrix function which is
$H$-monogenic in $\Omega^{+}$ if and only if
$\Gmat\in\operatorname{im}P^{+}_{\Gamma}$, equivalently
$\Hilb[\Gmat]=\Gmat$.
\item $\Gmat$ is the boundary value of a matrix function which is
$H$-monogenic in $\Omega^{-}$ and vanishes at infinity if and only if
$\Gmat\in\operatorname{im}P^{-}_{\Gamma}$, equivalently
$\Hilb[\Gmat]=-\Gmat$.
\end{enumerate}
Analogous statements hold for right $H$-monogenic functions.
\end{theorem}

\subsection{Two-sided \texorpdfstring{$H$}{H}-monogenicity}

The Hilbert transform also detects the two-sidedness of
$H$-monogenicity, a phenomenon with no counterpart in the scalar
Euclidean theory.

\begin{theorem}[\cite{A21}]\label{thm:twosided1}
Let $\Gmat\in C^{0,\nu}(\overline{\Omega})$, $0<\nu<1$, with
$\Dop[\Gmat]=\Zero$ in $\Omega$.  Then the following are equivalent:
\begin{enumerate}[label=\upshape(\roman*)]
\item $\Gmat$ is two-sided $H$-monogenic in $\Omega$;
\item $\Hilb\bigl[\Gmat|_{\Gamma}\bigr]
=\bigl[\Gmat|_{\Gamma}\bigr]\Hilb$.
\end{enumerate}
\end{theorem}

\begin{theorem}[\cite{A21}]\label{thm:twosided2}
Let $\Gmat:\overline{\Omega}\to\CM$ be a continuous matrix function
whose trace $g^{1}_{2}=\Gmat|_{\Gamma}$ belongs to $C^{0,\nu}(\Gamma)$.  Then
$\Gmat$ is two-sided $H$-monogenic in $\Omega$ if and only if it is
harmonic in $\Omega$ and
\[
\Hilb\bigl[\Gmat|_{\Gamma}\bigr]=\Gmat|_{\Gamma}
=\bigl[\Gmat|_{\Gamma}\bigr]\Hilb .
\]
\end{theorem}

\begin{remark}
Theorems~\ref{thm:twosided1} and \ref{thm:twosided2} are the
$(\ell)$AD-regular ancestors of Theorems~\ref{thm:fractalcrit}--\ref{thm:fractalconservation}: in both
geometries $H$-monogenicity is equivalent to harmonicity plus the
identity $\Hilb[g^{1}_{2}]=g^{1}_{2}$, and two-sidedness to the
commutation $\Hilb[g^{1}_{2}]=g^{1}_{2}\Hilb$, only the meaning of the
operator being adjusted to the regularity of $\Gamma$.
\end{remark}

\subsection{The Hilbert transform on fractal boundaries}
\label{subsec:fractalhilbert}

None of the ingredients of Subsections~\ref{subsec:PS}--%
\ref{subsec:hardy} survives verbatim when $\dimH(\Gamma)>2n-1$: the
singular integrals $S_{\Gamma}$, $S^{|}_{\Gamma}$ have no meaning and
the non-tangential limit procedure loses its support.  The way out,
proposed in \cite{ABBM2012}, is again to take a known formula as a
definition: one declares the Hilbert transform to be the trace on
$\Gamma$ of the continuous extension of the fractal Cauchy integral of
Definition~\ref{def:fractalcauchy}, i.e.
\begin{equation}\label{eq:fractalhilbertdef}
\Hilb^{*}\bigl[\Gmat\bigr](U)
=\frac{2}{\lam}\bigl(\Cauchy^{*}\bigl[\Gmat\bigr]\bigr)^{+}(U)
-\Gmat(U),\qquad U\in\Gamma,
\end{equation}
whenever that continuous extension exists.  By the Borel--Pompeiu
formula this is equivalent to
\begin{equation}\label{eq:fractalhilbert}
\Hilb^{*}\bigl[\Gmat\bigr](U)
=\Gmat(U)-\frac{2}{\lam}\,
\Teo\bigl[\Dop\varUpsilon[{G}^{1}_{2}]\bigr](U),\qquad U\in\Gamma,
\end{equation}
which is a natural generalization of the more conventional expression
\[
\Hilb\bigl[\Gmat\bigr](U)=\Gmat(U)-\frac{2}{\lam}\,\int_{\Gamma}
\Ekernel(Z-V)\,N_{(Z,\Zd)}\bigl(\Gmat(U)-\Gmat(X)\bigr)
\,d\mathcal{H}^{2n-1}.
\]
\begin{theorem}[\cite{ABBM2012,A49}]\label{thm:fractalhilbertwell}
Let $\Omega$ be a Jordan domain in $\R^{2n}$ with $d$-summable boundary
$\Gamma$, $d\in(2n-1,2n)$, let $d-2n+1<\nu<1$ and
$\Gmat\in C^{0,\nu}(\Gamma)$.  If moreover
\begin{equation}\label{eq:alphacond}
\nu>\frac{d}{2n},
\end{equation}
then $\Cauchy^{*}[\Gmat]$ admits a continuous extension to
$\overline{\Omega}$, so that $\Hilb^{*}$ is well defined by
\eqref{eq:fractalhilbertdef}, and
$\Hilb^{*}[\Gmat]\in C^{0,\beta}(\Gamma)$ whenever
\begin{equation}\label{eq:betacond}
0<\beta<\frac{2n\nu-d}{2n-d}<1 .
\end{equation}
\end{theorem}

\begin{remark}[\cite{ABBM2012}]\label{rem:righthanded}
Right-handed versions of all the operators above are obtained from the
alternative definitions
\[
\bigl(\Gmat\,\Cauchy^{*}\bigr)(Y)
=\lam\,\boldsymbol{\chi}_{\Omega}(Y)\varUpsilon[{G}^{1}_{2}](Y)
-\bigl(\varUpsilon[{G}^{1}_{2}]\,\Dop\,\Teo\bigr)(Y),
\]
\[
\bigl(\Gmat\,\Hilb^{*}\bigr)(U)
=\Gmat(U)-\frac{2}{\lam}
\bigl(\varUpsilon[{G}^{1}_{2}]\,\Dop\,\Teo\bigr)(U),
\]
where $(\Fmat\,\Teo)(Y)
=-\int_{\Omega}\Fmat(X)\,\Ekernel(Z-V)\,dW(Z,\Zd)$, and they enjoy the
same properties.  All the criteria of
Subsection~\ref{subsec:fractalcriteria} come in a left- and a
right-handed version accordingly.
\end{remark}

The main quantitative result of \cite{A49} shows that $\Hilb^*$ is in fact bounded and makes the dependence on the geometry of $\Gamma$ explicit through
the $d$-sum $s(d)$ of Lemma~\ref{lem:dsum}.

\begin{theorem}[H\"older norm estimate, $d$-summable case, \cite{A49}]
\label{thm:normestimate}
Let $\Gamma$ be $d$-summable and let \eqref{eq:betacond} hold.  Then
$\Hilb^{*}$ is bounded from $C^{0,\nu}(\Gamma)$ into
$C^{0,\beta}(\Gamma)$ and
\[
\bigl\|\Hilb^{*}\bigr\|
\;\le\;
1+|\Gamma|^{\nu-\beta}
+c_{1}\bigl(s(\tilde d)\bigr)^{\frac{1-\beta}{2n}}
+c_{2}\bigl(s(\tilde d)\bigr)^{\frac{1-\beta}{2n}}|\Gamma|^{\beta},
\]
where $|\Gamma|=\operatorname{diam}\Gamma$,
$\tilde d=2n\frac{\nu-\beta}{1-\beta}$, and $c_{1},c_{2}$ depend
only on $\nu,\beta$ and $n$.
\end{theorem}

Since, by Lemma~\ref{lem:ADdsum}, $(\ell)$AD-regularity implies
$d$-summability for every $d>2n-1$, the previous theorem contains the
classical situation, in which the estimate can be phrased purely in
terms of the $(2n-1)$-dimensional Hausdorff measure of $\Gamma$.

\begin{theorem}[H\"older norm estimate, $(\ell)$AD-regular case,
\cite{A49}]\label{thm:normestimateAD}
Let $\Gamma$ be $(\ell)$AD-regular.  Then $\Hilb$ is bounded from
$C^{0,\nu}(\Gamma)$ into $C^{0,\beta}(\Gamma)$ whenever
$0<\beta<2n\nu+1-2n<1$, and
\[
\|\Hilb\|\;\le\;
1+c_{1}\bigl[\mathcal{H}^{2n-1}(\Gamma)\bigr]^{\frac{\nu-\beta}{2n-1}}
+c_{2}\bigl[\mathcal{H}^{2n-1}(\Gamma)\bigr]^{\frac{1-\beta}{2n}}
+c_{3}\bigl[\mathcal{H}^{2n-1}(\Gamma)\bigr]^{
\frac{2n-1+\beta}{2n(2n-1)}},
\]
where $c_{1},c_{2},c_{3}$ depend only on $\nu,\beta,n$ and on the
regularity constant $c_{\Gamma}$ of \eqref{eq:AD}.
\end{theorem}

\subsection{Criteria for \texorpdfstring{$H$}{H}-monogenicity on fractal
domains}\label{subsec:fractalcriteria}

The fractal Hilbert transform is not merely a formal substitute for
$\Hilb$: it detects $H$-monogenicity, exactly as its smooth ancestor
does through Theorem~\ref{thm:traces}.  The following results,
established in \cite{ABBM2012}, are the fractal counterparts of
Theorem~\ref{thm:traces}, \ref{thm:twosided1} and
\ref{thm:twosided2}. Throughout, $\Omega$ is a Jordan domain with
$d$-summable boundary $\Gamma$, $d\in(2n-1,2n)$, and
$d-2n+1<\nu\le1$ is such that \eqref{eq:alphacond} holds.

We first record a reproducing property which is of independent
interest, since it relaxes the regularity assumed in
Theorem~\ref{thm:fractalcauchyformula}.

\begin{theorem}[\cite{ABBM2012}]\label{thm:fractalreproduce}
Let $\Gmat\in C^{0,\nu}(\overline{\Omega})$ be $H$-monogenic in
$\Omega$ and let $g^{1}_{2}=\Gmat|_{\Gamma}$ be its trace.  Then
\[
\Cauchy^{*}\bigl[g^{1}_{2}\bigr](Y)=\lam\,\Gmat(Y),
\qquad Y\in\Omega .
\]
\end{theorem}

\begin{remark}[\cite{ABBM2012}]
Theorem~\ref{thm:fractalreproduce} is a small but useful improvement of
Theorem~\ref{thm:fractalcauchyformula}, where the hypothesis
$\Gmat\in C^{1}(\overline{\Omega})$ was made a priori; only H\"older
continuity up to the boundary is in fact needed.
\end{remark}

\begin{theorem}[\cite{ABBM2012}]\label{thm:fractalcrit}
Let $\Gmat\in C^{0,\nu}(\overline{\Omega})$ with trace
$g^{1}_{2}=\Gmat|_{\Gamma}$.  Then $\Gmat$ is left, resp.\ right,
$H$-monogenic in $\Omega$ if and only if $\Gmat$ is harmonic in
$\Omega$ and
\[
\Hilb^{*}\bigl[g^{1}_{2}\bigr]=g^{1}_{2},
\qquad\text{resp.}\qquad
g^{1}_{2}\,\Hilb^{*}=g^{1}_{2}.
\]
\end{theorem}

\begin{theorem}[Conservation law on a fractal boundary,
\cite{ABBM2012}]\label{thm:fractalconservation}
Let $\Gmat\in C^{0,\nu}(\overline{\Omega})$ be left $H$-monogenic in
$\Omega$, with trace $g^{1}_{2}$.  Then the following are equivalent:
\begin{enumerate}[label=\upshape(\roman*)]
\item $\Gmat$ is two-sided $H$-monogenic in $\Omega$;
\item $\Hilb^{*}\bigl[g^{1}_{2}\bigr]=g^{1}_{2}\,\Hilb^{*}$.
\end{enumerate}
\end{theorem}

The pattern is worth underlining: harmonicity, which is the
``unconditional'' part of $H$-monogenicity, is separated from the
boundary condition, and the whole of the boundary information is
carried by a single identity involving $\Hilb^{*}$; two-sidedness, in
turn, corresponds exactly to the commutation of the trace with the
Hilbert transform.  Neither statement makes any reference to a surface
integral over $\Gamma$, which is what allows them to survive in the
fractal regime.

% ==================================================================
\section{Jump and decomposition problems}\label{sec:jump}
% ==================================================================

\subsection{The scalar jump problem for \texorpdfstring{$h$}{h}-monogenic
functions}

Historically, the first boundary value problem attacked in the Hermitian
setting was the following \cite{A19}.

\begin{problem}[Jump problem]\label{prob:jump}
Given $f\in C^{0,\nu}(\Gamma)$, $0<\nu<1$, decompose
\begin{equation}\label{eq:jump}
f=f^{+}-f^{-}\quad\text{on }\Gamma,
\end{equation}
where $f^{\pm}\in C^{0,\nu}(\Gamma)$ are extendable to
$h$-monogenic functions $F^{\pm}$ in $\Omega^{\pm}$ with
$F^{-}(\infty)=0$.
\end{problem}

If a solution exists it is unique, as follows from the Painlev\'e and
Liouville theorems in Clifford analysis.  Two independent criteria were
established in \cite{A19}, the first one integral, the second one
differential.

\begin{theorem}[Integral criterion for $h$-monogenicity, \cite{A19}]
\label{thm:integralcriterion}
Let $\Gamma$ be $(\ell)$AD-regular.  A function
$f\in C^{0,\nu}(\Gamma)$ has an $h$-monogenic extension $F^{\pm}$
to $\Omega^{\pm}$, with $F^{-}(\infty)=0$, if and only if
\[
S_{\Gamma}f=\pm f=S^{|}_{\Gamma}f .
\]
\end{theorem}

\begin{theorem}[\cite{A19}]\label{thm:jumpequiv}
Under the hypotheses of Theorem~\ref{thm:integralcriterion}, the following statements are equivalent:
\begin{enumerate}[label=\upshape(\roman*)]
\item $f$ admits the decomposition \eqref{eq:jump};
\item $S_{\Gamma}f=S^{|}_{\Gamma}f$;
\item ${C}_{\Gamma}f=\mathcal{C}^{|}_{\Gamma}f$;
\item ${C}_{\Gamma}f$ is $h$-monogenic in
$\R^{2n}\setminus\Gamma$.
\end{enumerate}
Moreover, if \eqref{eq:jump} is solvable, its unique solution is
\[
f^{\pm}={C}^{\pm}_{\Gamma}f
=\tfrac12\bigl(S_{\Gamma}f\pm f\bigr)
={C}^{|,\pm}_{\Gamma}f
=\tfrac12\bigl(S^{|}_{\Gamma}f\pm f\bigr).
\]
\end{theorem}

The second approach replaces the integral criterion by a
\emph{conservation law} on $\Gamma$.  Assume $\Gamma$ to be
$C^{1}$-smooth, let $\eta$ denote the unit normal and decompose
\[
\dX=-\eta\bigl(\eta\dX\bigr)=\eta\,\partial_{\eta}+\partial^{\Gamma}_{X},
\qquad
\dXt=-\eta^{|}\bigl(\eta^{|}\dXt\bigr)
=\eta^{|}\partial_{\eta}+\partial^{\Gamma}_{\Xt},
\]
where $\partial_{\eta}=\langle\eta,\dX\rangle$ and
$\partial^{\Gamma}_{X}=-\eta(\eta\wedge\dX)$; the restrictions of
$\partial^{\Gamma}_{X}$ and $\partial^{\Gamma}_{\Xt}$ to $\Gamma$ are
denoted $\partial_{\omega}$ and $\partial_{\omega^{|}}$.

\begin{theorem}[Conservation law, \cite{A19}]\label{thm:conservation}
Let $F^{\pm}\in C^{1}(\overline{\Omega^{\pm}})$ be monogenic in
$\Omega^{\pm}$ with respect to $\dX$, with $F^{-}(\infty)=0$.  Then
$F^{\pm}$ is $h$-monogenic in $\Omega^{\pm}$ if and only if
$g=F^{\pm}|_{\Gamma}$ satisfies
\begin{equation}\label{eq:conservation}
\eta^{|}\eta\,\partial_{\omega}g+\partial_{\omega^{|}}g=0
\qquad\text{on }\Gamma .
\end{equation}
\end{theorem}

\begin{theorem}[\cite{A19}]\label{thm:jumpconservation}
The jump problem \eqref{eq:jump} is solvable if and only if
\[
\eta^{|}\eta\,\partial_{\omega}f+\partial_{\omega^{|}}f=0
\qquad\text{on }\Gamma .
\]
\end{theorem}

Theorem~\ref{thm:jumpconservation} is remarkable in that it converts an
integral solvability condition into a first-order differential
condition on the data, intrinsic to the hypersurface.

\subsection{Matrix decomposition on Ahlfors--David regular hypersurfaces}

The matrix version of Problem~\ref{prob:jump} reads: given
$\Gmat\in C^{0,\nu}(\Gamma)$, decompose
\begin{equation}\label{eq:matrixjump}
\Gmat=\bigl(\Gmat\bigr)^{+}-\bigl(\Gmat\bigr)^{-},
\end{equation}
where the components are extendable to $H$-monogenic matrix functions
in $\Omega^{\pm}$, the one in $\Omega^-$ vanishing at infinity.  Uniqueness again follows
from the Painlev\'e and Liouville theorems.  Theorem~\ref{thm:PS}
immediately gives existence and the explicit formula
\begin{equation}\label{eq:matrixjumpsol}
\bigl(\Gmat\bigr)^{\pm}=\frac{1}{\lam}\,
\Cauchy^{\pm}\bigl[\Gmat\bigr],
\end{equation}
so that in the matrix setting the jump problem is \emph{always}
solvable.  The subtlety lies entirely in returning to the $h$-monogenic
picture, which requires the Cauchy integral of a diagonal matrix to be
again diagonal.

\begin{theorem}[\cite{A21}]\label{thm:diagjump}
Let $g\in C^{0,\nu}(\Gamma)$, $0<\nu<1$, and let $G_{0}$ be the
associated diagonal matrix function.  Then the jump problem
\eqref{eq:matrixjump} is solvable in terms of $h$-monogenic functions
if and only if
\[
\Cauchy\bigl[G_{0}\bigr]
=\frac{\lam}{2}
\begin{pmatrix}
{C}_{\Gamma}g+{C}^{|}_{\Gamma}g & 0\\[1mm]
0 & {C}_{\Gamma}g+{C}^{|}_{\Gamma}g
\end{pmatrix}.
\]
Equivalently, in terms of the singular operators,
$S_{\Gamma}g=S^{|}_{\Gamma}g$.
\end{theorem}

\begin{theorem}[\cite{A21}]\label{thm:twosidedjump}
Let $\Gmat\in C^{0,\nu}(\Gamma)$, $0<\nu<1$.  The following are
equivalent:
\begin{enumerate}[label=\upshape(\roman*)]
\item $\Gmat$ admits a decomposition \eqref{eq:matrixjump} whose
components are extendable to \emph{two-sided} $H$-monogenic functions
in $\Omega^{\pm}$, the one in $\Omega^-$ vanishing at infinity;
\item $\Cauchy[\Gmat]$ is two-sided $H$-monogenic in
$\R^{2n}\setminus\Gamma$.
\end{enumerate}
\end{theorem}

\subsection{Decomposition on fractal boundaries}

On fractal boundaries the Cauchy integral is unavailable in its
classical form, and the decomposition must be produced directly by the
Teodorescu transform.  The scalar prototype, due to  \cite{A25} and formulated 
in the language of the Teodorescu transform, is the following.
Under
\begin{equation}\label{eq:nuboxcond}
1>\nu>\frac{\alpha(\Gamma)}{2n}
\end{equation}
the function $\dX \varUpsilon[f]$ is integrable in $\Omega$ with any exponent
less than $(2n-\alpha(\Gamma))/(1-\nu)$; a decomposition
$f=\Phi^{+}-\Phi^{-}$ is said to be of \emph{class $C^{0,\mu}$} when
$\Phi$ is $\mu$-H\"older continuous on $\overline{\Omega}$ and on
$\R^{2n}\setminus\Omega$ and
\[
\mu<\frac{2n\nu-\alpha(\Gamma)}{2n-\alpha(\Gamma)},
\qquad
\mu>\dimH(\Gamma)-2n+1 .
\]

\begin{theorem}[Monogenic decomposition on fractal surfaces,
\cite{A25}]\label{thm:monodecomp}
Let $\alpha(\Gamma)<2n$ and assume \eqref{eq:nuboxcond}.  Then, for
every $f\in C^{0,\nu}(\Gamma)$, there is a unique
monogenic decomposition of class $C^{0,\mu}$, given by
\[
\Phi(X)=\chi_{\Omega}(X)\,\varUpsilon[f](X)+T_{\Omega}\bigl[\dX \varUpsilon[f]\bigr](X).
\]
\end{theorem}

\begin{theorem}[Hermitian decomposition on a fractal surface,
\cite{A25}]\label{thm:hermdecomp}
Let $\alpha(\Gamma)<2n$, assume \eqref{eq:nuboxcond} and let
$f\in C^{0,\nu}(\Gamma)$.  Then there exists a unique
decomposition $f=\Phi^{+}-\Phi^{-}$ of class $C^{0,\mu}$ with
$\Phi^{\pm}$ \emph{Hermitian} monogenically extendable to $\Omega$ and
to $\R^{2n}\setminus\overline{\Omega}$ respectively, if and only if
\begin{equation}\label{eq:hermdecompcond}
T_{\Omega}\bigl[\dX \varUpsilon[f]\bigr]\Big|_{\Gamma}
=T^{|}_{\Omega}\bigl[\dXt \varUpsilon[f]\bigr]\Big|_{\Gamma},
\end{equation}
where $T^{|}_{\Omega}$ is the Teodorescu transform associated with
the fundamental solution $E^{|}$ of $\dXt$.
\end{theorem}

\begin{corollary}[\cite{A25}]
If $\Gamma$ is $(2n-1)$-rectifiable, Theorem~\ref{thm:hermdecomp}
remains true under the weaker assumption
$1>\nu>\frac{2n-1}{2n}$, with $0<\mu<2n(\nu-1)+1$.
\end{corollary}

The matricial version, which is at the same time simpler and more
powerful, is the content of \cite{A36}. Throughout, assume
\begin{equation}\label{eq:nucond36}
\nu>\frac{\alpha(\Gamma)}{2n},
\end{equation}
and let $\mu$ satisfy
$\dimH(\Gamma)-2n+1<\mu<\frac{2n\nu-\alpha(\Gamma)}{2n-\alpha(\Gamma)}$.

\begin{theorem}[Hermitian Teodorescu decomposition, \cite{A36}]
\label{thm:teodecomp}
Let $\Gmat\in C^{0,\nu}(\Gamma)$ with $\nu$ as in \eqref{eq:nucond36}.
Then $\Gmat$ admits the $H$-monogenic decomposition
\eqref{eq:matrixjump}, the components being explicitly given by
\[
\bigl(\Gmat\bigr)^{+}
=\varUpsilon[{G}^{1}_{2}]-\frac{1}{\lam}\,
\Teo\bigl[\Dop\varUpsilon[{G}^{1}_{2}]\bigr],
\qquad
\bigl(\Gmat\bigr)^{-}
=-\frac{1}{\lam}\,
\Teo\bigl[\Dop\varUpsilon[{G}^{1}_{2}]\bigr].
\]
\end{theorem}

\begin{theorem}[\cite{A36}]\label{thm:teodecompdiag}
Let $g\in C^{0,\nu}(\Gamma)$ with $\nu$ as in \eqref{eq:nucond36}, and
let $G_{0}$ be the corresponding diagonal matrix function.  Then
$G_{0}$ admits the decomposition \eqref{eq:matrixjump} in terms of
$h$-monogenic functions if and only if
\[
\frac{1}{\lam}\,\Teo\bigl[\Dop\varUpsilon[{G}_{0}] \bigr]
\quad\text{is a diagonal matrix function,}
\]
that is, if and only if
\[
T_{\Omega}\bigl[\dZ\varUpsilon[{g}]-\dZd\varUpsilon[{g}]\bigr]
=iT^{|}_{\Omega}\bigl[\dZ\varUpsilon[{g}]+\dZd[{g}]\bigr],
\]
which, after a direct computation, can be rewritten as
\[
T_{\Omega}\bigl[\dX\varUpsilon[{g}]\bigr]
=T^{|}_{\Omega}\bigl[\dXt\varUpsilon[{g}]\bigr].
\]
\end{theorem}

\begin{remark}[\cite{A36}]
Although the decomposition \eqref{eq:matrixjump} is not unique in
general, it \emph{is} unique in the class
$C^{0,\mu}=\{\Gmat\mid g_{1},g_{2}\in C^{0,\mu}\}$ for
$\dimH(\Gamma)-2n+1<\mu<\frac{2n\nu-\alpha(\Gamma)}{2n-\alpha(\Gamma)}$; 
this rests on a removability theorem of Dolzhenko type, namely: if
$\dimH(K)<2n+\mu-1$, then every $g\in C^{0,\mu}(\R^{2n})$ which is
monogenic in $\R^{2n}\setminus K$ is monogenic in the whole of
$\R^{2n}$.  When $\Gamma$ is sufficiently regular the formulae of
Theorem~\ref{thm:teodecomp} reduce to \eqref{eq:matrixjumpsol}.
\end{remark}

% ==================================================================
\section{Dirichlet-type problems}\label{sec:dirichlet}
% ==================================================================

\subsection{Dirichlet problems in the Sobolev scale}

In \cite{A59} two Dirichlet-type problems are studied in the Hermitian
Hilbert--Sobolev scale.  They read
\begin{align}
\Dop\,\Xmat&=\Fmat \quad\text{in }\Omega, &
\Xmat&=\Gmat\quad\text{on }\Gamma,
\label{eq:BVP1}\\[1mm]
\Dop^{2}\,\Xmat&=\Fmat\quad\text{in }\Omega, &
\Xmat&=\Gmat\quad\text{on }\Gamma .
\label{eq:BVP2}
\end{align}
In this section $\mathcal C_\Gamma$ and $\mathcal T_\Omega$ stand for $\lam^{-1}$ times the operators (\ref{eq:cauchyint}) and (\ref{eq:teo}) and have the expected mapping properties, namely
$\Dop:W^{k}_{2}(\Omega)\to W^{k-1}_{2}(\Omega)$,
$\Teo:W^{k}_{2}(\Omega)\to W^{k+1}_{2}(\Omega)$ and
$\Cauchy:W^{k-\frac12}_{2}(\Gamma)\to W^{k}_{2}(\Omega)$, and
Theorems~\ref{thm:BP}--\ref{thm:teoinverse} together with the
Plemelj--Sokhotski formulae remain valid in this scale (with $\Teo$
inverting $\Dop$ already on $L^{2}(\Omega)$).  The key structural
result is an orthogonal splitting of the Sobolev spaces; writing
$\mathring{W}^{1}_{2}(\Omega)$ for the subspace of $W^{1}_{2}(\Omega)$
of functions with vanishing trace on $\Gamma$.

\begin{proposition}[Orthogonal decomposition, \cite{A59}]
\label{prop:orthdecomp}
With respect to the inner product \eqref{eq:innerprod},
\[
W^{k}_{2}(\Omega)
=\bigl(W^{k}_{2}(\Omega)\cap\ker\Dop\bigr)
\oplus
\Dop\bigl(\mathring{W}^{1}_{2}(\Omega)\cap W^{k+1}_{2}(\Omega)\bigr).
\]
We denote by $P$ and $Q=\mathbf{1}-P$ the corresponding orthoprojections.
\end{proposition}

\begin{proposition}[\cite{A59}]\label{prop:trTeo}
Let $\Fmat\in L^{2}(\Omega)$.  Then $\tr_{\Gamma}\Teo[\Fmat]=0$ if and
only if $\Fmat\in\operatorname{im}Q$.
\end{proposition}

\begin{theorem}[Solvability of \eqref{eq:BVP1}, \cite{A59}]
\label{thm:BVP1}
Let $\Fmat\in W^{k}_{2}(\Omega)$ and
$\Gmat\in W^{k+\frac12}_{2}(\Gamma)$.  Then \eqref{eq:BVP1} has a
solution if and only if
\begin{equation}\label{eq:BVP1cond}
P^{-}_{\Gamma}\bigl[\Gmat\bigr]=\tr_{\Gamma}\Teo\bigl[\Fmat\bigr],
\end{equation}
and in that case the unique solution $\Xmat\in W^{k+1}_{2}(\Omega)$ is
\[
\Xmat=\Cauchy\bigl[\Gmat\bigr]+\Teo\bigl[\Fmat\bigr].
\]
\end{theorem}

\begin{theorem}[Solvability of \eqref{eq:BVP2}, \cite{A59}]
\label{thm:BVP2}
Let $\Fmat\in W^{k}_{2}(\Omega)$ and
$\Gmat\in W^{k+\frac32}_{2}(\Gamma)$, $k\ge 0$.  Then \eqref{eq:BVP2}
has the unique solution $\Xmat\in W^{k+2}_{2}(\Omega)$ given by
\[
\Xmat=\Cauchy\bigl[\Gmat\bigr]
+\Teo\,P\bigl[\Dop\Hmat\bigr]
+\Teo\,Q\,\Teo\bigl[\Fmat\bigr],
\]
where $\Hmat$ is an appropriate $W^{k+2}_{2}$-extension of the
boundary data associated with the auxiliary homogeneous problem.
\end{theorem}

Finally, the orthoprojections of Proposition~\ref{prop:orthdecomp}
admit a closed algebraic representation in terms of the boundary
operators, a result of independent interest.

\begin{proposition}[\cite{A59}]\label{prop:iso}
For every $k\in\N\cup\{0\}$ the operator
\[
\tr_{\Gamma}\Teo\,\Cauchy:\;
W^{k+\frac12}_{2}(\Gamma)\cap\operatorname{im}P^{+}_{\Gamma}
\longrightarrow
W^{k+\frac32}_{2}(\Gamma)\cap\operatorname{im}P^{-}_{\Gamma}
\]
is an isomorphism.
\end{proposition}

\begin{theorem}[\cite{A59}]\label{thm:projrep}
Let $\Xmat\in W^{k}_{2}(\Omega)$, $k\ge 1$.  Then
\[
P\bigl[\Xmat\bigr]
=\Cauchy\bigl(\tr_{\Gamma}\Teo\,\Cauchy\bigr)^{-1}
\tr_{\Gamma}\Teo\bigl[\Xmat\bigr],
\qquad
Q\bigl[\Xmat\bigr]
=\Bigl(\mathbf{1}-\Cauchy\bigl(\tr_{\Gamma}\Teo\,\Cauchy\bigr)^{-1}
\tr_{\Gamma}\Teo\Bigr)\bigl[\Xmat\bigr],
\]
both belonging to $W^{k}_{2}(\Omega)$.
\end{theorem}

% ==================================================================
\section{Duality: from harmonic forms to Hermitian systems}
\label{sec:duality}
% ==================================================================

The circulant matrix formalism has a further, more algebraic, payoff:
it allows a clean description of the topological dual of the space of
germs of $H$-monogenic functions, in the spirit of the classical
K\"othe--Grothendieck duality for holomorphic functions
\cite{Kothe1953}.  The theory has three layers, and
we present them in the order in which they clarify one another.  The
first is the Euclidean model \cite{ABDS2011}: there the
dual of the module of germs of monogenic functions is represented by a
single monogenic function outside the compact set, whereas the dual of
the space of harmonic multivector fields or, equivalently, of
harmonic differential forms, is represented by a \emph{pair}
subject to a Moisil--Teodorescu-type system, the representation
being merely surjective with an explicitly described kernel.  The
second layer is the Hermitian matricial one \cite{A47,CSS2012}, where
the duality is again a clean isomorphism, algebraic as well as
topological.  The third is the Hermitian scalar one
\cite{A47,CSS2012}, where the representing object is once more a pair,
now subject to a system built out of $\dZ$ and $\dZd$: the Euclidean
pattern repeats itself verbatim one refinement higher.

\subsection{The Euclidean model: monogenic germs, harmonic multivector
fields and their annihilators}\label{subsec:euclideandual}

Throughout this subsection we work, following
\cite{ABDS2011}, in the real Clifford algebra
$\R_{0,n}$ over $\R^{n}$ and write $\mathcal{M}_{L}(\Omega)$ and
$\mathcal{M}_{R}(\Omega)$ for the spaces of left, resp.\ right,
monogenic functions in an open set $\Omega$, and
$\mathcal{M}_{L}(\Omega,0)$, $\mathcal{M}_{R}(\Omega,0)$ the
subspaces of functions vanishing at infinity when $\Omega$ is
unbounded. Equipped with the topology of uniform convergence on
compact sets, $\mathcal{M}_{L}(\Omega)$ is a right Fr\'echet
$\R_{0,n}$-module.  For $\mathcal{K}\subset\R^{n}$ compact with
$\R^{n}\setminus \mathcal{K}$ connected one forms, exactly as in
Subsection~\ref{subsec:matrixdual} below, the module
$\mathcal{M}_{L}(\mathcal{K})$ of germs.

\begin{theorem}[\cite{ABDS2011}]\label{thm:eucdual}
Let $\Lambda$ be a continuous right $\R_{0,n}$-linear functional on
$\mathcal{M}_{L}(\mathcal{K})$.  Then there exists a unique
$\omega\in\mathcal{M}_{R}(\R^{n}\setminus \mathcal{K})$ with $\omega(\infty)=0$
such that, for every germ $\dot u\in\mathcal{M}_{L}(\mathcal{K})$,
\[\Lambda[\dot u]=\int_{\Gamma_{u}}\omega(y)\,\eta(y)\,u(y)\, d\mathcal{H}^{n-1}(y),\]
$u$ being any representative of the germ and $\Gamma_{u}$ a smooth
surface surrounding $\mathcal{K}$.  In other words, $\mathcal{M}_{L}(\mathcal{K})^{*}$ and
$\mathcal{M}_{R}(\R^{n}\setminus \mathcal{K},0)$ are isomorphic.
\end{theorem}

The situation changes when one replaces monogenic functions by
harmonic multivector fields.  Recall that for a $k$-vector-valued
function $F_{k}$ the action of $\partial_{X}$ splits as
$\partial_{X}F_{k}=\partial_{X}\bullet F_{k}+\partial_{X}\wedge F_{k}$ into a $(k-1)$- and a
$(k+1)$-vector part, so that $F_{k}$ is monogenic if and only if
\begin{equation}\label{eq:hodgeclifford}
\partial_{X}\bullet F_{k}=0,\qquad \partial_{X}\wedge F_{k}=0 .
\end{equation}
A $k$-vector field satisfying \eqref{eq:hodgeclifford} in $\Omega$ is
called a \emph{harmonic $k$-grade multivector field}, and the space of
these is denoted $\mathcal{HF}_{k}(\Omega)$.  Identifying $F_{k}$ with
the $k$-form $\omega_{k}$ having the same components, one checks that
\eqref{eq:hodgeclifford} is precisely the Hodge--de Rham system
$d\omega_{k}=0$, $d^{*}\omega_{k}=0$, so that harmonic multivector
fields and harmonic differential forms are in one-to-one
correspondence \cite{ABDS2011}; for $k=1$ one recovers the Riesz
system, the simplest of the Stein--Weiss systems of conjugate harmonic
functions.

Dually, one considers pairs $(F_{k-1},F_{k+1})$ satisfying in $\Omega$
the Moisil--Teodorescu-type system
\begin{equation}\label{eq:MTsystem}
\partial_{X}\wedge F_{k-1}+\partial_{X}\bullet F_{k+1}=0,
\qquad
\partial_{X}\bullet F_{k-1}=0,\qquad \partial_{X}\wedge F_{k+1}=0,
\end{equation}
equivalently: pairs for which $f=F_{k-1}+F_{k+1}$ is left monogenic.
Such a pair is called a \emph{left harmonic pair}; right harmonic pairs
are defined by the corresponding right-hand system, and the space of
right harmonic pairs in $\R^{n}\setminus \mathcal{K}$ vanishing at infinity is
denoted $\mathcal{HP}^{r}_{k}(\R^{n}\setminus \mathcal{K},0)$.  In the language
of forms, \eqref{eq:MTsystem} is the system
\[
dv_{k-1}-d^{*}v_{k+1}=0,\qquad d^{*}v_{k-1}=0,\qquad dv_{k+1}=0,
\]
which arises in the study of the annihilators of harmonic $k$-forms
\cite{ABDS2011}.

\begin{theorem}[\cite{ABDS2011}]\label{thm:harmpairdual}
Let $\Lambda$ be a continuous real-linear functional on
$\mathcal{HF}_{k}(\mathcal{K})$.  Then there exists a pair
$\omega=(\omega_{k-1},\omega_{k+1})$ satisfying the right-hand version
of \eqref{eq:MTsystem} in $\R^{n}\setminus \mathcal{K}$, with $\omega(\infty)=0$,
such that for every germ $\dot F_{k}\in\mathcal{HF}_{k}(\mathcal{K})$
\[\Lambda[\dot F_{k}]
=\int_{\Gamma_{F_{k}}}
\bigl\langle\omega_{k-1}(y),\eta(y)\bullet F_{k}(y)\bigr\rangle
\,d\mathcal{H}^{n-1}(y)
+\int_{\Gamma_{F_{k}}}
\bigl\langle\omega_{k+1}(y),\eta(y)\wedge F_{k}(y)\bigr\rangle
\,d\mathcal{H}^{n-1}(y),\]
$\Gamma_{F_{k}}$ being a smooth surface surrounding $\mathcal{K}$.
\end{theorem}

Theorem~\ref{thm:harmpairdual} produces a surjective map
\[
J:\mathcal{HP}^{r}_{k}\bigl(\R^{n}\setminus \mathcal{K},0\bigr)
\longrightarrow \mathcal{HF}_{k}(\mathcal{K})^{*},
\]
so that $\mathcal{HF}_{k}(\mathcal{K})^{*}$ and
$\mathcal{HP}^{r}_{k}(\R^{n}\setminus \mathcal{K},0)/\ker J$ are isomorphic.  The
kernel is described as follows.

\begin{theorem}[\cite{ABDS2011}]\label{thm:kerJeuclid}
If $\omega=(\omega_{k-1},\omega_{k+1})
\in\mathcal{HP}^{r}_{k}(\R^{n}\setminus \mathcal{K},0)$ belongs to $\ker J$, then
$\omega_{k-1}$ and $\omega_{k+1}$ are \emph{both} harmonic multivector
fields in $\R^{n}\setminus \mathcal{K}$.  The converse does not hold in general;
it does hold when $\mathcal{K}$ is cohomologically trivial, and in that case
these pairs characterize $\ker J$ completely.
\end{theorem}

\begin{theorem}[\cite{ABDS2011}]\label{thm:formsdual}
Let $\mathcal{K}$ be a compact cohomologically trivial subset of $\R^{n}$.  Then
the dual of the space of germs of harmonic differential $k$-forms on
$\mathcal{K}$ is isomorphic to the space of solutions $(v_{k-1},v_{k+1})$ in
$\R^{n}\setminus \mathcal{K}$ of
\[
dv_{k-1}-d^{*}v_{k+1}=0,\qquad d^{*}v_{k-1}=0,\qquad dv_{k+1}=0,
\]
vanishing at infinity, modulo the subspace of those pairs whose two
components are both harmonic in $\R^{n}\setminus \mathcal{K}$.
\end{theorem}

The architecture displayed by Theorems~\ref{thm:eucdual}--\ref{thm:formsdual}---a clean isomorphism
at the level of the ``full'' Clifford-algebra-valued theory, and a
representation by pairs, surjective but not injective, at the level of
the refined one, with a kernel described by an extra annihilation
requirement---is exactly the architecture we shall meet again in
Subsections~\ref{subsec:matrixdual} and \ref{subsec:scalardual}, with
circulant matrices in the role of Clifford-algebra-valued functions and
$h$-monogenic functions in the role of harmonic multivector fields.

\subsection{Algebraic analysis: why the Hermitian duality must be
proved directly}\label{subsec:algebraic}

When studying duality for systems of partial differential equations the
standard route is algebraic: one analyses the module associated with
the system, proves the vanishing of suitable cohomology groups, and
deduces the duality from general principles \cite{SabSom2025}.  This
route is closed in the Hermitian setting, and it is worth recording
why, since the failure is precisely what makes the matrix formalism
valuable.

By the last part of Subsection~\ref{subsec:hermitian}, the Hermitian
Dirac operators shift the spinor grading, whence they admit a block
representation.

\begin{proposition}[\cite{DES2009,CSS2012,SabSom2025}]
\label{prop:blockmatrices}
With respect to the decomposition
$\mathbb{S}=\bigoplus_{j=0}^{n}\mathbb{S}^{j}$, and with entries in
$\C[Z,\Zd]$, the operator $\dZ$ is represented by the lower-shift block matrix 
whose only non-zero blocks, on the subdiagonal are $d_{i}=\pi_{\mathbb{S}^{i}}\circ[\dZ]|_{\mathbb{S}^{i-1}}$,
$1\le i\le n$, and $\dZd$ by the upper-shift block matrix whose only non-zero blocks lie on the subdiagonal and are $\delta_{i}=\pi_{\mathbb{S}^{i}}\circ[\dZd]|_{\mathbb{S}^{i+1}}$,
$0\le i\le n-1$.  One has $d_{i}d_{i-1}=0$ and $\delta_{i}\delta_{i+1}=0$ whenever these compositions are defined.
\end{proposition}

The following theorem makes the structure of this resolution explicit and describes the maps attached to the first syzygies, writing an element of the first syzygy module as $(a,b)$, with $a$ and $b$ corresponding to $\dZ$ and $\dZd$, see \cite[Theorem~9]{SabSom2025}.
\begin{theorem}[\cite{DES2009,SabSom2025}]\label{thm:resolution}
The free resolution of the module associated with the Hermitian Dirac system in dimension $2n$ is linear of length $n$.  The first syzygies are given by the relations
\[
d_{i}\,a|_{\mathbb{S}^{i-1}}=0\quad(2\le i\le n),
\qquad
\delta_{i}\,b|_{\mathbb{S}^{i+1}}=0\quad(0\le i\le n-2),
\]
which amount to $2^n-n-1$ complex relations for each of $a$ and $b$, i.e.\ $2(2^n-n-1)=\beta_2$ in total\footnote{In \cite[Theorem~3.3]{DES2009}, and in \cite{SabSom2025}, the count $2^n-n-1$  must be read as referring to each of $a$ and $b$ separately. The relations for $a$ take values in $\bigoplus_{i=2}^n\mathbb S^i$ and those for $b$ in $\bigoplus_{i=0}^{n-2}\mathbb S^i$, so each family contributes $2^n-n-1$ scalar relations, in agreement with $\beta_2=2(2^n-n-1)$ \cite[Corollary~3.5]{DES2009}.}

The Betti numbers of the module are
$\beta_{0}=2^{n}$ and
$\beta_{i}=2\sum_{j=i}^{n}\binom{n}{j}$, $1\le i\le n$.
\end{theorem}

\begin{remark}[\cite{DES2009,SabSom2025}]
Unlike what happens for the Cauchy--Riemann or the Cauchy--Fueter systems, the compatibility conditions of
Theorem~\ref{thm:resolution} never involve $a$ and $b$ simultaneously.
\end{remark}

\begin{corollary}[Failure of the Hartogs phenomenon,
\cite{DES2009,SabSom2025}]\label{cor:noHartogs}
Let $\mathcal{K}$ be a compact convex subset of an open set $\Theta\subseteq\R^{2n}$,
$n>1$, and let $f$ be $h$-monogenic in $\Theta\setminus \mathcal{K}$.  Then $f$ cannot
in general be extended to an $h$-monogenic function on $\Theta$.
\end{corollary}

The cohomology vanishing needed for a formal derivation of duality
theorems does not hold for this module \cite{DES2009}; every duality
statement in the Hermitian setting must therefore be established
directly.  At the matrix level, however, the picture is as simple as
one could wish.

\begin{proposition}[\cite{CSS2012}]\label{prop:trivialsyz}
The circulant matrix Dirac operator $\Dop$ has trivial first syzygies:
the matrix of its first syzygies is the zero matrix.
\end{proposition}

\begin{corollary}[Poincar\'e lemma for circulant matrices,
\cite{CSS2012}]\label{cor:poincare}
Let $\Theta\subseteq\R^{2n}$ be a convex open set and let $\Gmat$ be a
circulant matrix with entries in $C^{\infty}(\Theta)$.  Then there exists a
circulant matrix $\Fmat$ with $C^{\infty}(\Theta)$ entries such that
$\Dop\,\Fmat=\Gmat$ in $\Theta$.
\end{corollary}

\begin{corollary}[Mittag--Leffler theorem, \cite{CSS2012}]
\label{cor:mittag}
Let $\mathcal{U}=\{\Theta_{\alpha}\}$ be a locally finite open covering of
an open sets and let
$F_{\alpha,\beta}\in\HM_{r}(\Theta_{\alpha}\cap \Theta_{\beta})$ be circulant
matrices satisfying the cocycle condition
$F_{\alpha,\beta}+F_{\beta,\gamma}+F_{\gamma,\alpha}=\Zero$ on
$\Theta_{\alpha}\cap \Theta_{\beta}\cap \Theta_{\gamma}$.  Then for every $\alpha$
there is $G_{\alpha}\in\HM_{r}(\Theta_{\alpha})$ with
$F_{\alpha,\beta}=G_{\beta}-G_{\alpha}$ on
$\Theta_{\alpha}\cap \Theta_{\beta}$.
\end{corollary}

\begin{remark}\label{rem:sheaf}
The assignment $\Theta\mapsto\HM(\Theta)$ is a sheaf of $\CM$-modules, so that
Corollary~\ref{cor:mittag} may be restated as
$H^{1}(\Theta,\HM_{r})=\Zero$ for every open $\Theta\subseteq\R^{2n}$ \cite[Remark~2.19]{CSS2012}. Moreover the
elements of this sheaf enjoy the unique continuation property, a
consequence of the vanishing of the zeroth $\mathrm{Ext}$-module
associated with $\HM_{r}$, see \cite[Remark~2.20]{CSS2012} and \cite[Theorem~3.9]{DES2009}. 
It also follows from the harmonicity of $H$-monogenic matrix functions; see Subsection~\ref{subsec:circulant}.
\end{remark}

\subsection{Duality for \texorpdfstring{$H$}{H}-monogenic matrix
functions}\label{subsec:matrixdual}

Equipped with the topology of uniform convergence on compact sets, for every open $\Theta\subseteq\R^{2n}$ the sets $\HM_{r}(\Theta)$ and
$\HM_{l}(\Theta)$ are, respectively, left and right Fr\'echet
$\CM$-modules for this topology, and in fact Montel modules, whence
$\HM_{l}(\mathcal{K})$ and $\HM_{r}(\mathcal{K})$ are $LF$-modules \cite{CSS2012};
moreover a right $\CM$-linear functional $\Lambda$ on $\HM_{l}(\mathcal{K})$ is
continuous precisely when $\Lambda\circ\rho_{\Theta}$ is continuous on
$\HM_{l}(\Theta)$ for $\Theta\supset \mathcal{K}$, $\rho_{\Theta}$
denoting the canonical projection onto germs \cite{A47}.

For a compact set $\mathcal{K}\subset\R^{2n}$ with $\R^{2n}\setminus \mathcal{K}$ connected, put
\[
\widetilde{\HM}_{l}(\mathcal{K})=\bigcup_{\Omega\supset \mathcal{K}}\HM_{l}(\Omega),
\qquad \Omega\text{ open},
\]
and let $\HM_{l}(\mathcal{K})$ be the quotient right $\CM$-module of germs, that
is, $\widetilde{\HM}_{l}(\mathcal{K})$ modulo the equivalence
$g^{1}_{2}\sim f^{1}_{2}$ if and only if $g^{1}_{2}=f^{1}_{2}$ on some
open neighbourhood of $\mathcal{K}$; the modules $\HM_{r}(\mathcal{K})$ are defined
symmetrically.  Finally,
$\HM^{\infty}_{l}(\R^{2n}\setminus \mathcal{K})$ denotes the module of left
$H$-monogenic matrix functions outside $E$ vanishing at infinity, and
similarly with $l$ and $r$ interchanged.

Before the duality statements, \cite{A47} establishes a
characterization of traces by annihilation, which is of independent
interest.

\begin{theorem}[\cite{A47}]\label{thm:annihilators}
Let $f^{1}_{2},g^{1}_{2}\in C^{0,\nu}(\Gamma)$, $0<\nu<1$, and consider
the pairing
\[
\bigl\langle f^{1}_{2},g^{1}_{2}\bigr\rangle_{\Gamma}
=\int_{\Gamma}f^{1}_{2}(X)\,N_{(Z,\Zd)}\,g^{1}_{2}(X)
\,d\mathcal{H}^{2n-1}.
\]
Then $\langle f^{1}_{2},g^{1}_{2}\rangle_{\Gamma}=\Zero$ for all
$f^{1}_{2}\in\HM_{r}(\Omega^{+})$ continuous up to $\Gamma$ if and only
if $g^{1}_{2}$ is the trace on $\Gamma$ of a matrix function in
$\HM_{l}(\Omega^{+})$; symmetrically,
$\langle f^{1}_{2},g^{1}_{2}\rangle_{\Gamma}=\Zero$ for all
$g^{1}_{2}\in\HM_{l}(\Omega^{+})$ continuous up to $\Gamma$ if and only
if $f^{1}_{2}$ is the trace of a matrix function in
$\HM_{r}(\Omega^{+})$.
\end{theorem}

\begin{theorem}[Integral representation of the dual, \cite{A47}]
\label{thm:dualmatrix}
Let $\Lambda$ be a continuous right $\CM$-linear functional on
$\HM_{l}(\mathcal{K})$.  Then there exists a unique $\CM$-valued function
$R^{1}_{2}\in\HM_{r}\bigl(\R^{2n}\cup\{\infty\}\setminus \mathcal{K}\bigr)$ with
$R^{1}_{2}(\infty)=\Zero$ such that, for all germs
$\dot{g}^{1}_{2}\in\HM_{l}(K)$,
\[
\Lambda\bigl[\dot{g}^{1}_{2}\bigr]
=\int_{\Gamma_{g}}R^{1}_{2}(X)\,N_{(Z,\Zd)}\,g^{1}_{2}(X)
\,d\mathcal{H}^{2n-1},
\]
where $g^{1}_{2}$ is any representative of the germ and $\Gamma_{g}$ is
a smooth surface surrounding $\mathcal{K}$.  Equivalently,
$\HM_{l}(\mathcal{K})'\cong\HM^{\infty}_{r}(\R^{2n}\setminus \mathcal{K})$.
\end{theorem}

The dual statement, together with the reflexivity of the modules
involved, was supplied in \cite{CSS2012}, where the representing matrix
function is constructed as the analogue of the Fantappi\`e indicatrix
of the functional, $\Lambda\mapsto\Lambda[\Ekernel(Z-V)]$.

\begin{theorem}[\cite{CSS2012}]\label{thm:CSSdual1}
Let $\mathcal{K}\subset\R^{2n}$ be compact.  Then there are $\CM$-module
isomorphisms
\[
\HM_{r}(\mathcal{K})'\cong\HM^{\infty}_{l}\bigl(\R^{2n}\setminus \mathcal{K}\bigr),
\qquad
\HM^{\infty}_{l}\bigl(\R^{2n}\setminus \mathcal{K}\bigr)'\cong\HM_{r}(\mathcal{K}).
\]
\end{theorem}

\begin{theorem}[\cite{CSS2012}]\label{thm:CSStop}
Both isomorphisms of Theorem~\ref{thm:CSSdual1} are topological.
\end{theorem}

\begin{theorem}[Cohomological form of the duality, \cite{CSS2012}]
\label{thm:CSScohom}
Let $\mathcal{K}\subset\R^{2n}$ be compact and $\Theta$ an open neighbourhood of $\mathcal{K}$.
Then there is a topological $\CM$-module isomorphism
\[
\HM_{r}(\mathcal{K})'\cong H^{1}\bigl(\Theta,\Theta\setminus \mathcal{K},\HM\bigr).  
\]
\end{theorem}

Theorem~\ref{thm:CSScohom} rests on the unique continuation property
and on the Mittag--Leffler theorem of
Subsection~\ref{subsec:algebraic}, which together identify
$H^{1}(\Theta,\Theta\setminus \mathcal{K},\HM)$ with the quotient
$\HM(\Theta\setminus \mathcal{K})/\HM(\Theta)$; it is the exact analogue, in the present
setting, of the classical cohomological formulation of
K\"othe--Grothendieck duality.

\subsection{The scalar level: pairs, differential forms and the dual of
\texorpdfstring{$h(\mathcal{K})$}{h(\mathcal{K})}}\label{subsec:scalardual}

Descending from the matricial to the scalar level produces a
description of the dual of the space $h_{l}(\mathcal{K})$ of germs of
$h$-monogenic functions, which is necessarily somewhat weaker: exactly
as in Theorem~\ref{thm:harmpairdual}, the representing object is now a
\emph{pair}.

\begin{theorem}[\cite{A47}]\label{thm:dualscalar}
Let $\Lambda$ be a continuous right $\Cln$-linear functional on
$h_{l}(\mathcal{K})$. Then there exists a pair $(\omega,\omega^{\dagger})$, vanishing at infinity and satisfying in $\R^{2n}\setminus \mathcal{K}$ the
conditions\footnote{This corrects an inaccuracy that appeared in \cite[System~(7), p. ~354]{A47}, where the domains of validity of the conditions were misstated.}
\begin{equation}\label{eq:dualsystem}
\begin{cases}
[\omega]\dZ+[\omega^{\dagger}]\dZd=0,\\
[\omega]\dZd=0,\quad [\omega^{\dagger}]\dZ=0
\end{cases}
\end{equation}
such that, for all $\dot{g}\in h_{l}(K)$,
\[\Lambda[\dot{g}]=\int_{\Gamma_g}
\bigl(\omega N-\omega^{\dagger}N^{\dagger}\bigr)\,g
\,d\mathcal{H}^{2n-1}.\]
\end{theorem}

Theorem~\ref{thm:dualscalar} produces a surjective $\Cln$-linear map
$J$ from the space of pairs $(\omega,\omega^{\dagger})$ solving
\eqref{eq:dualsystem} and vanishing at infinity onto the dual
$h^{'}_{l}(\mathcal{K})$.  Its kernel is identified as follows.

\begin{theorem}[\cite{A47}]\label{thm:kernelJ}
If $(\omega,\omega^{\dagger})\in\ker J$, then $\omega$ and
$\omega^{\dagger}$ are both right $h$-monogenic in $\R^{2n}\setminus \mathcal{K}$.
\end{theorem}

The comparison with Theorems~\ref{thm:harmpairdual}--\ref{thm:kerJeuclid} is immediate: in
both cases the dual of the refined object is represented by a pair
subject to a first-order system coupling its two components, and the
ambiguity of the representation is measured by the pairs whose two
components separately satisfy the refined equation.  In \cite{CSS2012}
this description is completed into an honest isomorphism by passing to
a quotient of spaces of differential forms; we conclude the section
with that statement.

\begin{proposition}[\cite{CSS2012}]\label{prop:CSSsystem}
Let $\mathcal{K}\subset\R^{2n}$ be compact and let $\Lambda\in h(\mathcal{K})'$.  Then, on
$\R^{2n}\setminus \mathcal{K}$,
\[
\dZ\Lambda\mathcal{E}(Z-V)+\dZd\Lambda\mathcal{E}^{\dagger}(Z-V)=0,
\qquad
\dZd\Lambda\mathcal{E}(Z-V)=0,
\qquad
\dZ\Lambda\mathcal{E}^{\dagger}(Z-V)=0,
\]
the last two relations being automatic consequences of
Remark~\ref{rem:conventions}.  Moreover the action of $\Lambda$ on
$g\in h(K)$ is given by
\[\Lambda[g]=\int_{\Gamma_{g}}g\bigl(Z,\Zd\bigr)
\Bigl(d\sigma_{Z}\bigl[\Lambda\mathcal{E}(Z-V)\bigr]
-d\sigma_{\Zd}\bigl[\Lambda\mathcal{E}^{\dagger}(Z-V)\bigr]\Bigr),\]
$\Gamma_{g}$ surrounding $\mathcal{K}$ inside a neighbourhood to which $g$
extends, while the companion integral obtained by interchanging
$\mathcal{E}$ and $\mathcal{E}^{\dagger}$ vanishes identically.
\end{proposition}

\begin{lemma}[\cite{CSS2012}]\label{lem:CSSdd}
Put $\partial=\sum_j dz_j\wedge\partial/\partial z_j$ and $\bar\partial=\sum_j d z_j^c\wedge\partial/\partial z_j^c,$
acting on the forms of Subsection~\ref{subsec:forms}.  Then
\[\partial\bigl(d\sigma_{Z}\,\omega\bigr)=-dW(Z,\Zd)\,\dZ\omega,
\qquad
\bar\partial\bigl(d\sigma_{Z}\,\omega\bigr)=0,\]
\[\partial\bigl(d\sigma_{\Zd}\,\omega\bigr)=0,
\qquad
\bar\partial\bigl(d\sigma_{\Zd}\,\omega\bigr)=dW(Z,\Zd)\,\dZd\omega.\]
Consequently
$\bigl(\partial+\bar\partial\bigr)\bigl(d\sigma_{Z}\omega
-d\sigma_{\Zd}\omega^{\dagger}\bigr)=0$ if and only if
$\dZ\omega+\dZd\omega^{\dagger}=0$.
\end{lemma}

Denote by $\mathcal{M}^{\dagger\infty}_{(n-1,n)}(\Theta)$ the space of forms
$d\sigma_{Z}\omega$ whose coefficient satisfies $\dZd\omega=0$
in $\Theta$ and vanishes at infinity, and by
$\mathcal{M}^{\infty}_{(n,n-1)}(\Theta)$ the space of forms
$d\sigma_{\Zd}\omega^{\dagger}$ with $\dZ\omega^{\dagger}=0$ and vanishing at
infinity; when the coefficients are moreover $h$-monogenic the
corresponding spaces are written $\HM^{\infty}_{(n-1,n)}(\Theta)$ and
$\HM^{\infty}_{(n,n-1)}(\Theta)$.

\begin{theorem}[\cite{CSS2012}]\label{thm:CSSscalar}
For every compact set $\mathcal{K}\subset\R^{2n}$ there is an algebraic
isomorphism
\[h(\mathcal{K})'\;\cong\;
\frac{\ker\Bigl\{\partial+\bar\partial:
\mathcal{M}^{\dagger\infty}_{(n-1,n)}\bigl(\R^{2n}\setminus \mathcal{K}\bigr)
+\mathcal{M}^{\infty}_{(n,n-1)}\bigl(\R^{2n}\setminus \mathcal{K}\bigr)
\longrightarrow C_{(n,n)}\bigl(\R^{2n}\setminus \mathcal{K}\bigr)\Bigr\}}
{\HM^{\infty}_{(n-1,n)}\bigl(\R^{2n}\setminus \mathcal{K}\bigr)
+\HM^{\infty}_{(n,n-1)}\bigl(\R^{2n}\setminus \mathcal{K}\bigr)}.\]
\end{theorem}

% ==================================================================
\section{Perturbed and iterated Hermitian Dirac operators}\label{sec:helmholtz}
% ==================================================================

\subsection{Factorization and perturbed Hermitian Dirac operators}\label{subsec:helmfact}

In Euclidean Clifford analysis the Helmholtz operator is factorized by
means of the perturbed, or $k$-Dirac, operators
$\dX\pm k$.  The Hermitian analogue, developed in \cite{A45}, splits
this factorization further.  Consider the Helmholtz equation for a
$\Cln$-valued function $g$,
\begin{equation}\label{eq:helmholtz}
\Delta_{2n}g+k^{2}g=0,\qquad k\in\R,
\end{equation}
and introduce the complex constants
\begin{equation}\label{eq:Kconst}
K=-\tfrac14(k-ik),\qquad K^{\dagger}=-\tfrac14(k+ik).
\end{equation}
A direct computation gives the twofold factorization
\begin{equation}\label{eq:helmfact}
\begin{split}
\Delta_{2n}+k^{2}
&=4\Bigl[(\dZ+K)\bigl(\dZd+K^{\dagger}\bigr)
+\bigl(\dZd-K^{\dagger}\bigr)(\dZ-K)\Bigr]\\[1mm]
&=-(\dX+k)(\dX-k)=-(\dXt+k)(\dXt-k),
\end{split}
\end{equation}
which for $k=0$ reduces to \eqref{eq:hermlaplace}.  At the matrix level
one sets
\[
\DopK=\begin{pmatrix}
\dZ+K & \dZd-K^{\dagger}\\[1mm]
\dZd-K^{\dagger} & \dZ+K
\end{pmatrix},
\]
so that $\boldsymbol{\Delta}+k^{2}\mathbf{1}
=4\bigl(\DopK\bigr)\bigl(\DopK\bigr)^{\dagger}
=4\bigl(\DopK\bigr)^{\dagger}\bigl(\DopK\bigr)$.

\begin{definition}[\cite{A45}]\label{def:Hhelmholtz}
The circulant matrix function $\Gmat$ is \emph{(left) Helmholtz
$H$-monogenic} in $\Omega$ if $\DopK[\Gmat]=\Zero$ there.  We write
$\HM^{K}(\Omega)$ for the corresponding space.  Explicitly,
\begin{align*}
\DopK\bigl[\Gmat\bigr]=\Zero
&\iff
\begin{cases}
(\dZ+K)g_{1}+\bigl(\dZd-K^{\dagger}\bigr)g_{2}=0,\\[1mm]
\bigl(\dZd-K^{\dagger}\bigr)g_{1}+(\dZ+K)g_{2}=0,
\end{cases}\\[3mm]
&\iff
\begin{cases}
(\dX+k)g_{1}-(\dX+k)g_{2}=0,\\[1mm]
(\dXt+k)g_{1}+(\dXt+k)g_{2}=0 .
\end{cases}
\end{align*}
A $\Cln$-valued function $g$ is \emph{Helmholtz $h$-monogenic} if
$(\dZ+K)g=0=(\dZd-K^{\dagger})g$, equivalently
$(\dX+k)g=0=(\dXt+k)g$; this is the diagonal case $G_{0}$.
\end{definition}

The relevant kernels are built from the fundamental solution
$\theta_{k}$ of $\Delta_{2n}+k^{2}$ and from the pair
$(E_{k},E^{|}_{k})$ of fundamental solutions of $\dX+k$ and
$\dXt+k$:
\begin{equation}\label{eq:helmkernels}
\mathcal{E}_{K}=-\bigl(E_{k}+iE^{|}_{k}\bigr)
=4\bigl(\dZd+K^{\dagger}\bigr)\theta_{k},
\qquad
\mathcal{E}^{\dagger}_{K}=E_{k}-iE^{|}_{k}
=4(\dZ-K)\theta_{k},
\end{equation}
which for $k=0$ reduce to $\mathcal{E}$ and $\mathcal{E}^{\dagger}$.
The matrix $\mathbf{E}_{K}$ built from $\mathcal{E}_{K}$ and
$\mathcal{E}^{\dagger}_{K}$ is a fundamental solution of $\DopK$, and the
perturbed Cauchy and Teodorescu operators
$\mathcal{C}^{K}_{\Gamma}$, $\mathcal{T}^{K}_{\Omega}$ are defined by
\eqref{eq:cauchyint}--\eqref{eq:teo} with $\Ekernel$ replaced by
$\Ekernel_{K}$.

\begin{remark}
By \eqref{eq:helmholtz}, $(\partial_Z+A)g=0=(\partial_{Z^\dagger}-A^\dagger)g$ amounts to $(\partial_X+k)g=0=(\partial_{X|}+k)g$. Since $\partial_X$ and $\partial_{X|}$ anticommute, applying $\partial_{X|}$ to the first equation and $\partial_X$ to the second and adding yields $2k^2g=0$. Hence for $k\neq0$ the only Helmholtz $h$-monogenic function is $g\equiv0$.
\end{remark}

\subsection{Integral formulae and the perturbed Hilbert transform}\label{subsec:helmformulae}

\begin{proposition}[Borel--Pompeiu, Cauchy and inversion formulae,
\cite{A45}]\label{prop:helmformulae}
Let $\Omega$ have Liapunov, piecewise smooth or Lipschitz boundary
$\Gamma$.
\begin{enumerate}[label=\upshape(\roman*)]
\item If $\Gmat\in C^{1}(\overline{\Omega})$, then
\[
\mathcal{C}^{K}_{\Gamma}\bigl[\Gmat\bigr](Y)
+\mathcal{T}^{K}_{\Omega}\bigl[\DopK\Gmat\bigr](Y)
=\begin{cases}\lam\,\Gmat(Y),&Y\in\Omega^{+},\\ 0,&Y\in\Omega^{-}.
\end{cases}
\]
\item If in addition $\Gmat\in\HM^{K}(\Omega)$, then
$\mathcal{C}^{K}_{\Gamma}[\Gmat](Y)=\lam\,\Gmat(Y)$ for
$Y\in\Omega^{+}$ and $=0$ for $Y\in\Omega^{-}$.
\item If $\Gmat\in C^{1}(\Omega)\cap C(\overline{\Omega})$, then
$\DopK\,\mathcal{T}^{K}_{\Omega}[\Gmat](Y)=\lam\,\Gmat(Y)$ for
$Y\in\Omega^{+}$ and $=0$ for $Y\in\Omega^{-}$.
\end{enumerate}
\end{proposition}

\begin{proposition}[Liouville-type theorem, \cite{A45}]
\label{prop:liouville}
Let $\Gmat\in\HM^{K}(\R^{2n})$ with
$\|\Gmat\|=o\bigl(|X|^{-\frac{2n-1}{2}}\bigr)$ as $|X|\to\infty$.
Then $\Gmat\equiv\Zero$ in $\R^{2n}$.
\end{proposition}

\begin{proposition}[Perturbed Plemelj--Sokhotski formulae, \cite{A45}]
\label{prop:helmPS}
Let $\Gmat\in C^{0,\nu}(\Gamma)$.  Then, at every $U\in\Gamma$,
\[
\lim_{\substack{Y\to U\\ Y\in\Omega^{\pm}}}
\mathcal{C}^{K}_{\Gamma}\bigl[\Gmat\bigr](Y)
= \lam  \left(\pm\tfrac12\Gmat(U)
+\tfrac12\mathcal{H}^{K}_{\Gamma}\bigl[\Gmat\bigr](U)\right)
=\pm\lam P^{\pm}_{K}\bigl[\Gmat\bigr](U),
\]
where $P^{\pm}_{K}:=\tfrac12\left(\mathbf{1} \pm \mathcal{H}^{K}_{\Gamma}\right).$
\end{proposition}

The heart of the matter is that the perturbation is compact:
the kernel $\Ekernel_{K}-\Ekernel$ behaves as $o(|X|^{-2n+2})$ as
$|X|\to 0$, hence is weakly singular, so that
$\mathcal{H}^{K}_{\Gamma}$ coincides with $\Hilb$ up to a compact
operator and \cite{A45} shows that the structural properties survive.

\begin{theorem}[\cite{A45}]\label{thm:helmhilbert}
The perturbed Hermitian Hilbert transform
$\mathcal{H}^{K}_{\Gamma}$ is a bounded operator on
$C^{0,\nu}(\Gamma)$ and an involution,
$\bigl(\mathcal{H}^{K}_{\Gamma}\bigr)^{2}=\mathbf{1}$.  Consequently
$P^{+}_{K}$ and $P^{-}_{K}$ are mutually complementary bounded
projection operators onto the traces of $\HM^{K}(\Omega^{\pm})$, and
there is $c_{\Gamma}>0$ with
\begin{equation}\label{eq:Pminusbound}
\bigl\|P^{-}_{K}\bigl[\Gmat\bigr]\bigr\|_{\nu,\Gamma}
\le c_{\Gamma}\bigl\|\Gmat\bigr\|_{\nu,\Gamma},
\qquad \Gmat\in C^{0,\nu}(\Gamma).
\end{equation}
\end{theorem}

\begin{theorem}[$L^{p}$ version, \cite{A45}]\label{thm:helmLp}
Let $1<p<\infty$ and let $\Gamma$ be a Lipschitz hypersurface.  Then
$\mathcal{H}^{K}_{\Gamma}$ and $P^{\pm}_{K}$ are bounded linear
operators on $L^{p}(\Gamma)$:
\[
\bigl\|\mathcal{H}^{K}_{\Gamma}[\Gmat]\bigr\|_{p}\le c_{p}\|\Gmat\|_{p},
\qquad
\bigl\|P^{\pm}_{K}[\Gmat]\bigr\|_{p}\le c_{p}\|\Gmat\|_{p} .
\]
\end{theorem}

Finally, one obtains a characterization of the solutions of the
matricial Helmholtz equation themselves.

\begin{theorem}[\cite{A45}]\label{thm:helmsolutions}
Let $\Gmat\in C^{1}(\overline\Omega)$ and $\Fmat\in C^{1}(\Omega)\cap\HM^{K}(\Omega)$. The circulant matrix
function $\Jmat=\mathcal{T}^{K}_{\Omega}\bigl[\Gmat\bigr]+\Fmat$
solves the matricial Helmholtz equation
\[
\boldsymbol{\Delta}\Jmat+k^{2}\Jmat=\Zero
\]
if and only if $\bigl(\mathcal D^K_{(Z,Z^\dagger)}\bigr)^\dagger\bigl[G_2^1\bigr]=\mathcal O$ in $\Omega$
\end{theorem}

\subsection{Riemann problems for the Hermitian matricial operator}\label{subsec:riemann}

The above operator calculus allows one to transfer the theory of
Riemann (transmission) problems to the Hermitian matricial setting.

\begin{problem}[Riemann problem, \cite{A45}]\label{prob:riemann}
Given $\Gmat,\Fmat\in C^{0,\nu}(\Gamma)$, find
$\Phimat\in\HM^{K}(\Omega^{\pm})$, vanishing at infinity, whose
boundary values satisfy the transmission condition
\begin{equation}\label{eq:riemann}
\bigl[\Phimat\bigr]^{+}(U)
=\Gmat(U)\bigl[\Phimat\bigr]^{-}(U)+\Fmat(U),
\qquad U\in\Gamma .
\end{equation}
\end{problem}

\begin{theorem}[Constant coefficient case, \cite{A45}]
\label{thm:riemannconst}
Let $\Gmat\in\CM$ be a complex constant invertible matrix.  Then
\eqref{eq:riemann} has a unique solution, given by
\[
\Phimat(Y)=\frac{1}{\lambda_n}\,\mathfrak{X}^{1}_{2}(Y)\int_{\Gamma}
\mathbf{E}_{K}(Z-V)\,N_{(Z,\Zd)}\bigl[\Gmat\bigr]^{-1}\Fmat(X)
\,d\mathcal{H}^{2n-1},
\qquad
\mathfrak{X}^{1}_{2}(Y)=\begin{cases}\Gmat,&Y\in\Omega^{+},\\
\mathbf{1},&Y\in\Omega^{-}.\end{cases}
\]
\end{theorem}

\begin{theorem}[Variable coefficient case, \cite{A45}]
\label{thm:riemannvar}
Let $\Gmat,\Fmat\in C^{0,\nu}(\Gamma)$ and assume the smallness
condition
\begin{equation}\label{eq:smallness}
c_{\Gamma}\bigl\|\Gmat-\mathbf{1}\bigr\|_{\nu,\Gamma}<1,
\end{equation}
with $c_{\Gamma}$ the constant of \eqref{eq:Pminusbound}.  Then the
Riemann boundary value problem \eqref{eq:riemann} has a unique
solution.
\end{theorem}

\begin{remark}[\cite{A45}]
By Theorem~\ref{thm:helmLp}, Theorems~\ref{thm:riemannconst} and
\ref{thm:riemannvar} translate verbatim to $L^{p}$ data, all formulae
being reinterpreted in the Lebesgue sense.
\end{remark}

\medskip
The same operator calculus supports a second family of Riemann-type
problems, in which the perturbation of the Hermitian Dirac operator is
replaced by its \emph{iteration}.  This line was opened by Gu and Fu
\cite{GF2014} and carried further by Gu \cite{Gu2019}; we describe it
here because both its tools and the mechanism of its proofs are exactly
those of Theorems~\ref{thm:riemannconst}--\ref{thm:riemannvar}.

\subsection{Riemann problems for H-2-monogenic matrix functions}

In this subsection, following \cite{GF2014,Gu2019}, we use the equivalent norm $\|G_2^1\|_2=\bigl(|g_1|^2+|g_2|^2\bigr)^{1/2}$, and the H\"older norm built on it.

\begin{definition}[\cite{GF2014,Gu2019}]\label{def:H2monogenic}
A circulant matrix function $\Gmat$ with entries of class $C^{2}$ is
called \emph{(left) $H$-$2$-monogenic} in $\Omega$ if
$\bigl(\Dop\bigr)^{2}\bigl[\Gmat\bigr]=\Zero$ in $\Omega$.
\end{definition}

Let
\begin{equation}\label{eq:E1kernel}
\mathcal{E}_{1}(Z-V)=\frac{1}{\sigma_{2n}(2-2n)}
\begin{pmatrix}
0 & \dfrac{4}{|X-Y|^{2n-2}}\\[3mm]
\dfrac{4}{|X-Y|^{2n-2}} & 0
\end{pmatrix},
\qquad X\ne Y,
\end{equation}
assembled from the fundamental solution of the Laplacian, which plays
for $\bigl(\Dop\bigr)^{2}$ the role played by $\Ekernel$ for $\Dop$.

A circulant matrix function $G_2^1$ of class $C^2$ is $H$-$2$-monogenic in $\Omega$ if and only if $g_1$ and $g_2$ are harmonic in $\Omega$. If $\varXi$ denotes the fundamental solution of $\Delta_{2n}$, the kernel (\ref{eq:E1kernel}) is $\mathcal E_1=4\varXi\,\Pi$, whence $(\mathcal D_{(Z,Z^\dagger)})^2\mathcal E_1=\boldsymbol\delta$ and $\mathcal D_{(Z,Z^\dagger)}\mathcal E_1=\mathcal E_1\mathcal D_{(Z,Z^\dagger)}=\Ekernel$. 

Throughout the remainder of this subsection we assume $n\geq2$; for $n=1$, $\varXi$ has to be replaced by $(2\pi)^{-1}\log|X|$. In this light, Proposition~\ref{prop:secondBP} is Green's representation formula written entrywise, and Theorem~\ref{thm:maxmodulus} and  \ref{thm:GuLiouville} are matrix versions of classical properties of harmonic functions.
\begin{proposition}[Second-order Borel--Pompeiu formula,
\cite{GF2014,Gu2019}]\label{prop:secondBP}
Let $\Sigma\subset\Omega$ be a $2n$-dimensional compact oriented
differentiable manifold with $C^{\infty}$ boundary $\partial\Sigma$ and
let $\Gmat$ have entries in $C^{2}(\Omega,\Cln)$.  Then, for
$Y\notin\partial\Sigma$,
\begin{multline*}
\int_{\partial\Sigma}\Ekernel(Z-V)\,\,N_{(Z,\Zd)}\,\Gmat(X)\,
d\mathcal{H}^{2n-1}
-\int_{\partial\Sigma}\mathcal{E}_{1}(Z-V)\,\,N_{(Z,\Zd)}\,
\Dop\Gmat(X)\,d\mathcal{H}^{2n-1}
\\
{}+\int_{\Sigma}\mathcal{E}_{1}(Z-V)
\bigl[\bigl(\Dop\bigr)^{2}\Gmat(X)\bigr]\,dW(Z,\Zd)
=\lam\,\boldsymbol{\chi}_{\Sigma}(Y)\,\Gmat(Y).
\end{multline*}
If $\Gmat$ is $H$-$2$-monogenic the volume term drops out, and an
exterior version of the formula holds for matrix functions vanishing at
infinity.
\end{proposition}
Out of Proposition~\ref{prop:secondBP} the elementary function theory
of $H$-$2$-monogenic matrix functions is obtained.

\begin{theorem}[Mean value and maximum modulus principles,
\cite{GF2014,Gu2019}]\label{thm:maxmodulus}
Let $\Gmat$ be $H$-$2$-monogenic in $\Omega$.
\begin{enumerate}[label=\upshape(\roman*)]
\item For every ball $B(Y,R)\subset\Omega$,
\[
\frac{2n}{\sigma_{2n}R^{2n}}\int_{B(Y,R)}\Gmat(X)\,dV(X)=\Gmat(Y).
\]
\item If $\Omega$ is connected and $\|G_2^1(X)\|_2\leq\|G_2^1(X_0)\|_2$ for some $X_0\in\Omega$ and all $X\in\Omega$, then $G_2^1$ is a constant circulant matrix in $\Omega$.\footnote{Statement (ii) depends on the strict convexity of $\|\cdot\|_2$, which is the norm used in \cite{GF2014,Gu2019}. With the norm $\|\cdot\|$ of Subsection~\ref{subsec:fspaces} it fails: on the unit ball, $g_1\equiv1$ and $g_2=\frac12 \mathfrak f_1^\dagger z_1$ define an $H$-monogenic, hence $H$-$2$-monogenic, matrix function with $\|G_2^1\|\equiv1$ which is not constant.}
\item If $\Omega$ is bounded with Liapunov boundary $\Gamma$ and
$\Gmat\in C^{2}(\Omega)\cap C^{1}(\overline{\Omega})$, then
$\max_{\overline{\Omega}}\|\Gmat\|=\max_{\Gamma}\|\Gmat\|$.
\end{enumerate}
\end{theorem}

\begin{theorem}[\cite{Gu2019}]\label{thm:GuLiouville}
Write $M[r,\Gmat]=\max_{|Y|=r}\|\Gmat(Y)\|$.
\begin{enumerate}[label=\upshape(\roman*)]
\item If $\Gmat$ is $H$-$2$-monogenic in $\R^{2n}$ and
$\liminf_{r\to\infty}M[r,\Gmat]<\infty$, then
\[
\liminf_{r\to\infty}M\bigl[r,\Dop[\Gmat]\bigr]=0 .
\]
\item The same conclusion holds when $\Gmat$ is $H$-$2$-monogenic in
$\Omega^{-}$ only, provided the traces on $\Gamma$ of $\Gmat$ and of
$\Dop[\Gmat]$ belong to $C^{0,\nu}(\Gamma)$ and
$C^{0,\beta}(\Gamma)$ respectively, and
$\liminf_{r\to\infty}M[r,\Gmat]<\infty.$
\end{enumerate}
\end{theorem}

Since $H$-$2$-monogenic matrix functions are harmonic, the hypothesis in (i) actually forces $G_2^1$ to be constant, by the maximum principle and \cite[Corollary~4.8]{GF2014}. Statement (ii), obtained from (i) by subtracting the Cauchy-type integrals built on the two traces, is what makes the exterior part of the following transmission problem well posed.

\begin{problem}[Riemann problem for $H$-$2$-monogenic functions,
\cite{Gu2019}]\label{prob:riemannH2}
Let $\Omega$ be a bounded domain with Liapunov boundary $\Gamma$.
Given $A^{1}_{2},\Fmat, \mathfrak{U}^{1}_{2}\in C^{0,\nu}(\Gamma)$,
$0<\nu<1$, and an invertible \emph{constant} circulant matrix
$B^{1}_{2}$, find $\Gmat$ such that
\[
\begin{cases}
\bigl(\Dop\bigr)^{2}\Gmat=\Zero
&\text{in}\, \R^{2n}\setminus\Gamma,\\[1mm]
\bigl[\Gmat\bigr]^{+}(U)
=\bigl[\Gmat\bigr]^{-}(U)\,A^{1}_{2}(U)+\Fmat(U),
&U\in\Gamma,\\[1mm]
\bigl[\Dop\Gmat\bigr]^{+}(U)
=\bigl[\Dop\Gmat\bigr]^{-}(U)\,B^{1}_{2}+\mathfrak{U}^{1}_{2}(U),
&U\in\Gamma,\\[1mm]
\Gmat(\infty)=\Zero.
\end{cases}
\]
\end{problem}

\begin{remark}[\cite{GF2014}]\label{rem:GF}
For $\Gamma=\partial B(a,R)$ and \emph{constant} invertible
coefficients $A^{1}_{2}$, $B^{1}_{2}$, Problem~\ref{prob:riemannH2}
was solved explicitly in \cite{GF2014} by means of the second-order
Borel--Pompeiu formula, the Plemelj--Sokhotski formulae and a Liouville-
type theorem.  
\end{remark}

\begin{theorem}[Variable coefficient case, \cite{Gu2019}]
\label{thm:GuRiemann}
Assume $A^{1}_{2},\Fmat,\mathfrak U^1_2\in C^{0,\nu}(\Gamma)$,
$0<\nu<1$, that $B^{1}_{2}$ is a constant invertible circulant
matrix, and that
\begin{equation}\label{eq:Gusmallness}
2^{2n-2}\,\bigl\|\mathbf 1-A_2^1\bigr\|_{\nu,\Gamma}\,\bigl(1+c(\nu,n,\Gamma)\bigr)<1,
\end{equation}
where $c(\nu,n,\Gamma)$ denotes the norm of $\Hilb$ on
$C^{0,\nu}(\Gamma)$ furnished by
Theorem~\ref{thm:hilbertprops}.  Then Problem~\ref{prob:riemannH2} has
a unique solution.
\end{theorem}
\begin{remark}
In \cite[(5.6)]{Gu2019} the factor $1+c(\alpha,n,\Gamma)$ appears as $c(\nu,n,\Gamma)$. However, the contraction estimate in the proof of \cite[Theorem~5.2]{Gu2019} produces $1+c$, so that the printed condition only bounds the contraction constant by $1+1/c\leq2$ (note that $\mathcal H_\Gamma^2=\mathbf 1$ forces $c\geq1$) and does not by itself yield a contraction. The constant $2^{2n-2}$ comes from the estimate of products of circulant matrix functions in the H\"older norm built on $\|\cdot\|_2$.
\end{remark}

The mechanism deserves to be isolated, since it is the one already met
in Theorem~\ref{thm:riemannvar}.  The second transmission condition becomes a Riemann problem with
\emph{constant} coefficient for an $H$-monogenic matrix function, which
is solved in closed form by the Cauchy transform $\Cauchy$; a
particular solution $\Jmat$ of the resulting inhomogeneous equation is
then written down by means of the kernel \eqref{eq:E1kernel}, that is, as a single layer potential with the Newtonian kernel, whose image under $\mathcal D_{(Z,Z^\dagger)}$ is the Cauchy integral of the same density; and the difference $\Gmat-\Jmat$ satisfies a first-order transmission problem
with a variable coefficient.  By the Plemelj--Sokhotski formulae the
latter is equivalent to a singular integral equation on
$C^{0,\nu}(\Gamma)$ whose right-hand side involves the data only
through $\mathbf{1}-A^{1}_{2}$ and $\Hilb$; condition
\eqref{eq:Gusmallness} turns it into a contraction of the Banach space
$C^{0,\nu}(\Gamma)$, and the Banach fixed point theorem does the
rest.  The step $\mathcal D_{(Z,Z^\dagger)}G_2^1(\infty)=\mathcal O$, needed to represent $W_2^1$ by a Cauchy integral, follows from the harmonicity of $G_2^1$ and $G_2^1(\infty)=\mathcal O$; uniqueness follows by applying the same reasoning to the homogeneous problem.\footnote{This is not stated explicitly in \cite[Theorem~5.2]{Gu2019}.}

Along the way \cite{Gu2019} reproves, in this notation, the
boundedness and the involutivity of $\Hilb$
(Theorem~\ref{thm:hilbertprops}) as well as the characterization of the
traces of $H$-monogenic functions in $\Omega^{\pm}$ by the equations
$\Hilb[\Gmat]=\pm\Gmat$ (Theorem~\ref{thm:traces}).

\subsection{Nonlinear Riemann problems and the Newton embedding
method}\label{subsec:nonlinear}

All the transmission conditions considered so far are linear in the
unknown.  The nonlinear case was taken up by Gu and Liu
\cite{GuLiu2023} for the perturbed operator $\DopK$ of
Subsection~\ref{subsec:helmfact}, that is, precisely in the Helmholtz
setting of \cite{A45}.

\begin{problem}[Nonlinear Riemann problem, \cite{GuLiu2023}]
\label{prob:nonlinearriemann}
Let $\Omega$ be a bounded domain with Liapunov boundary $\Gamma$ and
let $F^{1}_{2}=F^{1}_{2}\bigl(U,\,\cdot\,,\,\cdot\,\bigr)$ be a given
circulant matrix-valued function of $U\in\Gamma$ and of two matrix
arguments.  Find $\Gmat$ with
\[
\begin{cases}
\DopK\Gmat=\Zero &\text{in }\R^{2n}\setminus\Gamma,\\[1mm]
\bigl[\Gmat\bigr]^{+}(U)=\bigl[\Gmat\bigr]^{-}(U)
+F^{1}_{2}\bigl(U,\bigl[\Gmat\bigr]^{+},\bigl[\Gmat\bigr]^{-}\bigr),
&U\in\Gamma,\\[1mm]
\|\Gmat\|=o\bigl(|U|^{-\frac{2n-1}{2}}\bigr), &|U|\to\infty .
\end{cases}
\]
\end{problem}

The difficulty is structural rather than technical: the Hermitian
Hilbert transform is \emph{not} a compact operator \cite{BDKDS2008}, so
that the usual functional-analytic route---writing the problem as a
Fredholm equation and invoking a compactness argument---is closed.
What replaces it is a continuation, or Newton embedding, argument
resting on two ingredients: the H\"older continuity of the perturbed
Cauchy transform up to the boundary from both sides, and an a priori
estimate for the linear problem.

\begin{lemma}[\cite{GuLiu2023}]\label{lem:GLcauchy}
Let $\Gmat\in C^{0,\nu}(\Gamma)$, $0<\nu<1$.  Then
$\mathcal{C}^{K}_{\Gamma}[\Gmat]$ extends H\"older continuously from
$\Omega^{+}$ to $\overline{\Omega^{+}}$ and from $\Omega^{-}$ to
$\overline{\Omega^{-}}$, with limit values given by
Proposition~\ref{prop:helmPS}, and there is
$c_{\nu,\Gamma}>0$, depending only on $n$, $\nu$ and $\Gamma$, with
\[
\bigl\|\mathcal{C}^{K}_{\Gamma}[\Gmat]\bigr\|
_{C^{0,\nu}(\overline{\Omega^{\pm}})}
\le c_{\nu,\Gamma}\,\bigl\|\Gmat\bigr\|_{\nu,\Gamma}.
\]
\end{lemma}

The proof rests on the weak singularity of $\Ekernel_{K}-\Ekernel$
already exploited in Theorem~\ref{thm:helmhilbert}.  For matrix
functions defined on both sides of $\Gamma$ we write
\[
\bigl\|\Gmat\bigr\|_{\nu}
=\bigl\|\Gmat\bigr\|_{C^{0,\nu}(\overline{\Omega^{+}})}
+\bigl\|\Gmat\bigr\|_{C^{0,\nu}(\overline{\Omega^{-}})},
\qquad
\tilde{c}_{\nu,\Gamma}=2\,c_{\nu,\Gamma}.
\]

\begin{theorem}[A priori estimate, \cite{GuLiu2023}]\label{thm:GLapriori}
Let $\Fmat\in C^{0,\nu}(\Gamma)$ and let $\Gmat$ solve the linear
problem
\[
\DopK\Gmat=\Zero \ \text{ in }\R^{2n}\setminus\Gamma,\qquad
\bigl[\Gmat\bigr]^{+}=\bigl[\Gmat\bigr]^{-}+\Fmat \ \text{ on }\Gamma,
\qquad
\|\Gmat\|=O\bigl(|Y|^{-\frac{2n-1}{2}}\bigr).
\]
Then $\Gmat$ is unique and given by
$\lam\,\Gmat=\mathcal{C}^{K}_{\Gamma}[\Fmat]$, and
\[
\bigl\|\Gmat\bigr\|_{\nu}
\le\tilde{c}_{\nu,\Gamma}\,\bigl\|\Fmat\bigr\|_{\nu,\Gamma}.
\]
\end{theorem}

The hypothesis imposed on the nonlinearity is a Lipschitz condition in
the H\"older norm, quantitatively coupled to the constant of
Theorem~\ref{thm:GLapriori}.

\medskip
\noindent
\textbf{(C1)} For all $\Gmat,\Hmat\in C^{0,\nu}(\Gamma)$ the function
$Y\mapsto F^{1}_{2}\bigl(Y,\Gmat(Y),\Hmat(Y)\bigr)$ belongs to
$C^{0,\nu}(\Gamma)$, and there is a constant $L\ge 0$ with
\begin{equation}\label{eq:GLsmallness}
\tilde{c}_{\nu,\Gamma}\,L<1
\end{equation}
such that, for all
$\Gmat,\widehat{G}^{1}_{2},\Hmat,\widehat{H}^{1}_{2}
\in C^{0,\nu}(\Gamma)$,
\[
\bigl\|F^{1}_{2}(\cdot,\Gmat,\Hmat)
-F^{1}_{2}(\cdot,\widehat{G}^{1}_{2},\widehat{H}^{1}_{2})
\bigr\|_{\nu,\Gamma}
\le L\Bigl[\bigl\|\Gmat-\widehat{G}^{1}_{2}\bigr\|_{\nu,\Gamma}
+\bigl\|\Hmat-\widehat{H}^{1}_{2}\bigr\|_{\nu,\Gamma}\Bigr].
\]

\begin{theorem}[\cite{GuLiu2023}]\label{thm:GLexistence}
Under condition \textbf{(C1)}, Problem~\ref{prob:nonlinearriemann} has
exactly one solution.
\end{theorem}

The proof is the Newton embedding procedure, and it is instructive to
see which of the results of Subsections~\ref{subsec:helmformulae} and
\ref{subsec:nonlinear} are consumed by it.  One embeds the problem in the
family obtained by replacing $F^{1}_{2}$ by $tF^{1}_{2}$,
$t\in[0,1]$.  For $t=0$ the unique solution is $\Gmat\equiv\Zero$, by
the Liouville-type theorem, Proposition~\ref{prop:liouville}.  Given a
solution for some $t_{0}<1$, one solves for $t$ close to $t_{0}$ by
successive approximations: $G^{1,m+1}_{2,t}$ is defined as the solution
of the \emph{linear} problem with datum
$F^{1}_{2}\bigl(\cdot,[G^{1,m}_{2,t}]^{+},[G^{1,m}_{2,t}]^{-}\bigr)$,
which exists and is unique by Theorem~\ref{thm:GLapriori}; the
differences of consecutive iterates then satisfy
\[\bigl\|G^{1,m+1}_{2,t}-G^{1,m}_{2,t}\bigr\|_{\nu}
\le\omega\,\bigl\|G^{1,m}_{2,t}-G^{1,m-1}_{2,t}\bigr\|_{\nu},
\qquad
\omega:=\tilde{c}_{\nu,\Gamma}N<1,\]
by \textbf{(C1)}, so that the sequence converges in $\|\cdot\|_{\nu}$.
Crucially, the admissible step length depends on $\omega$ and on
$\omega_{0}:=\tilde{c}_{\nu,\Gamma}
\|F^{1}_{2}(\cdot,\Zero,\Zero)\|_{\nu,\Gamma}$ but not on $t_{0}$;
hence finitely many steps reach $t=1$.  Uniqueness follows from
Theorem~\ref{thm:GLapriori} and \textbf{(C1)} again.

The same estimates measure how fast the approximation converges, which
is what makes the procedure of practical interest.

\begin{theorem}[Error estimate, \cite{GuLiu2023}]\label{thm:GLerror}
With $\omega$ and $\omega_{0}$ as above and
$c=\omega_{0}\dfrac{2-\omega}{(1-\omega)^{2}}$, the solution $\Gmat$ of
Problem~\ref{prob:nonlinearriemann} and its approximations
$G^{1,m}_{2,t}$ satisfy
\[
\bigl\|G^{1,m}_{2,t}-\Gmat\bigr\|_{\nu}
\le c\bigl[(t\omega)^{m+1}+(1-t)\bigr].
\]
\end{theorem}

\begin{remark}\label{rem:GLcompare}
Three comments.  (i) For $F^{1}_{2}$ independent of the last two
arguments one recovers the linear jump problem, and
Theorem~\ref{thm:GLapriori} is then nothing but the explicit solution
by the perturbed Cauchy transform.  (ii) Condition
\eqref{eq:GLsmallness} plays exactly the role played by
\eqref{eq:smallness} in Theorem~\ref{thm:riemannvar} and by
\eqref{eq:Gusmallness} in Theorem~\ref{thm:GuRiemann}: in all three
cases solvability is obtained by a contraction argument, and no index
theory is available which would dispense with the smallness assumption.
(iii) The whole argument uses the geometry of $\Gamma$ only through the
constant $c_{\nu,\Gamma}$ of Lemma~\ref{lem:GLcauchy}; this is what
makes it plausible that the scheme survives on the $d$-summable
boundaries of the next subsection, where, however, it has not been
carried out.
\end{remark}

\subsection{Perturbed Hermitian Dirac equations on fractal domains}

The paper \cite{A61} carries the Helmholtz theory into the fractal
regime, combining the perturbed operators of \cite{A45} with the
Whitney extension mechanism of Subsection~\ref{subsec:geometry}.  Write
$\mathcal{D}^{A}_{(Z,\Zd)}$ and $\mathcal{T}^{A}_{\Omega}$ for the
perturbed matrix Dirac and Teodorescu operators with the notation $A=-\tfrac14(\alpha-i\alpha)$, $\alpha\in\R$.

\begin{theorem}[Mapping properties, \cite{A61}]\label{thm:A61maps}
Let $\Omega\subset\R^{2n}$ and $p>2n$.  Then
$\mathcal{T}^{A}_{\Omega}:L^{p}(\Omega)\to W^{1}_{p}(\Omega)$ is
continuous, and for $\Gmat\in L^{p}(\Omega)$,
\[
\mathcal{D}^{A}_{(Z,\Zd)}\,\mathcal{T}^{A}_{\Omega}\bigl[\Gmat\bigr](Y)
=\begin{cases}\lam\,\Gmat(Y),&Y\in\Omega,\\
\Zero,&Y\in\R^{2n}\setminus\overline{\Omega}.\end{cases}
\]
\end{theorem}

\begin{proposition}[\cite{A61}]\label{prop:A61integrability}
Let $\Omega$ have $d$-summable boundary $\Gamma$, $d\in(2n-1,2n)$, and
let $\Gmat\in C^{0,\nu}(\Gamma)$.  Then
$\mathcal{D}^{A}_{(Z,\Zd)}\varUpsilon[{G}^{1}_{2}]\in L^{p}(\Omega)$ for every
$0<p<\frac{2n-d}{1-\nu}$.  If moreover $\nu>\frac{d}{2n}$, then
$\mathcal{T}^{A}_{\Omega}\mathcal{D}^{A}_{(Z,\Zd)}\varUpsilon[{G}^{1}_{2}]
\in C^{0,\mu}(\Gamma)$ whenever
$\mu<\frac{2n\nu-d}{2n-d}$.
\end{proposition}

The two problems treated in \cite{A61} are the jump problem
\begin{equation}\label{eq:A61jump}
\bigl[\Xmat\bigr]^{+}-\bigl[\Xmat\bigr]^{-}=\Gmat\ \text{ on }\Gamma,
\qquad \Xmat(\infty)=\Zero,
\end{equation}
and the Dirichlet-type problem
\begin{equation}\label{eq:A61dirichlet}
\mathcal{D}^{A}_{(Z,\Zd)}\Xmat=\Fmat\ \text{ in }\Omega,
\qquad
\Xmat=\Gmat\ \text{ on }\Gamma .
\end{equation}

\begin{theorem}[\cite{A61}]\label{thm:A61jump}
Let $\Gmat\in C^{0,\nu}(\Gamma)$ with $\nu>\frac{d}{2n}$.  Then the
jump problem \eqref{eq:A61jump} is solvable, and a solution is
explicitly given by
\[
\Xmat=\varUpsilon[{G}^{1}_{2}]
-\frac{1}{\lam}\,
\mathcal{T}^{A}_{\Omega}\mathcal{D}^{A}_{(Z,\Zd)}\varUpsilon[{G}^{1}_{2}]
\ \ \text{in }\Omega^{+},
\qquad
\Xmat=-\frac{1}{\lam}\,
\mathcal{T}^{A}_{\Omega}\mathcal{D}^{A}_{(Z,\Zd)}\varUpsilon[{G}^{1}_{2}]
\ \ \text{in }\Omega^{-}.
\]
\end{theorem}

\begin{theorem}[Homogeneous Dirichlet problem, \cite{A61}]
\label{thm:A61dirichlet}
Let $\Gmat\in C^{0,\nu}(\Gamma)$ with $\nu>\frac{d}{2n}$.  If there
exists $\Xmat\in C^{0,\nu}(\overline{\Omega})$ with
$\mathcal{D}^{A}_{(Z,\Zd)}\Xmat=\Zero$ in $\Omega$ and $\Xmat=\Gmat$
on $\Gamma$, then
\begin{equation}\label{eq:A61cond}
\Bigl[\mathcal{T}^{A}_{\Omega}
\mathcal{D}^{A}_{(Z,\Zd)}\varUpsilon[{G}^{1}_{2}]\Bigr]\Big|_{\Gamma}=\Zero .
\end{equation}
Conversely, if \eqref{eq:A61cond} holds, then there exists
$\Xmat\in C^{0,\mu}(\overline{\Omega})$, $\mu<\nu$, solving that
problem.
\end{theorem}

\begin{corollary}[Reproducing formula, \cite{A61}]\label{cor:A61repro}
Let $\Gmat\in C^{0,\nu}(\overline{\Omega})$ with $\nu>\frac{d}{2n}$ be
Helmholtz $H$-monogenic in $\Omega$ with trace $g^{1}_{2}$ on $\Gamma$. Then
\[
\Gmat(X)=\varUpsilon[{G}^{1}_{2}](X)
-\frac{1}{\lam}\,\mathcal{T}^{A}_{\Omega}\mathcal{D}^{A}_{(Z,\Zd)}
\varUpsilon[{g}^{1}_{2}](X),\qquad X\in\Omega .
\]
\end{corollary}

\begin{theorem}[Inhomogeneous problem, \cite{A61}]
\label{thm:A61inhom}
Let $\Gmat\in C^{0,\nu}(\Gamma)$ and $\Fmat\in L^{p}(\Omega)$ with
$\nu>\frac{d}{2n}$ and $\frac{p-2n}{p}>\frac{d}{2n}$.  If
\eqref{eq:A61dirichlet} admits a solution
$\Xmat\in C^{0,\nu}(\overline{\Omega})$, then the traces of
$\mathcal{T}^{A}_{\Omega}\mathcal{D}^{A}_{(Z,\Zd)}\varUpsilon[{g}^{1}_{2}]$ and
of $\mathcal{T}^{A}_{\Omega}\Fmat$ on $\Gamma$ coincide; conversely,
this coincidence yields a solution in a slightly smaller H\"older
class.
\end{theorem}

\subsection{The \texorpdfstring{$\bar\partial$}{dbar}-problem and
Maxwell's equations}\label{subsec:dbar}

The most recent contribution of the series, \cite{A152}, addresses the
inhomogeneous, or $\bar\partial$-type, problem associated with the perturbed
Hermitian Helmholtz operator. Write $A,A^{\dagger}$ for the constants
\eqref{eq:Kconst} replacing $k$ by a real parameter $\alpha$. 

\begin{definition}[\cite{A152}]\label{def:conjugate}
If $\Gmat$ satisfies $\mathcal{D}^{A}_{(Z,\Zd)}[\Gmat]=\Zero$, the pair
$(g_{1},g_{2})$ is called a pair of \emph{Helmholtz harmonic conjugate}
functions.
\end{definition}

\begin{theorem}[Solvability of the $\bar\partial$-problem, \cite{A152}]
\label{thm:dbar}
Let $\Fmat=\begin{pmatrix}f_{1}&f_{2}\\ f_{2}&f_{1}\end{pmatrix} \in C^{1}(\Omega)\cap C(\overline{\Omega})$.  The inhomogeneous system
\begin{equation}\label{eq:dbar-system}
\begin{cases}
(\partial_{Z^\dagger}-A^\dagger)g=f_1,\\
(\partial_Z+A)g=f_2,
\end{cases}
\iff
\begin{cases}
(\partial_X+\alpha)g=2(f_1-f_2),\\
(\partial_{X^{|}}+\alpha)g=-2i(f_1+f_2),
\end{cases}
\end{equation}
has a solution if and only if the function
\begin{equation}\label{eq:Pfunction}
\mathcal{P}:=-\frac{1}{\lam}\int_{\Omega}
\Bigl(\mathcal{E}_{A}f_{1}+\mathcal{E}^{\dagger}_{A}f_{2}\Bigr)\,dW(Z,\Zd)
\end{equation}
either
\begin{enumerate}[label=\upshape(\Alph*)]
\item is identically zero, or
\item possesses a Helmholtz harmonic conjugate function.
\end{enumerate}
In case \textup{(A)} the general solution of \eqref{eq:dbar-system} is
\[
g=-\frac{1}{\lam}\int_{\Omega}
\Bigl(\mathcal{E}_{A}f_{2}+\mathcal{E}^{\dagger}_{A}f_{1}\Bigr)\,dW(Z,\Zd)
+\hat{h},
\]
and in case \textup{(B)}
\[
g=-\frac{1}{\lam}\int_{\Omega}
\Bigl(\mathcal{E}_{A}f_{2}+\mathcal{E}^{\dagger}_{A}f_{1}\Bigr)\,dW(Z,\Zd)
+\hat{h}+\tilde{h},
\]
where $\tilde{h}$ is a Helmholtz harmonic conjugate of $-\mathcal{P}$ and $\hat{h}$ is an arbitrary Helmholtz $h$-monogenic function.
\end{theorem}

\begin{remark}
The Euclidean form of \eqref{eq:dbar-system} corrects \cite[System (10)]{A152}, where it is printed as $-(\partial_{X^{|}}+\alpha)g=f_1$, $(\partial_X+\alpha)g=f_2$. With the conventions of Subsection~\ref{subsec:hermitian}
and \eqref{eq:helmholtz}, each of $\partial_{Z^\dagger}-A^\dagger$ and $\partial_Z+A$ involves both $\partial_X+\alpha$ and $\partial_{X^{|}}+\alpha$, and the
printed system coincides with the Hermitian one only in the homogeneous case. The solvability criterion of the theorem refers to
the Hermitian form and is not affected.
\end{remark}
\begin{remark}
In \cite[Theorem~3]{A152} $\hat h$ is described as an arbitrary $h$-monogenic function (and as an arbitrary monogenic function in its proof). Since $\hat h$ arises as the second entry of a Helmholtz $H$-monogenic matrix whose first entry vanishes, it solves the homogeneous system associated with (\ref{eq:dbar-system}), i.e.\ it is Helmholtz $h$-monogenic; for $\alpha\neq0$ an $h$-monogenic function does not.
\end{remark}
The formal interpretation below, to two-dimensional electromagnetism given $\alpha\in\R$ is given as follows. Let one takes $n=1$, so that $\Omega\subset\mathbb R^2$, $\lambda_1=-2i$ and, by \eqref{eq:dVdW}, $dW(Z,Z^\dagger)=-2i\,dV$. For a vector field $F=F_1e_1+F_2e_2$ we write $\nabla\cdot F=\partial_{x_1}F_1+\partial_{x_2}F_2$ and $\nabla\times F=e_{12}\,(\partial_{x_1}F_2-\partial_{x_2}F_1)$, the latter being the bivector $\nabla\wedge F$. A direct computation gives
\begin{equation}\label{eq:divcurl}
(\partial_X+\alpha)F=-\nabla\cdot F+\nabla\times F+\alpha F,\qquad
(\partial_{X|}+\alpha)F=e_{12}\,\nabla\cdot F-e_{12}\,\nabla\times F+\alpha F. \footnote{Note that this system coincides with the homogeneous case of (\ref{eq:dbar-system}).}
\end{equation}
A magnetic field $\mathbf H$ and an electric field $\mathbf E$, both real vector fields in $\Omega$, are encoded into the single Clifford function
\begin{equation}\label{eq:Mfield}
M:=\mathbf H-i\,\mathbf E. 
\end{equation}

\begin{proposition}[\cite{A152}]\label{prop:left}
The function $M$ of (\ref{eq:Mfield}) is (left) $h$-monogenic in the Helmholtz sense if and only if $\mathbf H$ and $\mathbf E$ satisfy
\begin{equation}\label{eq:maxwellsys}
-i\nabla\times\mathbf E+\alpha\mathbf H=0,\qquad \nabla\times\mathbf H-i\alpha\mathbf E=0,\qquad \nabla\cdot\mathbf E=0,\qquad \nabla\cdot\mathbf H=0. 
\end{equation}
\end{proposition}
\begin{theorem}[\cite{A152}]\label{thm:maxwell}
Let $\alpha \in \R$, and let $f_1,f_2\in C^1(\Omega)\cap C(\overline\Omega)$ be $\mathbb C_2$-valued. If the system
\[
\begin{cases}
-e_{12}\bigl(\nabla\cdot\mathbf H-i\nabla\cdot\mathbf E\bigr)-\alpha\mathbf H+i\alpha\mathbf E+e_{12}\bigl(\nabla\times\mathbf H-i\nabla\times\mathbf E\bigr)=f_1,\\[0.2em]
-\bigl(\nabla\cdot\mathbf H-i\nabla\cdot\mathbf E\bigr)+\alpha\mathbf H-i\alpha\mathbf E+\nabla\times\mathbf H-i\nabla\times\mathbf E=f_2,
\end{cases}
\]
has a solution, then the function $\mathcal{P}$ of (\ref{eq:Pfunction}), here $P=\frac1{2i}\int_\Omega(\mathbf E_Af_1+\mathbf E_A^*{\dagger}f_2)\,dW(Z,Z^\dagger)$, either vanishes identically or admits a Helmholtz harmonic conjugate. Conversely, under either condition the formulae of Theorem~\ref{thm:dbar} with $n=1$ produce a $\mathbb C_2$-valued solution $M$ of (\ref{eq:dbar-system}); it corresponds to a pair of fields $(\mathbf H,\mathbf E)$ precisely when it is vector-valued.
\end{theorem}

\begin{remark}
(i) The solution of (\ref{eq:maxwellsys}), when it exists, is unique (the function $\hat h$ of Theorem~\ref{thm:dbar} vanishes), and in Proposition~\ref{prop:left} both conditions reduce to $M\equiv0$; indeed, with the above conventions, the vector parts of (\ref{eq:maxwellsys}) read $\alpha\mathbf H=0$ and $\alpha\mathbf E=0$. The Maxwell interpretation is therefore formal: system (49) mirrors the three-dimensional time-harmonic Maxwell equations, which lie outside the even-dimensional Hermitian setting.

(ii) For the same reason, conditions (A)--(B) do not by themselves guarantee a vector-valued solution: if $\varphi$ is a non-vanishing real scalar function and $f_1=-(\partial_{X|}+\alpha)\varphi$, $f_2=(\partial_X+\alpha)\varphi$, then (\ref{eq:dbar-system}) is solved by $\varphi$ and by no vector field.

(iii) For $\tau=\sum_A\tau_Ae_A\in\mathbb C_2$ put $\operatorname{Re}\tau=\sum_A(\operatorname{Re}\tau_A)e_A$ and $\operatorname{Im}\tau=\sum_A(\operatorname{Im}\tau_A)e_A$. If $M=\mathbf H-i\mathbf E$ solves the system of Proposition~\ref{prop:left}, then $\mathbf H=\operatorname{Re}M$, $\mathbf E=-\operatorname{Im}M$, and Theorem~\ref{thm:dbar} with $n=1$ gives
\[
\mathbf H=-\operatorname{Re}\int_\Omega\bigl(\mathcal{E}_Af_2+\mathcal{E}_A^{\dagger}f_1\bigr)\,dV+\operatorname{Re}\tilde h,\qquad
\mathbf E=\operatorname{Im}\int_\Omega\bigl(\mathcal{E}_Af_2+\mathcal{E}_A^{\dagger}f_1\bigr)\,dV-\operatorname{Im}\tilde h,
\]
with $\tilde h=0$ in case (A). The kernels $\mathcal{E}_A,\mathcal{E}_A^{\dagger}$ are complex-valued (the leading singularity of $E_\alpha$ is the real kernel $E$), so $\operatorname{Re}$ and $\operatorname{Im}$ cannot be distributed over the integrand term by term. The explicit formulae of \cite[Remark~1]{A152}, which rest on the assertion that $E_\alpha$ and $E^{|}_\alpha$ are purely imaginary, are not reproduced here.
\end{remark}

% ==================================================================
\section{The quaternionic Hermitian setting}\label{sec:quaternionic}
% ==================================================================

\subsection{Four Dirac operators and circulant \texorpdfstring{$(4\times
4)$}{(4x4)} matrices}

Replacing the complex structure by a quaternionic one leads to
\emph{quaternionic Hermitian Clifford analysis}.  Here $m=4n$ and one
works in $\Hqn=\Hq\otimes_{\R}\R_{0,4n}$.

\begin{definition}[Quaternionic Witt basis, \cite{PSS2007}]
\label{def:qwitt}
The quaternionic Witt basis of $\Hqn$ is the family
$\{f_{\ell},f^{\alpha}_{\ell},f^{\beta}_{\ell},f^{\gamma}_{\ell}\}$,
$\ell=1,\dots,n$, where
\[
\begin{aligned}
f_{\ell}&=e_{1+4(\ell-1)}-ie_{2+4(\ell-1)}-je_{3+4(\ell-1)}
-ke_{4+4(\ell-1)},\\
f^{\alpha}_{\ell}&=e_{1+4(\ell-1)}-ie_{2+4(\ell-1)}+je_{3+4(\ell-1)}
+ke_{4+4(\ell-1)},\\
f^{\beta}_{\ell}&=e_{1+4(\ell-1)}+ie_{2+4(\ell-1)}-je_{3+4(\ell-1)}
+ke_{4+4(\ell-1)},\\
f^{\gamma}_{\ell}&=e_{1+4(\ell-1)}+ie_{2+4(\ell-1)}+je_{3+4(\ell-1)}
-ke_{4+4(\ell-1)} .
\end{aligned}
\]
\end{definition}

To $(x_{1},\dots,x_{4n})\in\R^{4n}$ one associates four real Clifford
vectors $X_{0},X_{1},X_{2},X_{3}$, obtained from one another by
multiplication by the quaternionic units: $X_{0}=X$ is the usual
Clifford vector
\[
\begin{aligned}
X_{0}&=\sum_{\ell=1}^{n}\bigl(e_{4\ell-3}x_{4\ell-3}
+e_{4\ell-2}x_{4\ell-2}+e_{4\ell-1}x_{4\ell-1}+e_{4\ell}x_{4\ell}\bigr),\\
X_{1}&=\sum_{\ell=1}^{n}\bigl(e_{4\ell-3}x_{4\ell-2}
-e_{4\ell-2}x_{4\ell-3}-e_{4\ell-1}x_{4\ell}+e_{4\ell}x_{4\ell-1}\bigr),
\end{aligned}
\]
and $X_{2},X_{3}$ are given by the analogous expressions in which the
roles of the indices are permuted by $j$ and $k$ respectively, for
which we refer to \cite{PSS2007,A51}.  The vectors $X_{1},X_{2},X_{3}$
play the role of the twisted vector variable of the complex Hermitian
case; the corresponding Dirac operators are denoted
$\partial_{X_{0}},\dots,\partial_{X_{3}}$ and their Hermitian
counterparts $\partial_{Z_{0}},\dots,\partial_{Z_{3}}$.

\begin{definition}[\cite{PSS2007,A51,A55}]\label{def:qmono}
A continuously differentiable function $f:\Omega\to\Hqn$, $\Omega$ open
in $\R^{4n}$, is \emph{(left) $q$-Hermitian monogenic} (or
$q$-monogenic) if
\[
\partial_{Z_{0}}f=\partial_{Z_{1}}f=\partial_{Z_{2}}f
=\partial_{Z_{3}}f=0,
\quad\text{equivalently}\quad
\partial_{X_{0}}f=\partial_{X_{1}}f=\partial_{X_{2}}f
=\partial_{X_{3}}f=0 .
\]
Right and two-sided $q$-monogenicity are defined analogously.  Every
$q$-monogenic function is monogenic, hence harmonic.
\end{definition}

The Euclidean Cauchy kernels are
$E_{r}(X)=-\frac{1}{\sigma_{4n}}\frac{X_{r}}{|X|^{4n}}$, $r=0,1,2,3$,
$\sigma_{4n}$ being the area of $S^{4n-1}$, so that
$\partial_{X_{r}}E_{r}=\delta$; the quaternionic Hermitian kernels
$\mathcal{E}_{r}$ are obtained from them by the linear combinations
dictated by the quaternionic units.  As before, the $\mathcal{E}_{r}$
are not fundamental solutions of the $\partial_{Z_{r}}$, but the
matrices assembled from them are.  For a $4$-tuple
$(a_{0},a_{1},a_{2},a_{3})$ we write $\circu(a_{0},a_{1},a_{2},a_{3})$
for the associated circulant $(4\times 4)$ matrix.

\begin{theorem}[\cite{A51}]\label{thm:qfundsol}
Put
\[
\mathbf{D}=\circu(\partial_{Z_{0}},\partial_{Z_{1}},
\partial_{Z_{2}},\partial_{Z_{3}}),\qquad
\mathbf{E}=\circu(\mathcal{E}_{0},\mathcal{E}_{1},\mathcal{E}_{2},
\mathcal{E}_{3}),\qquad
\boldsymbol{\delta}=\mathrm{diag}(\delta,\delta,\delta,\delta).
\]
Then
\[
\mathbf{D}^{T}\mathbf{E}=\mathbf{E}\,\mathbf{D}^{T}
=\boldsymbol{\delta}.
\]
\end{theorem}

\begin{definition}[\cite{A51,A55}]\label{def:Qmono}
The circulant $(4\times 4)$ matrix function
$\mathbf{G}=\circu(g_{0},g_{1},g_{2},g_{3})$ is \emph{(left)
$Q$-Hermitian monogenic} ($Q$-monogenic) in $\Omega$ if
$\mathbf{D}^{T}\mathbf{G}=\Zero$ in $\Omega$, i.e.
\[
\begin{aligned}
\partial_{Z_{0}}g_{0}+\partial_{Z_{1}}g_{1}
+\partial_{Z_{2}}g_{2}+\partial_{Z_{3}}g_{3}&=0,\\
\partial_{Z_{3}}g_{0}+\partial_{Z_{0}}g_{1}
+\partial_{Z_{1}}g_{2}+\partial_{Z_{2}}g_{3}&=0,\\
\partial_{Z_{2}}g_{0}+\partial_{Z_{3}}g_{1}
+\partial_{Z_{0}}g_{2}+\partial_{Z_{1}}g_{3}&=0,\\
\partial_{Z_{1}}g_{0}+\partial_{Z_{2}}g_{1}
+\partial_{Z_{3}}g_{2}+\partial_{Z_{0}}g_{3}&=0 .
\end{aligned}
\]
Right $Q$-monogenicity is defined by $\mathbf{G}\mathbf{D}^{T}=\Zero$.
The diagonal matrix $\mathbf{G}_{0}=\mathrm{diag}(g,g,g,g)$ is
$Q$-monogenic if and only if $g$ is $q$-monogenic.
\end{definition}

\subsection{Cauchy integral formulae and the Martinelli--Bochner formula}

\begin{theorem}[Matricial Stokes formula, \cite{A51}]\label{thm:qstokes}
Let $\mathbf{F},\mathbf{G}$ be circulant matrix functions with entries
in $C^{1}(\Omega,\Hqn)$, and let $\Sigma\subset\Omega$ be a
$4n$-dimensional compact differentiable oriented manifold with
$C^{\infty}$ boundary $\partial\Sigma$.  Then
\[
\int_{\Sigma}\Bigl[\bigl(\mathbf{F}\mathbf{D}^{T}\bigr)\mathbf{G}
+\mathbf{F}\bigl(\mathbf{D}^{T}\mathbf{G}\bigr)\Bigr]dV
=\int_{\partial\Sigma}\mathbf{F}\,\mathbf{N}^{T}\mathbf{G}\,d\mathcal{H}^{4n-1}.
\]
\end{theorem}

From now on $\Sigma^{+}$ will stand for the interior of $\Sigma$ and $\Sigma^{-}$ for $\Omega\setminus \Sigma$.

\begin{theorem}[$Q$-Hermitian Borel--Pompeiu formula, \cite{A51}]
\label{thm:qBP}
Under the hypotheses of Theorem~\ref{thm:qstokes},
\[
\int_{\partial\Sigma}\mathbf{E}(Z-V)\mathbf{N}^{T}(Z)\mathbf{G}(X)d\mathcal{H}^{4n-1}(X)
-\int_{\Sigma}\mathbf{E}(Z-V)
\bigl(\mathbf{D}^{T}\mathbf{G}(X)\bigr)dV(X)
=\begin{cases}\mathbf{G}(Y),&Y\in\Sigma^{+},\\
\Zero,&Y\in\Sigma^{-}.\end{cases}
\]
\end{theorem}

\begin{theorem}[$Q$- and $q$-Hermitian Cauchy integral formulae,
\cite{A51,A56}]\label{thm:qcauchy}
If $\mathbf{G}$ is $Q$-monogenic in $\Omega$, then
\[
\int_{\partial\Sigma}\mathbf{E}(Z-V)\mathbf{N}^{T}(Z)\mathbf{G}(X)d\mathcal{H}^{4n-1}(X)
=\begin{cases}\mathbf{G}(Y),&Y\in\Sigma^{+},\\
\Zero,&Y\in\Sigma^{-},\end{cases}
\]
and in particular, if $g$ is $q$-monogenic, the same formula holds with
$\mathbf{G}$ replaced by the diagonal matrix $\mathbf{G}_{0}$.  The
matrix $\mathbf{E}$ is therefore called the \emph{quaternionic
Hermitian Cauchy kernel}.
\end{theorem}

\begin{remark}[Martinelli--Bochner formula, \cite{A51}]\label{rem:qMB}
Taking $\mathbf{G}=\mathbf{G}_{0}$ in Theorem~\ref{thm:qBP} and comparing entries produces four scalar identities, the first of which reads 
\[
\sum_{s=0}^{3}\left[\int_{\partial\Sigma}
\mathcal{E}_{s}(Z-V)N_{s}(Z)g(X)d\mathcal{H}^{4n-1}(X)
-\int_{\Sigma}\mathcal{E}_{s}(Z-V)
\bigl(\partial_{Z_{s}}g(X)\bigr)dV(X)\right]=g(Y).
\] 
Applying the Stokes formula of Theorem~\ref{thm:qstokes} to the terms with $s=1,2,3$ on $\Sigma\setminus B(Y,R)$, using $\sum_{s=0}^3E_s(Z-V)\partial_{Z_s}=0$ for $Z\neq V$ and letting $R\to0$, it can be rearranged into a formula involving only $E_0$, $N_0$ and $\partial_{Z_0}$ together with a principal value integral; upon the translations $E_0N_0=\frac14E_qn_q$ and $E_0\partial_{Z_0}=\frac14E_q\partial_q$, where $E_q$, $n_q$ and $\partial_q$ are the kernel, the quaternionic normal and the quaternionic Hermitian vector derivative defined in \cite{PSS2007}, this is the Martinelli--Bochner-type formula of \cite[Theorem~4.3]{PSS2007} and \cite[Section~5]{A51}. The three
remaining identities are genuinely new and mix the kernels
$\mathcal{E}_{r}$ with the normal components $N_{s}$, $r\neq s$.
\end{remark}

\subsection{Cauchy transform, Hilbert transform and projections}

From now on $\Gamma$ is a compact Liapunov surface bounding a domain
$\Omega=\Omega^{+}\subset\R^{4n}$.  The \emph{$Q$-Hermitian Cauchy
transform} is
\begin{equation}\label{eq:qcauchytransform}
\mathcal{C}[\mathbf{G}](Y)=\int_{\Gamma}
\mathbf{E}(Z-V)\mathbf{N}^{T}(Z)\mathbf{G}(X)\,d\mathcal{H}^{4n-1}(X),
\qquad Y\notin\Gamma,
\end{equation}
which is $Q$-monogenic in $\R^{4n}\setminus\Gamma$ and vanishes at
infinity.  It decomposes into the Euclidean Cauchy-type integrals
\begin{equation}\label{eq:Crs}
C_{r,s}g(Y)=\int_{\Gamma}E_{r}(X-Y)n_{s}(X)g(X)\,d\mathcal{H}^{4n-1}(X),
\qquad r,s=0,1,2,3,
\end{equation}
by means of
\begin{equation}\label{eq:Cdecomp}
\mathcal{C}[\mathbf{G}]=\tfrac14\,
\circu\bigl(\mathcal{A}_{0}[C],\mathcal{A}_{1}[C],
\mathcal{A}_{2}[C],\mathcal{A}_{3}[C]\bigr)[\mathbf{G}],
\end{equation}
where, for a family $T=\{T_{r,s}\}$ of operators,
\[
\begin{aligned}
\mathcal{A}_{0}[T]&=T_{0,0}+T_{1,1}+T_{2,2}+T_{3,3},
&\quad
\mathcal{A}_{1}[T]&=T_{0,0}-T_{2,2}+j\bigl(T_{1,3}+T_{3,1}\bigr),\\
\mathcal{A}_{2}[T]&=T_{0,0}-T_{1,1}+T_{2,2}-T_{3,3},
&\quad
\mathcal{A}_{3}[T]&=T_{0,0}-T_{2,2}-j\bigl(T_{1,3}+T_{3,1}\bigr).
\end{aligned}
\]

\begin{lemma}[\cite{A56}]\label{lem:qCauchyprops}
In $\R^{4n}\setminus\Gamma$ one has, for $r,s,t=0,1,2,3$,
\[
\partial_{Y_{r}}C_{r,s}g=0,
\qquad
\partial_{Y_{r}}C_{s,t}g+\partial_{Y_{s}}C_{r,t}g=0 .
\]
\end{lemma}

\begin{proposition}[Diagonality criterion, \cite{A56}]
\label{prop:diagonality}
$\mathcal{C}[\mathbf{G}_{0}]$ is a diagonal matrix, or, equivalently,
its entries are $q$-monogenic, if and only if
\begin{equation}\label{eq:diagcond}
C_{0,0}g=C_{2,2}g,\qquad
C_{1,3}g=-C_{3,1}g,\qquad
2C_{0,0}g=C_{1,1}g+C_{3,3}g,
\end{equation}
in which case
\[
\mathcal{C}[\mathbf{G}_{0}]
=\mathrm{diag}\bigl(C_{0,0}g\bigr)
=\mathrm{diag}\bigl(C_{2,2}g\bigr)
=\tfrac12\,\mathrm{diag}\bigl(C_{1,1}g+C_{3,3}g\bigr).
\]
\end{proposition}

Introducing the singular integrals
\begin{equation}\label{eq:Hrs}
H_{r,s}g(U)=2\int_{\Gamma}E_{r}(X-U)n_{s}(X)g(X)\,d\mathcal{H}^{4n-1}(X),
\qquad U\in\Gamma,
\end{equation}
understood in the principal value sense, one defines the \emph{matrix
$Q$-Hermitian Hilbert transform} by the very same combination,
\begin{equation}\label{eq:qhilbert}
\mathcal{H}[\mathbf{G}]=\tfrac14\,
\circu\bigl(\mathcal{A}_{0}[H],\mathcal{A}_{1}[H],
\mathcal{A}_{2}[H],\mathcal{A}_{3}[H]\bigr)[\mathbf{G}].
\end{equation}

\begin{theorem}[Plemelj--Sokhotski formulae, \cite{A56,A55}]
\label{thm:qPS}
Let $\mathbf{G}\in C^{0,\nu}(\Gamma)$, $0<\nu\le 1$.  Then the
continuous limit values of $\mathcal{C}[\mathbf{G}]$ exist and are
given by
\[
\mathcal{C}^{\pm}[\mathbf{G}](U)
=\tfrac12\bigl(\mathcal{H}[\mathbf{G}](U)\pm\mathbf{G}(U)\bigr),
\qquad U\in\Gamma .
\]
The same formula holds for $\mathbf{G}_{0}$.
\end{theorem}

\begin{theorem}[\cite{A56,A55}]\label{thm:qhilbertprops}
The operator $\mathcal{H}$ satisfies:
\begin{enumerate}[label=\upshape(\roman*)]
\item $\mathcal{H}$ is a bounded linear operator on
$\bigl(C^{0,\nu}(\Gamma),\|\cdot\|_{\nu}\bigr)$, $0<\nu<1$;
\item $\mathcal{H}$ is an involution on $C^{0,\nu}(\Gamma)$, i.e.
$\mathcal{H}^{2}=\mathbf{I}$, the $(4\times 4)$ identity matrix
operator.
\end{enumerate}
Consequently $P^{+}_{\Gamma}=\tfrac12(\mathbf{1}+\mathcal{H})$ and
$P^{-}_{\Gamma}=\tfrac12(\mathbf{1}-\mathcal{H})$ satisfy
$P^{+}_{\Gamma}+P^{-}_{\Gamma}=\mathbf{1}$, $P^{+}_{\Gamma}-P^{-}_{\Gamma}=\mathcal{H}$, $({P^{+}_{\Gamma}})^{2}=P^{+}_{\Gamma}$, $({P^{-}_{\Gamma}})^{2}=P^{-}_{\Gamma}$,
$P^{+}_{\Gamma}P^{-}_{\Gamma}=P^{-}_{\Gamma}P^{+}_{\Gamma}=\Zero$, and give rise to the direct sum decomposition
\[C^{0,\alpha}(\Gamma)=P^{+}_{\Gamma}\bigl[C^{0,\nu}(\Gamma)\bigr]
\oplus P^{-}_{\Gamma}\bigl[C^{0,\nu}(\Gamma)\bigr].\]
\end{theorem}
Right-hand versions of all these operators are obtained from the
alternative definitions
\[[\mathbf{G}]\mathcal{C}(Y)=\int_{\Gamma}\mathbf{G}(X)
\mathbf{N}^{T}(Z)\mathbf{E}(Z-V)d\mathcal{H}^{4n-1}(X)\] 
and $g H_{r,s}(U)=2\int_{\Gamma}g(X)n_{s}(X)E_{r}(X-U)d\mathcal{H}^{4n-1}(X)$.

\subsection{Boundary value problems on Liapunov surfaces}
The following theorem solves the jump problem
\begin{equation}\label{eq:qjump}
\psi^{+}(U)-\psi^{-}(U)=g(U),\qquad U\in\Gamma,
\end{equation}
for $q$-monogenic functions.  Note how the solvability criteria are
exactly the diagonality conditions \eqref{eq:diagcond}, and how they
collapse into the strikingly symmetric form (iii).

\begin{theorem}[\cite{A55}]\label{thm:qjumpleft}
Let $g\in C^{0,\nu}(\Gamma)$.  The following are equivalent:
\begin{enumerate}[label=\upshape(\roman*)]
\item the jump problem \eqref{eq:qjump} is solvable in terms of left
$q$-monogenic functions;
\item $C_{0,0}g=C_{2,2}g$, $C_{1,3}g=-C_{3,1}g$,
$2C_{0,0}g=C_{1,1}g+C_{3,3}g$;
\item $C_{0,0}g=C_{1,1}g=C_{2,2}g=C_{3,3}g$.
\end{enumerate}
The statement for right $q$-monogenic functions is obtained by writing
all the operators on the right of $g$.
\end{theorem}

\begin{theorem}[Dirichlet problem, matrix version, \cite{A55}]
\label{thm:qdirichletmatrix}
Let $\mathbf{G}\in C^{0,\nu}(\Gamma)$.  The Dirichlet problem
\[
\mathbf{D}^{T}\mathbf{F}=\Zero\ \ \bigl(\text{resp.\ }
\mathbf{F}\mathbf{D}^{T}=\Zero\bigr)\ \text{in }\Omega,
\qquad
\mathbf{F}=\mathbf{G}\ \text{on }\Gamma,
\]
has a solution if and only if $\mathcal{H}[\mathbf{G}]=\mathbf{G}$
(resp.\ $[\mathbf{G}]\mathcal{H}=\mathbf{G}$).
\end{theorem}

\begin{theorem}[Dirichlet problem, scalar version, \cite{A55}]
\label{thm:qdirichletscalar}
Let $g\in C^{0,\nu}(\Gamma)$.  The following are equivalent:
\begin{enumerate}[label=\upshape(\roman*)]
\item the Dirichlet problem
$\partial_{Z_{0}}f=\partial_{Z_{1}}f=\partial_{Z_{2}}f
=\partial_{Z_{3}}f=0$ in $\Omega$, $f=g$ on $\Gamma$, has a solution;
\item $H_{0,0}g=H_{2,2}g=g$, $H_{1,3}g=-H_{3,1}g$,
$H_{1,1}g+H_{3,3}g=2g$;
\item $H_{0,0}g=H_{1,1}g=H_{2,2}g=H_{3,3}g=g$.
\end{enumerate}
The corresponding statement for right $q$-monogenic functions is
obtained by writing the operators on the right of $g$.
\end{theorem}

\begin{theorem}[Two-sidedness, \cite{A55}]\label{thm:qtwosided}
Let $\mathbf{G}\in C^{0,\nu}(\overline{\Omega})$ with
$\mathbf{D}^{T}\mathbf{G}=\Zero$ in $\Omega$.  Then $\mathbf{G}$ is
two-sided $Q$-monogenic in $\Omega$ if and only if
$\mathcal{H}[\mathbf{G}|_{\Gamma}]
=[\mathbf{G}|_{\Gamma}]\mathcal{H}$.
Correspondingly, a left $q$-monogenic $g\in
C^{0,\nu}(\overline{\Omega})$ is two-sided $q$-monogenic in $\Omega$
if and only if
\[
H_{r,r}g=gH_{r,r},\qquad r=0,1,2,3 .
\]
\end{theorem}

\subsection{Boundary value problems on fractal hypersurfaces}

Finally, \cite{A65} carries the whole quaternionic programme to Jordan
domains $\Omega\subset\R^{4n}$ with $d$-summable boundary $\Gamma$,
$d\in(4n-1,4n)$, combining the operators of the previous subsection
with the Whitney technique of Subsection~\ref{subsec:geometry}, with
$2n$ replaced by $4n$ throughout.  One defines the fractal
$Q$-Hermitian Cauchy integral $\mathcal{C}^{*}$ exactly as in
Definition~\ref{def:fractalcauchy}, and the associated Hilbert
transform $\mathcal{H}^{*}$, which admits the representation
\begin{equation}\label{eq:qfractalhilbert}
\mathcal{H}^{*}[\mathbf{G}](U)=\mathbf{G}(U)
-2\,\mathcal{T}_{\Omega}\mathbf{D}^{T}\varUpsilon[{\mathbf{G}}](U),
\qquad U\in\Gamma .
\end{equation}

\begin{theorem}[\cite{A65}]\label{thm:qfractalhilbertwell}
Let $\Omega\subset\R^{4n}$ be a Jordan domain with $d$-summable
boundary $\Gamma$, $d\in(4n-1,4n)$, and let
$\mathbf{G}\in C^{0,\nu}(\Gamma)$ with $d-4n+1<\nu\le 1$.  If
moreover
\begin{equation}\label{eq:qalphacond}
\nu>\frac{d}{4n},
\end{equation}
then $\mathcal{C}^{*}\mathbf{G}$ admits a continuous extension to
$\overline{\Omega}$ and
$\mathcal{H}^{*}\mathbf{G}\in C^{0,\mu}(\Gamma)$ for every
\begin{equation}\label{eq:qmucond}
\mu<\frac{4n\nu-d}{4n-d}.
\end{equation}
\end{theorem}

\begin{theorem}[Plemelj--Sokhotski formulae, fractal case, \cite{A65}]
\label{thm:qfractalPS}
Under the hypotheses of Theorem~\ref{thm:qfractalhilbertwell}, the
continuous limit values of $\mathcal{C}^{*}[\mathbf{G}]$ exist and
\[
\bigl(\mathcal{C}^{*}\bigr)^{\pm}[\mathbf{G}](U)
=\tfrac12\bigl(\mathcal{H}^{*}[\mathbf{G}](U)\pm\mathbf{G}(U)\bigr),
\qquad U\in\Gamma .
\]
\end{theorem}

\begin{theorem}[Cauchy formula, fractal case, \cite{A65}]
\label{thm:qfractalcauchy}
Under the same hypotheses, if $\mathbf{G}\in
C^{0,\nu}(\overline{\Omega})$ is $Q$-Hermitian monogenic in
$\Omega$, then
\[
\mathcal{C}^{*}\mathbf{G}(Y)=
\begin{cases}\mathbf{G}(Y),&Y\in\Omega,\\
\Zero,&Y\in\R^{4n}\setminus\overline{\Omega}.\end{cases}
\]
\end{theorem}

Writing $C^{*}_{r,s}$ and $H^{*}_{r,s}$ for the fractal analogues of
\eqref{eq:Crs} and \eqref{eq:Hrs}, one recovers all the solvability
criteria of the Liapunov case, under the additional dimensional
restriction
\begin{equation}\label{eq:qdimcond}
\dimH(\Gamma)-4n+1<\frac{4n\nu-d}{4n-d}.
\end{equation}

\begin{theorem}[Jump problem, fractal case, \cite{A65}]
\label{thm:qfractaljump}
Let $\Omega\subset\R^{4n}$ be a Jordan domain with $d$-summable
boundary $\Gamma$, assume \eqref{eq:qdimcond} and let
$g\in C^{0,\nu}(\Gamma)$.  The following are equivalent:
\begin{enumerate}[label=\upshape(\roman*)]
\item the jump problem \eqref{eq:qjump} for left $q$-Hermitian
monogenic functions has a solution of class $C^{0,\mu}$, $\mu$ as in
\eqref{eq:qmucond};
\item $C^{*}_{0,0}g=C^{*}_{2,2}g$, $C^{*}_{1,3}g=-C^{*}_{3,1}g$,
$2C^{*}_{0,0}g=C^{*}_{1,1}g+C^{*}_{3,3}g$;
\item $C^{*}_{0,0}g=C^{*}_{1,1}g=C^{*}_{2,2}g=C^{*}_{3,3}g$.
\end{enumerate}
The right-handed analogue holds with the operators acting on the right
of $g$.
\end{theorem}

\begin{theorem}[Dirichlet problem, fractal case, \cite{A65}]
\label{thm:qfractaldirichlet}
Let $\mathbf{G}\in C^{0,\nu}(\Gamma)$ with $\alpha$ as in
\eqref{eq:qalphacond}.  The Dirichlet problem
\[
\mathbf{D}^{T}\mathbf{F}=\Zero\ \bigl(\text{resp.\ }
\mathbf{F}\mathbf{D}^{T}=\Zero\bigr)\ \text{in }\Omega,
\qquad \mathbf{F}=\mathbf{G}\ \text{on }\Gamma,
\]
has a solution if and only if $\mathcal{H}^{*}[\mathbf{G}]=\mathbf{G}$
(resp.\ $[\mathbf{G}]\mathcal{H}^{*}=\mathbf{G}$).  Moreover, for
$g\in C^{0,\alpha}(\Gamma)$ and under \eqref{eq:qdimcond}, the scalar
Dirichlet problem
$\partial_{Z_{0}}f=\dots=\partial_{Z_{3}}f=0$ in $\Omega$, $f=g$ on
$\Gamma$, is solvable if and only if
\[
H^{*}_{0,0}g=H^{*}_{2,2}g=g,\quad
H^{*}_{1,3}g=-H^{*}_{3,1}g,\quad
H^{*}_{1,1}g+H^{*}_{3,3}g=2g,
\]
equivalently if and only if
$H^{*}_{0,0}g=H^{*}_{1,1}g=H^{*}_{2,2}g=H^{*}_{3,3}g=g$.
\end{theorem}

\begin{theorem}[Two-sidedness, fractal case, \cite{A65}]
\label{thm:qfractaltwosided}
Let $\alpha$ be as in \eqref{eq:qalphacond}.
\begin{enumerate}[label=\upshape(\roman*)]
\item If $\mathbf{G}\in C^{0,\nu}(\overline{\Omega})$ satisfies
$\mathbf{D}^{T}\mathbf{G}=\Zero$ in $\Omega$, then $\mathbf{G}$ is
two-sided $Q$-Hermitian monogenic in $\Omega$ if and only if
$\mathcal{H}^{*}[\mathbf{G}|_{\Gamma}]
=[\mathbf{G}|_{\Gamma}]\mathcal{H}^{*}$.
\item If $g\in C^{0,\nu}(\overline{\Omega})$ is $q$-Hermitian
monogenic in $\Omega$, then $g$ is two-sided $q$-Hermitian monogenic in
$\Omega$ if and only if
\[
\begin{aligned}
H^{*}_{0,0}g&=gH^{*}_{0,0}, &\qquad
H^{*}_{2,2}g&=gH^{*}_{2,2},\\[1mm]
H^{*}_{1,3}g+H^{*}_{3,1}g&=gH^{*}_{1,3}+gH^{*}_{3,1}, &\qquad
H^{*}_{1,1}g+H^{*}_{3,3}g&=gH^{*}_{1,1}+gH^{*}_{3,3},
\end{aligned}
\]
equivalently if and only if $H^{*}_{r,r}g=gH^{*}_{r,r}$ for
$r=0,1,2,3$.
\end{enumerate}
\end{theorem}

% ==================================================================
\section{Concluding remarks}\label{sec:remarks}
% ==================================================================

In light of the results detailed above, a few structural
features deserve to be singled out.

\medskip
\noindent
\textbf{The matrix formalism as an organizing principle.}
Every difficulty specific to the Hermitian setting---the absence of a
fundamental solution for the individual operators, the overdetermined
nature of the system, the failure of a naive Cauchy formula---is
absorbed by the passage to circulant matrices.  The price is that the
matrix theory is genuinely richer than the scalar one: $H$-monogenicity
of $\Gmat$ does not imply $h$-monogenicity of $g_{1},g_{2}$. Consequently, nearly every theorem in the survey comes in two versions,
a matricial one---clean, unconditional, structurally parallel to the
Euclidean theory---and a scalar one, whose hypotheses are exactly the
conditions ensuring that the relevant matrix object stays diagonal.
Theorem~\ref{thm:diagjump}, Theorem~\ref{thm:teodecompdiag} and
Proposition~\ref{prop:diagonality} are the three prototypes of this
phenomenon.

\medskip
\noindent
\textbf{The algebra behind the matrix device.}
The same phenomenon has an exact algebraic shadow, which
Section~\ref{sec:duality} makes visible.  The module associated with the
scalar Hermitian Dirac system has a free resolution which is linear of
length $n$, with $2(2^n-n-1)$ first syzygies, with Betti numbers 
$\beta_i=2\sum_{j\ge i}\binom{n}{j}$ \cite{DES2009}, and the cohomology vanishing that
would produce duality theorems by general principles simply fails
(Theorem~\ref{thm:resolution}); the circulant matrix operator $\Dop$,
by contrast, has \emph{trivial} first syzygies
(Proposition~\ref{prop:trivialsyz}), whence a Poincar\'e lemma, a
Mittag--Leffler theorem and the vanishing $H^{1}(U,\HM_{r})=\Zero$ follow.  What
appears analytically as ``the matrix theory has all the classical
integral formulae'' is thus, algebraically, the statement that the
overdetermination of the system has been absorbed into the matrix
structure.  The failure of a Hartogs-type extension phenomenon
(Corollary~\ref{cor:noHartogs}) belongs to the same list of scalar
pathologies.

\medskip
\noindent
\textbf{Fractality is a matter of definitions, not only of estimates.}
The passage from Ahlfors--David regular to $d$-summable boundaries
required redefining the Cauchy integral itself
(Definition~\ref{def:fractalcauchy}), after which the classical
formulae---Borel--Pompeiu, Koppelman, Cauchy, Plemelj--Sokhotski---
all survive, and the classical case is recovered as a corollary
(Lemma~\ref{lem:ADdsum}).  What genuinely changes is the range of
admissible H\"older exponents: the conditions $\nu>\frac{d}{2n}$ and
$\beta<\frac{2n\nu-d}{2n-d}$ of Theorem~\ref{thm:fractalhilbertwell}
and  reduce to those of Theorem~\ref{thm:normestimateAD} as $d\downarrow 2n-1$.

\medskip
\noindent
\textbf{Uniformity of the solvability criteria.}
Whatever the geometry, the differential operator or the algebraic
structure, the answer to every boundary value problem considered is
formulated in terms of the associated Hilbert transform or of the
Euclidean Cauchy-type integrals out of which it is built.  The
comparison of Theorems~\ref{thm:qjumpleft}, \ref{thm:qdirichletscalar},
\ref{thm:qfractaljump} and \ref{thm:qfractaldirichlet} is particularly
telling: the conditions $C_{0,0}g=C_{1,1}g=C_{2,2}g=C_{3,3}g$ and
$H_{r,r}g=g$ persist unchanged from Liapunov to fractal boundaries,
only the meaning of the operators being adjusted.

\medskip
\noindent
\textbf{Duality repeats a Euclidean pattern.}
The three-layered structure of Section~\ref{sec:duality} is not an
accident of the Hermitian setting.  Already in Euclidean Clifford
analysis, the dual of the module of germs of monogenic functions is
represented by a single monogenic function outside the compact set $\mathcal{K}$ on which the germs live (Theorem~\ref{thm:eucdual}), whereas the dual of the space of harmonic
multivector fields, that is, of harmonic differential forms, is
represented only by a \emph{pair} subject to a Moisil--Teodorescu-
type system, surjectively and with a kernel consisting of the pairs
both of whose components are themselves harmonic
(Theorems~\ref{thm:harmpairdual}--\ref{thm:formsdual}).  Passing to the
Hermitian refinement, circulant matrices take the place of Clifford-
algebra-valued-functions and $h$-monogenic functions that of harmonic
multivector fields, and the same dichotomy reappears verbatim in
Theorems~\ref{thm:dualmatrix}, \ref{thm:dualscalar},
\ref{thm:kernelJ} and \ref{thm:CSSscalar}.

\medskip
\noindent
\textbf{Some directions left open.}
Several natural questions are explicitly raised in the sources or
suggest themselves.  (i) The H\"older norm estimates of \cite{A49} have
no known analogue for the perturbed (Helmholtz) Hilbert transform on
fractal boundaries, nor in the quaternionic setting; obtaining them
would complete the picture of Sections~\ref{sec:helmholtz} and
\ref{sec:quaternionic}.  (ii) The duality theory of \cite{A47} is so
far confined to the unperturbed complex Hermitian case; a quaternionic
counterpart, presumably in terms of quadruples rather than pairs, is
missing.  (iii) The variable-coefficient and nonlinear Riemann problems of \ldots\ are 
solved under the smallness assumptions \eqref{eq:smallness},
\eqref{eq:Gusmallness} and \eqref{eq:GLsmallness} respectively, and
by one and the same contraction mechanism. 

A genuine index theory for the matrix operator $\Phi^{1}_{2}\mapsto P^{+}_{K}[\Phi^{1}_{2}]+G^{1}_{2}\,P^{-}_{K}[\Phi^{1}_{2}]$, which
would dispense with it, would be the natural next step.  Nothing is
known either about nonlinear transmission conditions on $d$-summable
boundaries, where the constant of Lemma~\ref{lem:GLcauchy} would have
to be replaced by the H\"older norm estimates of
Theorem~\ref{thm:normestimate}.  (iv) The $L^{p}$ theory is well developed on
Lipschitz hypersurfaces but essentially untouched in the fractal
regime.  (v) The $\bar\partial$-problem of \cite{A152} and its
electromagnetic interpretation invite the study of the corresponding
boundary value problems---jump, Dirichlet, Riemann---for two-dimensional Maxwell-type 
fields in the Hermitian Helmholtz framework, in each of the geometric
settings surveyed here.  (vi) The duality theory of
Section~\ref{sec:duality} is formulated for germs on compact sets and
knows nothing of the geometry of $\Gamma$. Combining it with the fractal machinery of Sections~\ref{sec:integral}--\ref{sec:hilbert}---for instance, describing the annihilators of $\HM_l(\overline\Omega)$ when $\Gamma$ is merely $d$-summable---appears to be entirely open, as does its quaternionic analogue. (vii) Finally, several nearby
developments would deserve to be brought into the present framework:
the Laurent expansions and $L^{2}$ boundary values recently obtained in
\cite{HH2024}, which are the natural point of departure for the $L^{p}$
theory mentioned in (iv); and the half Dirichlet problems of Ku, K\"ahler, Wang \cite{KKW2013}.

\medskip
\noindent
\textbf{A last word on scope.}
The line of research surveyed here is one strand of a much broader
fabric.  Hermitian Clifford analysis has been developed in parallel
along algebraic and representation-theoretic lines---Fischer
decompositions, Cauchy--Kovalevskaya extensions and the submonogenic
system, radial algebras, invariant syzygies---and it has been
transplanted to superspace and to the symplectic, discrete and split-quaternionic settings. 
For a systematic account of that side of the
theory, and for the many results we have had to leave aside, the reader
is referred to the survey of Sabadini and Sommen \cite{SabSom2025}.

% ==================================================================
\section*{Statements and Declarations}
% ==================================================================

% ==================================================================
\subsection*{Funding} The third author thanks Secretar\'ia de Investigaci\'on y Posgrado, Instituto Polit\'ecnico Nacional (grant number SIP20260491) for partial financial assistance.
% ==================================================================

% ==================================================================
\subsection*{Author contributions} All authors provided critical feedback, contributed to the drafting of the paper, and approved the final manuscript. 
% ==================================================================

% ==================================================================
\subsection*{Competing interests} 
% ==================================================================
The authors declare no conflicts of interest or personal ties that could have influenced the present work.
% ==================================================================

% ==================================================================
\subsection*{Data availability}
% ==================================================================
Data sharing is not applicable to this article, as no data sets were generated or analysed.

\subsection*{ORCID}
% ==================================================================
\noindent
Ricardo Abreu-Blaya: \url{https://orcid.org/0000-0003-1453-7223}\\
Juan Bory-Reyes: \url{https://orcid.org/0000-0002-7004-1794}\\
Carlos Daniel Tamayo-Castro: \url{https://orcid.org/0000-0002-5584-8274}

% ==================================================================
\subsection*{Generative AI in the writing process}
% ==================================================================
During the preparation of this work the authors used Claude Opus 5 (Anthropic) to improve the language, grammar and readability of sections they had drafted, to help synthesize sources, and to unify notation. After using this tool, the authors reviewed and edited the content as needed and take full responsibility for the content of the article.

% ==================================================================

\end{document}